\documentclass[11pt,a4paper]{amsart}
\usepackage[T1]{fontenc}
\usepackage{lmodern}
\usepackage[textwidth=158mm,textheight=243mm,centering]{geometry}
\usepackage{amsmath,amssymb,amsthm,mathtools,mathrsfs}
\usepackage{booktabs,microtype,needspace}
\usepackage[dvipsnames]{xcolor}
\usepackage[hidelinks]{hyperref}
\hypersetup{pdftitle={Multivariate diagonal harmonics: character polynomials and proofs of two conjectures of Bergeron},pdfauthor={Eren Ercan}}
\numberwithin{equation}{section}
\newtheorem{theorem}{Theorem}[section]
\newtheorem{proposition}[theorem]{Proposition}
\newtheorem{lemma}[theorem]{Lemma}
\newtheorem{corollary}[theorem]{Corollary}
\newtheorem{maintheorem}{Theorem}

\theoremstyle{definition}
\newtheorem{definition}[theorem]{Definition}
\theoremstyle{remark}
\newtheorem{remark}[theorem]{Remark}
\newtheorem{example}[theorem]{Example}
\DeclareMathOperator{\GL}{GL}
\DeclareMathOperator{\Sym}{Sym}
\DeclareMathOperator{\Tor}{Tor}
\DeclareMathOperator{\Hom}{Hom}
\DeclareMathOperator{\sgn}{sgn}
\DeclareMathOperator{\Frob}{Frob}
\DeclareMathOperator{\ch}{ch}
\DeclareMathOperator{\SYT}{SYT}
\DeclareMathOperator{\Des}{Des}
\DeclareMathOperator{\des}{des}
\DeclareMathOperator{\maj}{maj}
\DeclareMathOperator{\Alt}{Alt}
\DeclareMathOperator{\Av}{Av}

\DeclareMathOperator{\cont}{cont}
\newcommand{\id}{\operatorname{id}}
\newcommand{\Harm}{\mathcal H}
\newcommand{\Pol}{\mathcal P}
\newcommand{\Schur}{\mathbb S}
\newcommand{\SH}{\mathit{SH}}
\newcommand{\Q}{\mathbb Q}
\newcommand{\C}{\mathbb C}
\newcommand{\m}{\mathfrak m}
\DeclareMathOperator{\Span}{span}
\newcommand{\N}{\mathbb N}
\newcommand{\Z}{\mathbb Z}
\newcommand{\eps}{\varepsilon}

\newcommand{\bq}{\mathbf q}
\newcommand{\dd}{d_{\mathrm{dR}}}

\title[Diagonal harmonics and character polynomials]{Multivariate diagonal harmonics: character polynomials and proofs of two conjectures of Bergeron}
\author{Eren Ercan}
\email{eren321@gmail.com}
\date{31 October 2026}
\subjclass[2020]{Primary 05E05, 13A50; Secondary 05E10, 13D02, 20C30}
\keywords{diagonal harmonics, diagonal coinvariants, hook components, coefficient lengths, tensor coinvariants, character polynomials, Vandermonde operators}

\microtypesetup{expansion=false}
\allowdisplaybreaks[2]
\begin{document}

\begin{abstract}
We prove Bergeron's Hook-Components and Coefficient-Length conjectures
for the Vandermonde polarization module and full diagonal harmonics in
every rank and number of alphabets. Hook-label characters are elementary
skewings of the alternating character, and the maximal Schur length in
label $\mu$ is $n-\mu_1$. We compute character and dimension polynomials
of diagonal coinvariants in every bosonic--fermionic superspace for
$D\ge3$, $n\ge\max\{D,2D-6\}$, answering Church--Ellenberg--Farb's
character and dimension questions.

For augmented commutative algebras in characteristic zero, we construct
natural two-Schur presentations and filtered free resolutions with one
Schur generator in each homological degree. For polynomial algebras,
specialization to orthogonal idempotents identifies the first unstable
correction for every nonempty tail $\tau$ and degree $D>|\tau|$.
The stable threshold $n\ge D+\tau_1$ is sharp when the row dimension
is at least $\ell(\tau)$. First parabolic algebras recover finite jets,
their multiplication and their maps.

Selected-label lifts give uniform Vandermonde generation through degree
$D$ for $D\ge3$, $n\ge2D-3$, and disjoint-label products generate every
first occurrence for positive integral row weights. Constructive degree
reflection determines all two-column moment homology and generates its
first syzygies from degree $n+1$; two-column Vandermonde generation reduces
to degree $n-1$. Invariants and harmonics separate through the character
range and through degree $n$ in two alphabets. The hook identities also
give the fully graded sign component of the $t=0$, $(1,2)$-superspace
Theta conjecture.
\end{abstract}
\maketitle

\section{Introduction}

Bergeron's Hook-Components and Coefficient-Length conjectures
\cite[Conjectures~2 and~3]{Bergeron} describe how the symmetric-group
isotypes of multivariate diagonal harmonics fit together.  The first
recovers every hook-label character from the alternating character by
elementary skewing.  The second gives the exact number of alphabets
needed to retain every Schur coefficient in an arbitrary label
isotype.  We prove both conjectures for the Vandermonde polarization
module in all ranks and for full diagonal harmonics, with derivative
maps and finite-alphabet reconstruction.

We also answer Church--Ellenberg--Farb's eventual-character and
dimension-polynomial questions \cite[Problems~1.10 and~1.12]{ChurchEllenbergFarb2015},
construct natural Schur presentations and resolutions for tensor
coinvariants of arbitrary augmented algebras, and prove uniform
Vandermonde generation ranges. An idempotent specialization identifies
every first unstable character correction.
The hook transfer gives the fully graded Theta sign component.

\subsection{The hook-component and coefficient-length conjectures}

Let $K$ be a field of characteristic zero, let $n,k\ge1$, and put $V=K^k$.
Define
\[
 P_{n,k}=\Sym(V\otimes K^n)
       =K[x_i^{(a)}:i\in[n],\ a\in[k]],\qquad
 R_{n,k}=P_{n,k}/\langle(P_{n,k}^{S_n})_+\rangle.
\]
The row action of $\GL(V)$ commutes with the label action of $S_n$.
We use $\lambda$ for row partitions and $\mu$ for label partitions,
with Schur modules $\mathbb S_\lambda V$ and Specht modules $S^\mu$:
\begin{equation}\label{eq:character-decomposition}
 R_{n,k}\cong
 \bigoplus_{\mu\vdash n}\ \bigoplus_{\ell(\lambda)\le k}
 c_{\lambda,\mu}^{(n)}\,
                 \mathbb S_\lambda V\otimes S^\mu.
\end{equation}
Let $\Harm_{n,k}$ be the apolar harmonic model of $R_{n,k}$, and
let $\Pol_{n,k}$ be generated from the Vandermonde determinant
$\Delta_n(x^{(1)})$ by all partial
derivatives and generalized polarizations
\[
 E_{a,b}^{(r)}=\sum_{i=1}^n
             x_i^{(a)}\partial_{x_i^{(b)}}^r,\qquad r\ge1.
\]
This is Bergeron's module $\mathcal E_n^{\langle k\rangle}$
\cite[Section~2]{Bergeron}, in the general polarization framework of
\cite[Section~3.1]{r5-Blandin}.  Haiman's operator theorem gives
$\Pol_{n,2}=\Harm_{n,2}$~\cite[Theorem~4.2]{Haiman}.
For $\mathcal M=\Harm$ or $\Pol$, let $C_\mu^{\mathcal M}$ be the
stable row character of its label-$\mu$ multiplicity and put
$A_n^{\mathcal M}=C_{(1^n)}^{\mathcal M}$.  Here stability refers to
adjoining alphabets, with $n$ fixed, and $e_j^\perp$ is adjoint to
multiplication by $e_j$ for the Hall inner product.

\begin{maintheorem}[Bergeron's Hook-Components and Coefficient-Length conjectures]
\label{thm:intro-hooks}
For every $n\ge1$ and $\mathcal M=\Harm,\Pol$,
\begin{align*}
 C_{(j+1,1^{n-j-1})}^{\mathcal M}
          &=e_j^\perp A_n^{\mathcal M} &&(0\le j<n),\\
 \max\{\ell(\lambda):[s_\lambda]C_\mu^{\mathcal M}\ne0\}
          &=n-\mu_1 &&(\mu\vdash n).
\end{align*}
In $k\ge1$ alphabets the largest length is
$\min\{k,n-\mu_1\}$.  The empty row partition has length zero.
\end{maintheorem}

The length bound is attained by generated row-hook multiplicity
spaces (Theorem~\ref{thm:length}). The hook identities are realized
by exteriorization and derivative maps
(Theorems~\ref{thm:full-hook} and~\ref{thm:r2-operator-hook-transfer}).

The same calculation proves the nontrivial hook-label cases of
Bergeron's Conjecture~4, with the hook--hook formula
\eqref{eq:hook-hook-formula}. It uses the Fields character theorem
\cite{MRW} and the superharmonic operator theorem \cite{RW}.

The identities also make Bergeron's finite-alphabet reconstruction
unconditional: for $n\ge3$, the hook-label characters in
$\lfloor(n-1)/2\rfloor$ alphabets determine $A_n^{\mathcal M}$
(Theorem~\ref{thm:reconstruction}). The support bound
$\lambda'_1+\lambda'_2\le n-1$ gives alternating and hook-label
generation for $n\le6$ (Corollary~\ref{cor:r2-defect-support}).

For one bosonic and two fermionic alphabets, the transfer proves
the fully graded sign component of the $t=0$ Theta conjecture
\cite[Conjecture~8.2]{DIV2021}
(Theorem~\ref{thm:one-boson-two-fermion-sign}). At $q=1$ its character
is $(1+u+v)^{n-1}$ and its dimension is $3^{n-1}$.
The known prediction \cite[Theorem~5.6]{INV2024},
recorded in \cite[Theorem~8.3]{Lentfer12}, is realized by the hook
identity and diagonal supersymmetry \cite{Lentfer}.

For a nontrivial label $\mu=(c,\tau)$, the maximal-length row-hook
multiplicities form $\mathbb S_{\tau'}K^c$, with principal grading
shifted by one. Explicit harmonic operators realize its auxiliary
$\GL_c$ action and Lefschetz maps
(Theorem~\ref{thm:maximal-hook-action}).

\subsection{The character questions of Church--Ellenberg--Farb}

Problem~1.10 of \cite{ChurchEllenbergFarb2015} asks for the
multigraded symmetric-group characters in every number of alphabets,
at least for sufficiently large $n$.  Problem~1.12 asks for explicit
dimension polynomials whose existence follows from FI-module theory.
Write $X_j(\sigma)$ for the number of $j$-cycles of $\sigma\in S_n$.

\begin{maintheorem}[Full characters and dimension polynomials]
\label{thm:intro-full-characters}
Let $\mathbf b=(b_1,\ldots,b_k)\in\mathbb Z_{\ge0}^k$ have total degree $D\ge3$.
If $n\ge\max\{D,2D-6\}$, then
\begin{equation}\label{eq:intro-full-character}
 \operatorname{tr}\!\left(\sigma\mid(R_{n,k})_{\mathbf b}\right)
 = [\mathbf q^{\mathbf b}]
   \prod_{\substack{\alpha\in\mathbb Z_{\ge0}^k\\0<|\alpha|\le D}}
        (1-\mathbf q^\alpha)
   \prod_{j=1}^{D}\prod_{a=1}^{k}
                         (1-q_a^j)^{-X_j(\sigma)}.
\end{equation}
The moment Koszul complex is exact in positive homological degrees
in this range.  The formula also holds for $D\le1$, $n\ge1$,
and for $D=2$, $n\ge2$.
\end{maintheorem}

At the identity this gives the dimension polynomial in every
multidegree, with a range uniform in the number of alphabets.
Church--Ellenberg--Farb's theory establishes
eventual polynomiality and representation stability; subsequent
quantitative bounds include \cite{Bahran2024,Borisov2026}.
This establishes Bergeron's proposed low-degree Euler product
\cite[Section~4, equation~(33)]{Bergeron2013} in the displayed range.

The quadratic model also computes the first-occurrence degree and
the next two degrees of every label in every rank, and the entire two-column row sector
(Corollary~\ref{cor:first-three-character-layers}). These degrees
may lie beyond the uniform stable range.

Orthogonal Littlewood theory \cite{SSW2013} and the Schur relation
confine positive homology to boundary rows. Exactness in the first
possible degree removes an entire boundary before any case calculation. Wang's dimension polynomial
\cite[Theorem~5.1]{Wang2026Polynomiality} and a free cover give bivariate
exactness; the Fields character treats row hooks. At $n=2D-6$,
one explicit nonzero cubic image finishes the proof, independently
of operator generation. Supersymmetry transfers the formula to every
bosonic--fermionic superspace \cite{Lentfer}.

In every rank, two-column moment homology is acyclic at degree $n$
and reflects degrees $D>n$ to $2n-D$ by an explicit cycle map
(Theorem~\ref{thm:two-column-reflection}). Under this reflection,
multiplication becomes symplectic contraction. Thus all two-column
first syzygies are generated from their degree-$(n+1)$ layer
(Corollary~\ref{cor:two-column-syzygy-generation}).
Invariants and harmonics separate through the character range,
and through degree $n$ on row hooks and partitions of length at most two
(Corollary~\ref{cor:uniform-separation}).

\subsection{Natural stable multiplicity spaces}

The underlying presentation applies to every augmented commutative
algebra $B=K\oplus\mathfrak m$ in characteristic zero. The following
stable consequence requires only a connected positive grading.  Define
\[
 Q_n(B)=B^{\otimes n}/
       \bigl\langle(B^{\otimes n})^{S_n}_+\bigr\rangle.
\]
For a partition $\tau\vdash s$, write $\tau[n]=(n-s,\tau)$
when this is a partition, and set $\tau_1=0$ for the empty partition.

\Needspace{9\baselineskip}
\begin{samepage}
\begin{maintheorem}[Natural stable multiplicities]
\label{thm:intro-stable-spaces}
For $D\ge s$ and $n\ge D+\tau_1$, there is an isomorphism
\[
 \Hom_{S_n}\!\left(S^{\tau[n]},Q_n(B)_D\right)
                   \cong(\mathbb S_\tau\mathfrak m)_D.
\]
After fixing the one-dimensional Pieri branching line, the
isomorphism is natural for graded augmented algebra maps and
auxiliary symmetries preserving the grading and augmentation.
\end{maintheorem}
\end{samepage}

Multiplication enters these multiplicity spaces through the canonical
relation map $\beta_{\mu,B}$, whose source for $\mu=(t,\tau)$ is
$\mathbb S_{(t+1,\tau_2,\ldots)}\mathfrak m$. In a connected positive
grading, its first possible contribution has degree $n+1-\tau_1$.
The map commutes with arbitrary augmented algebra maps and vanishes
when $\mathfrak m^2=0$ (Theorem~\ref{thm:two-schur-presentation}).
For $B=\Sym V$, the multiplicity character is
$[s_\tau[\sum_{j\ge1}h_j(V)]]_D$.

For nonempty $\tau$, the first correction is exact in \emph{every} degree $D>s$.
At $n=D+\tau_1-1$, the kernel of the canonical Schur cover is
\[
                \mathbb S_{(D-s+\tau_1,\tau_2,\ldots)}V.
\]
It is nonzero whenever $\dim V\ge\ell(\tau)$, proving sharpness of
the stable threshold in every such degree
(Theorem~\ref{thm:first-unstable-specht-wall}). The proof specializes
to orthogonal idempotents: their tensor moment complex contracts on
every nonzero color. Naturality then detects the polynomial relation.

The relation extends to a natural filtered resolution with one
Schur generator in each homological degree
(Theorem~\ref{thm:full-schur-resolution}, proved in
Appendix~\ref{app:schur-resolution}). The idempotent argument gives
injectivity at every first possible polynomial degree except the
invariant first relation (Theorem~\ref{thm:polynomial-schur-fronts}).

The first parabolic algebra retains multiplication:
\[
                  Q_n(B)^{S_{n-1}}\cong B/\mathfrak m^n.
\]
Equivariant maps recover the augmented jet maps
(Theorem~\ref{thm:finite-jet-full-faithfulness}), and compatible towers
recover the formal completion. Thus these algebras retain the
multiplication that stable multiplicities forget.

\subsection{Constructive Vandermonde generation}

Selected-label evaluation and disjoint products construct harmonic
lifts. Write $\mathcal M_{n,k}[D]_\mu$ for the label-$\mu$ isotype
in total degree $D$.

\begin{maintheorem}[Uniform multivariate Vandermonde generation]
\label{thm:intro-linear}
In every number of alphabets, either condition
\[
 \mu_1\ge D-1
 \qquad\text{or}\qquad
 D\le\mu_1+\mu_3+\mu_5+\cdots
\]
implies $\Pol_{n,k}[D]_\mu=\Harm_{n,k}[D]_\mu$.
For $D\ge3$,
\begin{equation}\label{eq:intro-linear-range}
 n\ge2D-3
 \quad\Longrightarrow\quad
 \Pol_{n,k}[d]=\Harm_{n,k}[d]\qquad(0\le d\le D).
\end{equation}
For every $n,k$, every label, and every strictly positive integral
row-weight vector, the full label isotype at its first weighted
occurrence is generated. The first occurrence
means the least weighted degree in which the label appears.
\end{maintheorem}

Thus every degree through $\lfloor(n+3)/2\rfloor$ is generated,
uniformly in the number of alphabets.
A colored determinant map turns the curvilinear Quot degeneration
of B\'erczi--Svendsen \cite{r5-BS} into disjoint-label multiplication
and every positive-weight first occurrence.
Pair kernels save auxiliary labels, column alternation gives the
odd-row criterion, and the hook identities leave one family,
$(D-2,D-2,1)$, settled by an explicit rank-five lift.

Two-column Vandermonde generation reduces to the single degree
$n-1$. For $n\ge6$, it suffices to test the remaining rows in
$n-3$ alphabets; the cases $n\le5$ are already generated
(Proposition~\ref{prop:two-column-middle-test}).

\subsection{Reading routes}
Selected elimination, harmonic forms and the commuting model lead to the
conjecture proofs in Section~\ref{sec:augmented}. The Theta application
follows immediately; the auxiliary action and tensor-coinvariant results
form the next independent branches.
The canonical Schur presentation and its stable and first-unstable
consequences occupy Section~\ref{sec:tensor-applications}. The character
proof in Section~\ref{sec:stable-characters} combines them with quadratic,
bivariate and row-hook exactness, independently of operator generation.
The finer rowwise trace arguments are kept in
Appendix~\ref{app:finer-character-bounds}. Disjoint products then feed the local lifts
of Section~\ref{sec:linear-generation}.

\setcounter{tocdepth}{1}
\tableofcontents
\section{Selected-label elimination}\label{sec:selected-quotients}

Let $K$ be a field, let $n!$ be invertible in $K$, and let
$\varepsilon:B\to K$ be an augmented commutative $K$-algebra, with
$\mathfrak m=\ker\varepsilon$.  Put
\[
 Q_n(B)=B^{\otimes n}/
  \big\langle(B^{\otimes n})^{S_n}\cap
                  \ker(\varepsilon^{\otimes n})\big\rangle.
\]
For $b\in B$, write $b^{(i)}$ for its image in the $i$th tensor
factor.  The notation $S_{n-r}$ below means the subgroup fixing the
first $r$ tensor factors.  We call these $r$ factors selected and
the remaining factors auxiliary.  When $B$ is the coordinate algebra
of a pointed affine scheme, $Q_n(B)$ is its diagonal coinvariant
algebra at the basepoint.

Traces on the auxiliary factors can be eliminated using the global invariant traces.
The relations left after elimination are complete symmetric
functions on the selected factors.  The proof starts from the
fundamental trace presentation \cite[Theorem~1]{Domokos2017}; its
orbit-sum form makes both the resulting relations and the
characteristic hypothesis transparent.

\begin{theorem}[Exact selected-label presentation]
\label{thm:universal-selected-trace-presentation}
For $1\le r\le n$, set $m=n-r+1$.  Inclusion of the first $r$
tensor factors induces a natural algebra isomorphism
\begin{equation}\label{eq:universal-selected-trace-presentation}
 Q_n(B)^{S_{n-r}}
 \cong
 \frac{B^{\otimes r}}
 {\big\langle
     h_m\bigl(b^{(1)},\ldots,b^{(r)}\bigr):b\in\mathfrak m
       \big\rangle}.
\end{equation}
Equivalently, one may use every multilinear polarization of the
displayed relation, allowing arbitrary inputs from $\mathfrak m$.
The isomorphism is equivariant for $S_r$ and for every group acting
on $B$ by augmentation-preserving algebra automorphisms.
\end{theorem}

\begin{proof}
Let $0\le k\le n$ and put $\operatorname{TS}^k(B)
=(B^{\otimes k})^{S_k}$.  In the polynomial algebra
$\operatorname{Sym}_K(\mathfrak m)$ write $T_b$ for the generator
depending $K$-linearly on $b\in\mathfrak m$.  For $j\ge1$ define
\[
 \Psi_j(b_1,\ldots,b_j)=
 \sum_{\pi\in\Pi_j}(-1)^{j-|\pi|}
       \prod_{C\in\pi}(|C|-1)!
       \prod_{C\in\pi}T_{\prod_{i\in C}b_i},
\]
where $\Pi_j$ is the set of set partitions of $\{1,\ldots,j\}$.
Under $T_b\mapsto\sum_{i=1}^k b^{(i)}$, inclusion--exclusion on
the partition lattice identifies $\Psi_j$ with
\[
 \sum_{\substack{i_1,\ldots,i_j\in\{1,\ldots,k\}\\
                  \text{pairwise distinct}}}
      b_1^{(i_1)}\cdots b_j^{(i_j)}.
\]
In particular $\Psi_{k+1}$ maps to zero.  Multiplication by one
further trace gives the recursion
\begin{equation}\label{eq:trace-recursion-universal}
 \Psi_{j+1}(b_1,\ldots,b_j,c)
 =T_c\Psi_j(b_1,\ldots,b_j)
  -\sum_{a=1}^j
      \Psi_j(b_1,\ldots,b_{a-1},b_ac,b_{a+1},\ldots,b_j)
\end{equation}
which shows that all $\Psi_j$ with $j>k$ follow from
$\Psi_{k+1}$ and that the quotient by these relations is spanned
by $\Psi_j$, $0\le j\le k$ (with $\Psi_0=1$).

Choose a $K$-basis of $\mathfrak m$ and append $1$ to obtain a
basis of $B$.  The images of these $\Psi_j$, indexed by multisets
of $j$ basis elements, are nonzero scalar multiples of the orbit
basis of $\operatorname{TS}^k(B)$.  Their scalar factors are
products of factorials at most $k!$, and hence are invertible.
They are therefore linearly independent.  This proves
\[
 \operatorname{TS}^k(B)
 \cong\operatorname{Sym}_K(\mathfrak m)/
       \langle\Psi_{k+1}(b_1,\ldots,b_{k+1})\rangle.
\]
In particular, symmetric tensors are generated by the traces of
elements of $\mathfrak m$.

Now take $k=n-r$.  Before imposing the full $S_n$-invariants, the
$S_k$-invariant algebra is
$B^{\otimes r}\otimes\operatorname{TS}^k(B)$.
Reynolds averaging shows that the intersection of the full
invariant-generated ideal with this algebra is generated there by
the same invariant traces.  Thus the quotient imposes
\[
 T_b=-p_r(b),\qquad
 p_r(b)=\sum_{i=1}^r b^{(i)}\qquad(b\in\mathfrak m).
\]
Under this substitution every set-partition term of $\Psi_m$ has
sign $(-1)^m$.  Removing that common sign leaves
\[
 \sum_{\pi\in\Pi_m}
       \prod_{C\in\pi}(|C|-1)!
       \prod_{C\in\pi}p_r\Bigl(\prod_{i\in C}b_i\Bigr).
\]
The complete-symmetric generating identity identifies this with
the coefficient of $t_1\cdots t_m$ in
\[
 h_m\left(\sum_{a=1}^m t_ab_a^{(1)},\ldots,
           \sum_{a=1}^m t_ab_a^{(r)}\right).
\]
Polarization is valid because $m!$ is invertible.  Elimination of
all the $T_b$ therefore gives exactly
\eqref{eq:universal-selected-trace-presentation}.  Every construction
commutes with augmented algebra maps and with the stated actions.
\end{proof}

\section{Harmonic models and row characters}
\label{sec:preliminaries}

For the remainder of the paper, except where a more general base is
specified, $K$ is a field of characteristic zero.  We retain the
polynomial row convention and the notation $P_{n,k}$ and $R_{n,k}$
from the introduction.

\subsection{Harmonics and the operator module}

For $\alpha=(\alpha_1,\ldots,\alpha_k)\in\N^k$, define
\[
 D_\alpha=\sum_{i=1}^n
       \prod_{a=1}^k\partial_{x_i^{(a)}}^{\alpha_a}.
\]
The full diagonal harmonic space is
\begin{equation}\label{eq:harmonic-definition}
 \Harm_{n,k}=\bigcap_{\alpha\ne0}\ker D_\alpha.
\end{equation}
Polarized power sums generate $P_{n,k}^{S_n}$ in characteristic
zero.  Their apolar equations identify
\eqref{eq:harmonic-definition} with the orthogonal complement of
the Hilbert ideal for the Fischer pairing
\[
 \langle f,g\rangle=f(\partial)g\big|_{X=0}.
\]
Over the reals this pairing is positive definite on each
homogeneous component.  The quotient projection therefore
restricts to an isomorphism
\begin{equation}\label{eq:harmonic-quotient}
 \Harm_{n,k}\xrightarrow{\ \sim\ }R_{n,k}.
\end{equation}
The spaces and maps are defined over $\Q$, so the assertion
holds over $K$ by extension of scalars.  The differential
equations are stable under row changes; consequently the
isomorphism is $\GL_k\times S_n$-equivariant with our
polynomial row convention.

The ordinary Vandermonde is harmonic.  Partial derivatives
preserve \eqref{eq:harmonic-definition}, and
\begin{equation}\label{eq:polarization-commutator}
 [D_\alpha,E_{a,b}^{(s)}]
       =\alpha_a D_{\alpha-e_a+s e_b},
\end{equation}
where a zero coefficient gives zero.  This proves
$\Pol_{n,k}\subseteq\Harm_{n,k}$.
The operators with $s=1$ give the row Lie algebra.  As every
homogeneous part is a polynomial representation, $\Pol_{n,k}$
is a $\GL_k\times S_n$-submodule.  It is graded: the initial
Vandermonde is homogeneous and all generating operators are
homogeneous for total degree.

The linear harmonic equations say that every harmonic is
translation-invariant in each row.  Write
\[
 W_n=\{(z_1,\ldots,z_n)\in K^n:\ \textstyle\sum_i z_i=0\}.
\]
Using the invariant dot product to identify $W_n$ with its
dual, we can regard $\Harm_{n,k}$ and $\Pol_{n,k}$ as
derivative-closed submodules of $\Sym(K^k\otimes W_n)$.
In particular every occurring row partition has length at
most $n-1$ by the Cauchy decomposition.

\subsection{Stable row coefficients}

\begin{lemma}[Trailing-zero stability]\label{lem:stability}
For $\ell(\lambda)\leq k$, the highest-weight multiplicity
space of row weight $\lambda$ in $\Harm_{n,k}$ is naturally
the highest-weight multiplicity space of
$(\lambda,0)$ in $\Harm_{n,k+1}$.
\end{lemma}

\begin{proof}
A polynomial with last row degree zero does not use the last
row.  Every differential equation involving that row vanishes
on it.  The other equations are exactly the equations in
$k$ rows.  The additional row raising operators also differentiate
in the last row and hence vanish.  The remaining highest-weight
conditions are unchanged.
\end{proof}

Thus adjoining an unused row preserves every full harmonic
multiplicity in \eqref{eq:character-decomposition}.  The stable
character is finite: its row lengths
are bounded by $n-1$, and $R_{n,n-1}$ is finite-dimensional:
the invariant ideal contains the ordinary invariant ideal in
each row, so the quotient is a quotient of a tensor product
of ordinary coinvariant rings.
We write its label coefficients as symmetric functions
\[
 C_\mu^{(n)}=C_\mu^\Harm
       =\sum_\lambda c_{\lambda,\mu}^{(n)}s_\lambda,
 \qquad A_n=A_n^\Harm=C_{(1^n)}^{(n)}.
\]
Evaluation in $k$ variables recovers the character in $k$ row
alphabets.  For symmetric functions, $f^\perp$ denotes the
adjoint of multiplication by $f$ for the Hall inner product, in
which the Schur functions are orthonormal.

For $\Pol$, we use the direct system obtained by adjoining unused
rows.  For a fixed row partition its highest-weight spaces form
an increasing sequence inside the fixed finite-dimensional
full-harmonic space of Lemma~\ref{lem:stability}; the sequence
therefore stabilizes.  This defines the stable characters
\[
 C_\mu^{\Pol}=\sum_\lambda c_{\lambda,\mu}^{\Pol}s_\lambda,
 \qquad A_n^{\Pol}=C_{(1^n)}^{\Pol}.
\]
The commuting normal form in
Theorem~\ref{thm:r3-normal-form} will identify these stable
characters with their actual specialization in every finite
number of rows.

\subsection{A parabolic invariant quotient}

For a set of labels $T\subseteq[n]$, write $S_T$ for the
permutations supported on $T$.  Set
$X_i=(x_i^{(1)},\ldots,x_i^{(k)})$.
Passing to tail invariants leaves only the selected labels:
every tail power sum becomes the negative of the selected power
sum.  This elementary elimination controls row length throughout
the paper.

\begin{lemma}\label{lem:parabolic}
For $0\leq r\leq n$, the natural map
\[
 K[X_1,\ldots,X_r]\longrightarrow R_{n,k}^{S_{\{r+1,\ldots,n\}}}
\]
is surjective.  Consequently every row partition occurring
in the target has length at most $r$.
\end{lemma}

\begin{proof}
For $r\ge1$, surjectivity is the polynomial case of
Theorem~\ref{thm:universal-selected-trace-presentation}; for $r=0$
it is $R_{n,k}^{S_n}=K$.  Finally,
\[
 K[X_1,\ldots,X_r]
  =\Sym(K^k\otimes K^r)
  \cong\bigoplus_\lambda
       \Schur_\lambda K^k\otimes\Schur_\lambda K^r,
\]
and $\Schur_\lambda K^r=0$ for $\ell(\lambda)>r$.
\end{proof}

\begin{proposition}[The length upper bound]\label{prop:length-upper}
If $c_{\lambda,\mu}^{(n)}\ne0$, then
\[
 \ell(\lambda)\leq n-\mu_1.
\]
The same bound holds in every $\Pol_{n,k}$.
\end{proposition}

\begin{proof}
The restriction of $S^\mu$ to a subgroup $S_{\mu_1}$ contains
the trivial representation.  By the branching rule, one may
remove every box below the first row, in a sequence of removable
corners, and reach the one-row partition $(\mu_1)$.
Thus an occurrence of $\Schur_\lambda K^k\otimes S^\mu$
gives an occurrence of $\Schur_\lambda K^k$ in
$R_{n,k}^{S_{\mu_1}}$.
Apply Lemma~\ref{lem:parabolic} with $r=n-\mu_1$.
The assertion for $\Pol$ follows from its inclusion in
$\Harm\cong R$.
\end{proof}

\section{Harmonic forms and the Coefficient-Length conjecture}
\label{sec:length}

\subsection{The Euler splitting of harmonic forms}

A row hook has an alternating block of one-box rows below its first
row.  Replacing that block by exterior variables turns its
highest-weight condition into contraction by the Euler vector
field.  The de Rham identity then separates precisely the primitive
forms corresponding to row hooks.

On $K[x_1,\ldots,x_n]\otimes\Lambda(\theta_1,\ldots,\theta_n)$,
put
\[
 D_r=\sum_i\partial_{x_i}^r,\qquad
 \kappa_r=\sum_i\partial_{x_i}^r\partial_{\theta_i},
 \qquad
 \iota=\sum_i x_i\partial_{\theta_i},\qquad
 \dd=\sum_i\theta_i\partial_{x_i}.
\]
The superharmonic space $\SH_n$ is the joint kernel of $D_r$
for $r\geq1$ and $\kappa_r$ for $r\geq0$.
These are the apolar equations for the coinvariant algebra
of differential forms.  The invariant differential-form
description goes back to Solomon \cite{Solomon}.
Write $\SH_n^{a,j}$ for polynomial degree $a$ and exterior
degree $j$.

For $a\geq1$, let $K_{a,j}$ be the $S_n$-module of row
highest-weight vectors of weight $(a,1^j)$ in
$\Harm_{n,j+1}$.  For the boundary terms put
$K_{0,j}=0$ for $j>0$ and $K_{0,0}=K$.

\begin{lemma}[The row-hook bridge]\label{lem:hook-bridge}
There are $S_n$-equivariant identifications
\begin{equation}\label{eq:hook-bridge}
 K_{a,j}=\ker\iota\cap\SH_n^{a,j}.
\end{equation}
Away from bidegree $(0,0)$, there is a direct sum
\begin{equation}\label{eq:euler-splitting}
 \SH_n^{a,j}=K_{a,j}\oplus\dd K_{a+1,j-1},
\end{equation}
with nonexistent exterior degrees interpreted as zero.
\end{lemma}

\begin{proof}
Replace an exterior monomial by a determinant in $j$ auxiliary
rows:
\begin{equation}\label{eq:determinant-replacement}
 \theta_{i_1}\cdots\theta_{i_j}
 \longmapsto
 \det(y_{i_s}^{(r)})_{1\leq r,s\leq j}.
\end{equation}
This identifies exterior degree $j$ with polynomials
multilinear and alternating in the auxiliary rows.
Raising operators between auxiliary rows vanish on these
polynomials.  Raising an auxiliary row into the $x$-row
is precisely contraction by $\iota$.
Conversely, the auxiliary highest-weight conditions in
weight $(1,\ldots,1)$ select the determinant representation,
so every row-hook highest vector has this form.

The harmonic equations using only the $x$-row are $D_r$.
Those using exactly one auxiliary row are $\kappa_r$.
Any equation with two auxiliary derivatives at the same
label vanishes on an alternating determinant, and an equation
differentiating twice in an auxiliary row vanishes by
multilinearity.  These are all possible equations, proving
\eqref{eq:hook-bridge}.

Both $\iota$ and $\dd$ preserve $\SH_n$, since
\[
 [D_r,\dd]=0,\quad [D_r,\iota]=r\kappa_{r-1},\quad
 \kappa_r\dd+\dd\kappa_r=D_{r+1},\quad
 \kappa_r\iota+\iota\kappa_r=0.
\]
Moreover $\iota^2=\dd^2=0$ and
\[
 \iota\dd+\dd\iota=(a+j)\operatorname{id}
       \quad\hbox{on bidegree }(a,j).
\]
Thus every $f\in\SH_n^{a,j}$ with $a+j>0$ has the decomposition
\[
 f=\frac{\iota\dd f}{a+j}
       +\dd\left(\frac{\iota f}{a+j}\right).
\]
The first term belongs to $K_{a,j}$, and the expression in
parentheses belongs to $K_{a+1,j-1}$.  On this latter primitive
space, $\iota\dd=(a+j)\operatorname{id}$, so $\dd$ is injective
with inverse $\iota/(a+j)$ on its image.  The same identity shows
that its image meets $\ker\iota$ trivially.  This proves
\eqref{eq:euler-splitting}.
When $a=0$ and $j>0$, $\dd$ vanishes on the source;
the same identity gives $\ker\iota=0$, as required by
the boundary convention.
\end{proof}

\subsection{A tableau formula and generation of all row hooks}

For a standard tableau $T$, a descent is an integer $i$ for
which $i+1$ lies in a lower row than $i$.
Let $\des(T)$ be the number of descents and $\maj(T)$ their sum.

\begin{theorem}[Row-hook multiplicities]\label{thm:row-hooks}
For every $\mu\vdash n$ and $j\geq0$,
\begin{equation}\label{eq:row-hook-tableaux}
 \sum_{a\geq1}[K_{a,j}:S^\mu]q^a
 =
 \sum_{\substack{T\in\SYT(\mu)\\t=\des(T)\geq1}}
 q^{\maj(T)}e_j(q^{-2},q^{-3},\ldots,q^{-t}).
\end{equation}
Equivalently, the coefficient indexed by $(a,1^j)$ and $\mu$
counts pairs $(T,J)$ satisfying
\[
 J\subseteq\{2,\ldots,\des(T)\},\qquad
 |J|=j,\qquad \maj(T)-\sum_{r\in J}r=a.
\]
\end{theorem}

\begin{proof}
The Fields theorem of Murai, Rhoades and Wilson
\cite[Theorem~5.2 and equation~(2.5)]{MRW} gives
\begin{equation}\label{eq:fields-character}
 \Frob(\SH_n;q,u)
 =\sum_{\mu\vdash n}s_\mu[Z]
   \sum_{T\in\SYT(\mu)}
      q^{\maj(T)}
      \prod_{r=1}^{\des(T)}(1+u q^{-r}).
\end{equation}
The generating series of \eqref{eq:euler-splitting} says
that the nonconstant character is $(1+u/q)$ times the
generating series of the $K_{a,j}$.
The unique tableau without descents is the one-row tableau,
and contributes only the constant trivial representation.
For every other tableau the factor $1+u/q$ cancels.
Extracting $u^j$ gives \eqref{eq:row-hook-tableaux}.
All remaining powers of $q$ are positive:
if $t=\des(T)$, then
$\maj(T)\geq t(t+1)/2$, whereas
$\sum_{r=2}^t r=t(t+1)/2-1$.
\end{proof}

The harmonic-form bridge extracts the row-hook multiplicities from
the Fields character.  For the alternating label, the resulting
specialization also
appears in \cite[Corollary~4.6]{Lentfer}.
The row-hook coefficients after summing symmetric-group
dimensions are studied in \cite{CL}.

\begin{proposition}[Generation of the row-hook sector]
\label{prop:row-hook-generation}
Every row-hook isotypic component of $\Harm_{n,k}$ is
contained in $\Pol_{n,k}$.
\end{proposition}

\begin{proof}
For $r\geq1$ set $\Theta_r=\sum_i\theta_i\partial_{x_i}^r$.
The operator theorem of Rhoades and Wilson
\cite[Theorem~5.1]{RW} gives
\begin{equation}\label{eq:rw-operators}
 \SH_n=\sum_{J\subseteq[n-1]}\Theta_J\Harm_{n,1},
\end{equation}
where the $\Theta_r$ anticommute and a fixed order is used
in each product.
Every ordinary harmonic is a linear combination of partial
derivatives of $\Delta_n$.
Under \eqref{eq:determinant-replacement}, the form
$\Theta_{r_1}\cdots\Theta_{r_j}h$ becomes
\begin{equation}\label{eq:theta-polarization}
 \sum_{\sigma\in S_j}\sgn(\sigma)
 E_{r_1}^{(\sigma(1))}\cdots E_{r_j}^{(\sigma(j))}h,
 \qquad
 E_s^{(a)}=\sum_i y_i^{(a)}\partial_{x_i}^s.
\end{equation}
The even operators in this formula commute.
Expansion of both sides verifies the formula term by term.
It belongs to $\Pol_{n,j+1}$.
Thus the entire subspace
$\ker\iota\cap\SH_n^{a,j}$ is represented inside
$\Pol_{n,j+1}$.  By Lemma~\ref{lem:hook-bridge}, this
is the full row-hook highest-weight multiplicity space.
Closing under $\GL_{j+1}$, and then adjoining unused rows
and closing under $\GL_k$, proves the proposition.
\end{proof}

For a nontrivial hook label $\mu=(a+1,1^b)$, $a+b=n-1$, $b\ge1$,
the tableau formula and row-hook generation give
\begin{equation}\label{eq:hook-hook-formula}
 \frac{C_\mu^{\Pol}[q\mid u]}{q+u}
  =\begin{bmatrix}n-1\\b\end{bmatrix}_q
                 \prod_{r=2}^{b}(q^r+u).
\end{equation}
Here the supersymmetric specialization is characterized by
$p_r\mapsto q^r+(-1)^{r-1}u^r$, and the bracket is the Gaussian
binomial coefficient.  To verify the formula, every tableau of
this hook has $b$ descents, and their possible sets are all
$b$-element subsets of $[n-1]$.  Thus their major-index generating
polynomial is $q^{b(b+1)/2}\bigl[\begin{smallmatrix}n-1\\b\end{smallmatrix}\bigr]_q$.
Now use $s_{(c,1^j)}[q\mid u]=(q+u)q^{c-1}u^j$ and
Theorem~\ref{thm:row-hooks}.  All other nonempty row shapes vanish
under this specialization.  This proves the nontrivial hook-label
cases of Bergeron's Conjecture~4
\cite[Section~4.2, equation~(4.4)]{Bergeron}.

\subsection{Maximal row hooks and exact coefficient lengths}

At the longest possible hook length, the tableau formula forces the
maximum descent number. Transposing the tableaux turns this into a
minimum-descent condition, which is counted by a finite principal
specialization. This gives the whole multiplicity series at once.

\begin{theorem}[The complete maximal-length row-hook series]
\label{thm:maximal-hook-series}
Let $\mu=(c,\tau)\vdash n$ and $r=|\tau|>0$. Then, in both the
harmonic and Vandermonde modules,
\[
 H_\mu(q):=\sum_{a\ge1}[K_{a,r-1}:S^\mu]q^a
                         =q\,s_{\tau'}(1,q,\ldots,q^{c-1}).
\]
Thus all maximal-length row-hook coefficients are given by one finite
principal specialization of a Schur polynomial.
\end{theorem}
\begin{proof}
A tableau of shape $\mu$ has at most $r$ descents. In \eqref{eq:row-hook-tableaux}, exterior degree $r-1$ therefore forces exactly
$r$ descents and selects the whole set $\{2,\ldots,r\}$.
Transposition complements the descent set. If $N=\binom n2$, then
\[
 H_\mu(q)=q^{N-r(r+1)/2+1}
 \sum_{\substack{U\in\operatorname{SYT}(\mu')\\
                         \operatorname{des}(U)=c-1}}
      q^{-\operatorname{maj}(U)}
 =q^{N-r(r+1)/2+1}s_{\mu'}(1,q^{-1},\ldots,q^{-(c-1)}).
\]
For the second equality, use the standard principal-specialization
identity on exactly $c=\ell(\mu')$ variables. Standardization gives a
sum of $q^{\operatorname{maj}(U)}$ times Gaussian coefficients; only
the minimum descent number $c-1$ survives, with Gaussian coefficient
one. One may standardize using comajor index and use tableau evacuation
to replace comajor by major index.

Now $\mu'=(1^c)+\tau'$, with $\tau'$ padded by zeros. Factoring the
determinant contributes $q^{-\binom c2}$. Reverse the remaining $c$
variables and use homogeneity of $s_{\tau'}$, of degree $r$, to obtain
an additional factor $q^{-r(c-1)}$. Since
\[
 N-\frac{r(r+1)}2-\binom c2-r(c-1)=0,
\]
the total remaining factor is $q$, proving the formula.
\end{proof}

\begin{corollary}[A positive contribution at every permitted hook length]
\label{cor:shorter-hook-witnesses}
With $\mu=(c,\tau)$ and $r=|\tau|>0$, fix $1\le\ell\le r$.
In the row-hook coefficient series of length $\ell$, the contribution
of tableaux with exactly $r$ descents is
\[
 q^{\binom{r-\ell+2}{2}}
 \begin{bmatrix}r-1\\\ell-1\end{bmatrix}_q
                       s_{\tau'}(1,q,\ldots,q^{c-1}).
\]
All other contributions have nonnegative coefficients. Hence, if
$b=\sum_i\binom{\tau_i}{2}$, the row hook
\[
       \left(\binom{r-\ell+2}{2}+b,\,1^{\ell-1}\right)
\]
occurs. Taking $\ell=\min(k,r)$ gives the sharp coefficient length
in $k$ alphabets directly.
\end{corollary}
\begin{proof}
The proof of Theorem~\ref{thm:maximal-hook-series} gives
$\sum_{\operatorname{des}(T)=r}q^{\operatorname{maj}(T)}
=q^{r(r+1)/2}s_{\tau'}(1,q,\ldots,q^{c-1})$.
Multiply by
\[
 e_{\ell-1}(q^{-2},\ldots,q^{-r})
 =q^{-r(\ell-1)+\binom{\ell-1}{2}}
             \begin{bmatrix}r-1\\\ell-1\end{bmatrix}_q
\]
and simplify the exponent. Its least coefficient is positive.
\end{proof}

\begin{theorem}[Exact coefficient lengths]\label{thm:length}
For $\mu\vdash n$, put $r=n-\mu_1$.
The largest length of a row partition in the $\mu$-isotype of either
$\Harm_{n,k}$ or $\Pol_{n,k}$ is $\min(k,r)$.
The corresponding stable largest length is $r$.
\end{theorem}

\begin{proof}
Proposition~\ref{prop:length-upper} gives the upper bound. If
$\mu=(c,\tau)$ and $r=|\tau|>0$, take
$\ell=\min(k,r)$ in Corollary~\ref{cor:shorter-hook-witnesses}.
Its positive contribution contains a row hook of length $\ell$;
Proposition~\ref{prop:row-hook-generation} places the whole component
in $\Pol_{n,k}$. If $r=0$, exactness of taking $S_n$-invariants
makes the trivial isotype of $R_{n,k}$ equal to $K$, and constants
belong to $\Pol_{n,k}$.
\end{proof}

\begin{corollary}[Palindromicity, unimodality, and the extremal hook]
\label{cor:maximal-hook-coefficients}
Put $b=\sum_i\binom{\tau_i}{2}$. The coefficients of $H_\mu(q)$
are palindromic and unimodal on the interval
\[
             a_{\min}=1+b,\qquad a_{\max}=1+r(c-1)-b,
\]
with both endpoint coefficients equal to one. Their center is
$1+r(c-1)/2$. In particular the least-degree row hook of maximal
possible length $r$ is
\begin{equation}\label{eq:explicit-length-witness}
 \lambda_\mu=(1+b,1^{r-1}),\qquad c_{\lambda_\mu,\mu}^{(n)}=1.
\end{equation}
Its total degree is
$r+b=\sum_i\binom{\tau_i+1}{2}$. The sum of all maximal-length
row-hook multiplicities is
\[
 H_\mu(1)=\dim\mathbb S_{\tau'}K^c
 =\prod_{1\le i<j\le c}
                 \frac{\tau'_i-\tau'_j+j-i}{j-i}.
\]
An exact product form is
\[
 H_\mu(q)=q^{1+b}
   \prod_{(i,j)\in\tau}\frac{1-q^{c+i-j}}{1-q^{h_\tau(i,j)}}.
\]
\end{corollary}
\begin{proof}
The least principal-specialization exponent is
$n(\tau')=\sum_i\binom{\tau_i}{2}=b$, with coefficient one;
reversal of the $c$ variables gives the greatest exponent
$r(c-1)-b$, also with coefficient one. This proves the endpoints
and palindromicity.

For unimodality, apply the standard $\mathfrak{sl}_2$ criterion
\cite{Proctor1982,Stanley1989}: regard $K^c$ as the irreducible
$\mathfrak{sl}_2$-module $\operatorname{Sym}^{c-1}K^2$. Then
\[
 t^{-r(c-1)}s_{\tau'}(1,t^2,\ldots,t^{2(c-1)})
\]
is the $\mathfrak{sl}_2$ character of $\mathbb S_{\tau'}K^c$.
Its decomposition into irreducibles expresses its coefficients as
a sum of indicator functions of centered weight intervals. The
coefficients therefore increase weakly up to the center and decrease
weakly thereafter. Section~\ref{sec:hook-action} constructs the
operators on harmonic multiplicities. The Weyl dimension and hook-content
formulas give the last two identities.
\end{proof}

For $\mu=(3,3)$ this hook is $(4,1,1)$; for the alternating label
it is $(1^{n-1})$.

\begin{corollary}[Small row support and operator recovery]
\label{cor:near-row}
The stable character $C_\mu^{(n)}$ is determined by its
specialization to $n-\mu_1$ row variables.
If $\mu_1\geq n-2$, then the $\mu$-isotypes of
$\Pol_{n,k}$ and $\Harm_{n,k}$ agree for every $k$.
\end{corollary}

\begin{proof}
Schur functions of length at most $r$ remain linearly
independent in $r$ variables.  The length bound therefore
proves the first assertion.
For the second, Haiman's bivariate operator theorem
\cite[Theorem~4.2]{Haiman} gives $\Pol_{n,2}=\Harm_{n,2}$.
Every highest-weight vector in the specified labels
uses at most two rows, by the length bound and
Lemma~\ref{lem:stability}.  It is therefore generated
in those rows.  The row group generates its full
isotypic component in any larger number of rows.
The one-row assertion is classical.
\end{proof}

\begin{corollary}\label{cor:three-row-defect}
In three alphabets, a nonzero component of
$\Harm_{n,3}/\Pol_{n,3}$ can occur only at a pair
$(\lambda,\mu)$ satisfying
\[
 \lambda=(a,b,c),\qquad a\geq b\geq2,\quad c\geq1,
 \qquad \mu_1\leq n-3.
\]
\end{corollary}

\begin{proof}
All row partitions of length at most two are generated
by the bivariate theorem and trailing-zero stability.
The remaining row hooks are generated by
Proposition~\ref{prop:row-hook-generation}.
Corollary~\ref{cor:near-row} excludes all label types
with $\mu_1\geq n-2$.
\end{proof}

\section{A commuting operator model}
\label{sec:r3-normal-form}

All generalized polarizations act as vector fields on a single
family of pointed jet determinants.  An invertible source jet can
be reparametrized to $t$; the other jet coefficients then become
commuting polarizations from that source row.  Polynomial dependence
extends this description across the singular jets.  The resulting
normal form makes specialization in the number of alphabets exact.

\begin{lemma}[Normal ordering]\label{lem:r2-normal-order}
Let $V$ be a finite-dimensional row space over a field of
characteristic zero.  Suppose $\mathcal M(V)$ is generated from alternating harmonic
polynomials by partial derivatives and generalized polarizations.
Let $\mathcal B(V)$ be the span of the words using only generalized
polarizations on those generators.  Then
\[
 \mathcal M(V)=K[\partial]\mathcal B(V),
 \qquad
 \mathcal M(V)^{\sgn}=\mathcal B(V).
\]
In particular, the alternating part of $\Pol$ is generated from the
ordinary Vandermonde by generalized polarizations alone.
\end{lemma}

\begin{proof}
For $E=\sum_i x_i^{(a)}(\partial_{x_i^{(b)}})^r$,
\[
 [E,\partial_{x_j^{(c)}}]
       =-\delta_{ac}(\partial_{x_j^{(b)}})^r.
\]
More generally, the commutator of $E$ with any constant coefficient
differential operator is again such an operator.  Commuting partial
derivatives to the left, by induction on the number of polarization
operators and on their positions in a word, proves the first identity.
Every element of $\mathcal B(V)$ is alternating, because the
polarizations commute with $S_n$.

Apply the sign projector to an expression $\sum_\ell D_\ell g_\ell$,
where $g_\ell\in\mathcal B(V)$.  Since each $g_\ell$ is alternating,
the result is $\sum_\ell\Av_{S_n}(D_\ell)g_\ell$.  The positive-degree
invariant differential operators annihilate every diagonal harmonic.
Only the constant term of each averaged differential operator remains.
This proves the second identity.
\end{proof}

Put $N=\binom n2$, write $x_i=x_i^{(1)}$, and normalize
\[
 \Delta_n=\det(x_i^j)_{\substack{1\leq i\leq n\\0\leq j<n}}.
\]
All fields in this section have characteristic zero.  For $a\geq2$ put
$Y_{a,s}=\sum_i x_i^{(a)}\partial_{x_i}^{s}$.
The operators $Y_{a,s}$ and the individual $\partial_{x_i}$ commute.
Let $\mathcal B_{n,k}$ be the span of all words in generalized
polarizations applied to $\Delta_n$.

\begin{theorem}[Commuting normal form]\label{thm:r3-normal-form}
For every $n,k\geq1$,
\begin{align}
 \mathcal B_{n,k}
   &=K[Y_{a,s}:2\leq a\leq k,\ 1\leq s<n]\,\Delta_n,
       \label{eq:r3-pure-normal-form}\\
 \Pol_{n,k}
   &=K[\partial_{x_1},\ldots,\partial_{x_n},
           Y_{a,s}:2\leq a\leq k,\ 1\leq s<n]\,\Delta_n.
       \label{eq:r3-full-normal-form}
\end{align}
Moreover, $\Pol_{n,k}^{\sgn}=\mathcal B_{n,k}$.
Thus one source row and a commuting family of operators suffice
for the module defined using arbitrary generalized polarizations.
\end{theorem}

\begin{proof}
The assertion for $n=1$ is immediate, so assume $n\geq2$.
Let $\mathfrak t=tK[t]/(t^n)$, and for
$f=(f_1,\ldots,f_k)\in\mathfrak t^k$ define
\begin{equation}\label{eq:r3-jet-determinant}
 G_f(X)=\left(\prod_{j=0}^{n-1}j!\right)
 \det\left([t^j]\exp\left(\sum_{a=1}^k
                              x_i^{(a)}f_a(t)\right)\right)_{i,j}.
\end{equation}
Only coefficients of order less than $n$ are used.  Thus this is a
polynomial in $X$ and in the finitely many coefficients of $f$,
of $X$-degree at most $N$.
Let $J$ be the span of its values as $f$ varies.

The one-particle identity
\[
 x^{(a)}\partial_{x^{(b)}}^s e^{\sum_c x^{(c)}f_c(t)}
      =x^{(a)}f_b(t)^s e^{\sum_c x^{(c)}f_c(t)}
\]
and differentiation of a determinant give
\begin{equation}\label{eq:r3-jet-vector-field}
 E_{a,b}^{(s)}G_f
     =\left.\frac{d}{d\epsilon}
               G_{f+\epsilon e_a f_b^s}\right|_{\epsilon=0}.
\end{equation}
The product $f_b^s$ is taken modulo $t^n$ and has zero constant
term.  A derivative of a polynomial family lies in the span of its
values, by polynomial interpolation over an infinite field.
Consequently $J$ is stable under every generalized polarization,
including those with output in the source row.
It contains $G_{(t,0,\ldots,0)}=\Delta_n$, and hence
$\mathcal B_{n,k}\subseteq J$.

Suppose that $f_1'(0)\ne0$, and let $g\in\mathfrak t$ be its
compositional inverse modulo $t^n$.
Composition by $g$ acts on coefficients of degrees $0,\ldots,n-1$
through a triangular matrix with diagonal
$1,g'(0),\ldots,g'(0)^{n-1}$.  Therefore
\begin{equation}\label{eq:r3-reparametrization}
 G_{f\circ g}=g'(0)^N G_f.
\end{equation}
Since $f_1\circ g=t$, write
$f_a\circ g=\sum_{s=1}^{n-1}c_{a,s}t^s$ for $a\geq2$.
Expanding the commuting operator exponential gives
\begin{equation}\label{eq:r3-normalized-jet}
 G_{f\circ g}
   =\exp\left(\sum_{a=2}^k\sum_{s=1}^{n-1}
                       c_{a,s}Y_{a,s}\right)\Delta_n.
\end{equation}
The exponential is finite on $\Delta_n$, since every $Y_{a,s}$
lowers its source-row degree by $s$.
Thus $G_f$ belongs to the right side of
\eqref{eq:r3-pure-normal-form} whenever $f_1'(0)\ne0$.
That condition defines a dense open subset of the jet parameter
space.  A polynomial family taking values in a fixed linear
subspace on a dense subset takes values in that subspace everywhere:
apply any linear functional vanishing on the subspace.
It follows that $J$ is contained in the right side of
\eqref{eq:r3-pure-normal-form}.  The reverse inclusion is immediate,
which proves the first identity.

Lemma~\ref{lem:r2-normal-order} gives
$\Pol_{n,k}=K[\partial_{x_i^{(a)}}]\mathcal B_{n,k}$.
Substitute \eqref{eq:r3-pure-normal-form} and move every
derivative in an auxiliary row to the right, using
\[
 [\partial_{x_i^{(a)}},Y_{b,s}]
           =\delta_{ab}\partial_{x_i}^{s}\qquad(a,b\geq2).
\]
Those derivatives kill $\Delta_n$ when they reach it.
The remaining operators are exactly those in
\eqref{eq:r3-full-normal-form}.

The same normal-ordering lemma identifies the alternating part
with $\mathcal B_{n,k}$.
\end{proof}

Positive-order polarizations preserve the pointed condition
$f_a(0)=0$ in \eqref{eq:r3-jet-vector-field}.  Individual
translations of the source variables encode partial derivatives
and give the following finite generating family.

\begin{corollary}[A finite generating determinant]
\label{cor:r3-generating-determinant}
The coefficient span, in the parameters $u_i$ and $c_{a,s}$, of
\begin{equation}\label{eq:r3-generating-determinant}
 \left(\prod_{j=0}^{n-1}j!\right)
 \det\left([t^j]\exp\left((x_i+u_i)t+
         \sum_{a=2}^k x_i^{(a)}
                  \sum_{s=1}^{n-1}c_{a,s}t^s\right)\right)_{i,j}
\end{equation}
is exactly $\Pol_{n,k}$.
Without the parameters $u_i$, the coefficient span is its alternating
part.  A spanning set consists of the words
\[
 \partial_x^\alpha\prod_{a,s}Y_{a,s}^{m_{a,s}}\Delta_n,
 \qquad |\alpha|+\sum_{a,s}s\,m_{a,s}\leq N.
\]
\end{corollary}

\begin{proof}
The determinant is
$\exp(\sum_i u_i\partial_{x_i}+\sum_{a,s}c_{a,s}Y_{a,s})\Delta_n$.
Use Theorem~\ref{thm:r3-normal-form} and compare coefficients.
The displayed inequality is the available source-row degree.
\end{proof}

\begin{corollary}[Exact specialization in the number of rows]
\label{cor:r3-exact-specialization}
If $r\geq k$, setting rows $k+1,\ldots,r$ equal to zero maps
$\Pol_{n,r}$ onto $\Pol_{n,k}$.  Moreover,
\[
 \Pol_{n,r}\cap K[x_i^{(a)}:a\leq k]=\Pol_{n,k}.
\]
The stable row character specializes to the actual finite-row module
for every $k$.
\end{corollary}

\begin{proof}
In the commuting normal form, the auxiliary-row degree of a word is
exactly the number of its corresponding $Y$ factors.  Setting a row
to zero kills precisely the words using a polarization into that row.
All remaining words belong to $\Pol_{n,k}$, and every such word
already occurs in $\Pol_{n,r}$.
If a polynomial is independent of the discarded rows, restriction
fixes it, proving the intersection assertion.
For a polynomial row representation, the zero-weight subspace in
the discarded rows is the usual specialization of its Schur character.
\end{proof}

Haiman's bivariate operator theorem
\cite[Theorem~4.2]{Haiman} describes the full bivariate harmonic
space using the commuting operators in
\eqref{eq:r3-full-normal-form}.  Bergeron's multivariate module
\cite[Section~2.1]{Bergeron} is defined by closure under derivatives
in all variables and polarizations between arbitrary pairs of rows.
Theorem~\ref{thm:r3-normal-form} identifies these operator
descriptions in every number of rows and makes their finite-row
specialization exact.

\section{The Hook-Components conjecture and reconstruction}
\label{sec:augmented}

The sign representation turns an exterior degree into a hook label.
To use this observation on coinvariants, we show that every
positive odd trace already vanishes on the sign component once the
even traces have been removed.  The argument uses multiplication
and augmentation alone, so we work at the level of an arbitrary
pointed algebra.

This section works over any commutative unital ring $K$ in which $n!$ is
invertible.  Let $A$ be a commutative unital $K$-algebra
with an augmentation $\eps:A\to K$, and let
$\mathfrak m=\ker\eps$.  Set
\[
 B=A^{\otimes n},\qquad
 \eps_B=\eps^{\otimes n}.
\]
The symmetric group permutes the tensor factors.
For $a\in A$, let $a^{(i)}$ be $a$ in the $i$th tensor
factor and $1$ in all other factors.
Let $N_n$ be the permutation representation with basis
$\theta_1,\ldots,\theta_n$, and set
\[
 S=B\otimes\Lambda N_n.
\]
Its augmentation $\eps_S$ extends $\eps_B$ and sends
every $\theta_i$ to zero.

\begin{definition}\label{def:augmented-coinvariants}
The ordinary and exterior coinvariant quotients at the
augmentation are
\[
 I_A=B\,\ker(\eps_B|_{B^{S_n}}),\qquad R_A=B/I_A,
\]
and
\[
 J_A=S\,\ker(\eps_S|_{S^{S_n}}),\qquad Q_A=S/J_A.
\]
\end{definition}

If $A$ is connected and graded, these are the usual
quotients by invariants with zero constant term.
The augmentation definition also makes sense for
ungraded and nonreduced algebras.
Put $\theta_\Sigma=\sum_i\theta_i$.
As $n$ is invertible, $N_n/K\theta_\Sigma$ identifies
canonically with the standard representation $W_n$.

\begin{theorem}[Augmented sign transfer]\label{thm:augmented-transfer}
The natural quotient induces an isomorphism
\begin{equation}\label{eq:augmented-transfer}
 (R_A\otimes\Lambda W_n)^{\sgn}
       \xrightarrow{\ \sim\ }Q_A^{\sgn}.
\end{equation}
It preserves exterior degree, is natural under augmented
algebra homomorphisms, commutes with arbitrary base change, and
commutes with all augmented
automorphisms of $A$.  It also preserves any additional
grading respected by the augmentation.
\end{theorem}

For a polynomial algebra in two generators, this is the
sign-transfer theorem of Gonz\'alez, Lentfer and Mularczyk
\cite[Theorem~4.21]{GLM}.  Their complementary-averaging and
determinant argument \cite[Lemmas~4.13--4.15]{GLM} extends to the
present setting because collisions multiply elements of the
algebra and positive traces remain in its augmentation ideal.

\subsection{Trace generators}

For $a\in A$ and $T\subseteq[n]$, write
\[
 t_a(T)=\sum_{i\in T}a^{(i)},\qquad
 q_a=\sum_{i=1}^n\theta_i a^{(i)},\qquad
 t_a=t_a([n]).
\]

\begin{lemma}\label{lem:trace-generation}
The algebra $B^{S_n}$ is generated by $1$ and $t_a$
for $a\in\mathfrak m$.
The superalgebra $S^{S_n}$ is generated by these even
elements together with $q_a$ for $a\in A$.
Consequently
\begin{equation}\label{eq:trace-ideals}
 I_A=\langle t_a:a\in\mathfrak m\rangle,\qquad
 J_A=\langle t_a:a\in\mathfrak m;\
                    q_a:a\in A\rangle.
\end{equation}
\end{lemma}

\begin{proof}
Use the split $K$-module decomposition $A=K1\oplus\mathfrak m$.
Pure tensors in $B$, together with exterior monomials,
span $S$ as a $K$-module; averaging spans its invariant submodule.
An averaged tensor can be described by finitely many
occupied labels, each carrying an even decoration from
$\mathfrak m$ or an odd decoration $\theta a$ with
$a\in A$.

Multiply the traces corresponding to those decorations.
The terms with pairwise distinct labels give the desired
orbit sum, up to the multiplicity factorials of repeated
even decorations.  These factorials divide $n!$.
A repeated odd decoration has zero distinct-label
antisymmetrization, so no division by a vanishing scalar
is required.  Every other term has a collision and uses
fewer occupied labels.  Two odd decorations at one label
give zero.  An even collision is multiplication in $A$,
and a collision involving an even decoration from
$\mathfrak m$ stays in $\mathfrak m$.
Induction on the number of occupied labels proves the
two generation statements.  All generators except $1$
have zero augmentation, giving \eqref{eq:trace-ideals}.
\end{proof}

For a set $T$ of labels, let $\Av_T$ and $\Alt_T$ denote
the normalized symmetric and alternating averaging
operators over $S_T$.
The averaging identity below proves sign transfer. A trace relation
moves a factor from the complementary labels onto an
alternating block, where the resulting determinant has a repeated
row and vanishes.

\begin{lemma}[Complementary averaging]\label{lem:complement-average}
Let $D\subseteq[n]$, let $i\notin D$, and put
$C=D\cup\{i\}$.  Suppose $f\in B$ is alternating on
$C$ and symmetric on $C^c$.
For every $a\in\mathfrak m$,
\begin{equation}\label{eq:complement-average}
 \Av_{D^c}\bigl(a^{(i)}f\bigr)\in I_A.
\end{equation}
\end{lemma}

\begin{proof}
Separate the tensor factors on $C$ and $C^c$.
The symmetric and alternating averaging operators are idempotents
on separate tensor factors.  Their images are split direct summands;
tensoring these splittings identifies the module with the indicated
symmetries as the tensor product of the alternating module on $C$
and the invariant module on $C^c$.  No flatness is required.
By Lemma~\ref{lem:trace-generation}, the latter is
generated by the traces $t_b(C^c)$ with $b\in\mathfrak m$.
Modulo $I_A$,
\[
 t_b(C^c)=-t_b(C).
\]
After making these substitutions, $f$ is congruent to
a sum of elements using only the tensor factors on
$C$ and alternating on $C$.  Multiplication by $a^{(i)}$
and averaging preserve congruence modulo $I_A$.

Write $c_1,\ldots,c_h$ for the labels of $C$, with
$c_h=i$.
Alternating tensors on these factors are spanned as a $K$-module by
determinants
\[
 F=\det\bigl(b_j^{(c_\ell)}\bigr)_{1\leq\ell,j\leq h},
        \qquad b_j\in A.
\]
Indeed, the alternating projector sends a pure tensor to its
determinant divided by $h!$.  Thus it suffices to prove \eqref{eq:complement-average}
for $F$.
Let $M_j$ be the cofactor of the entry in its last row
and $j$th column.  For $r\in D^c$, let $F_r$ be the
same determinant with its last row evaluated at label
$r$.  The cofactor $M_j$ uses labels in $D$ only.
It follows that
\begin{align*}
 |D^c|\,\Av_{D^c}(a^{(i)}F)
   &=\sum_{r\in D^c}a^{(r)}F_r\\
   &=\sum_{j=1}^h M_j\,t_{ab_j}(D^c)\\
   &\equiv-\sum_{r\in D}a^{(r)}
                     \sum_{j=1}^h M_j b_j^{(r)}
                          \pmod{I_A}\\
   &=0.
\end{align*}
The congruence uses $ab_j\in\mathfrak m$.
Each final inner sum is a determinant with a repeated
row.  Since $|D^c|$ is invertible, the assertion follows.
\end{proof}

\subsection{Odd invariants on the sign component}

\begin{lemma}\label{lem:odd-annihilation}
For $a\in\mathfrak m$, multiplication by $q_a$
annihilates the sign submodule of
$R_A\otimes\Lambda N_n$.
\end{lemma}

\begin{proof}
Sign tensors are spanned by the images of
\[
 \Alt_{[n]}(\theta_J f),
 \qquad
 f\text{ symmetric on }J
       \text{ and alternating on }J^c,
\]
where $\theta_J$ is an exterior monomial with support $J$.
Indeed, average first over the stabilizer
$S_J\times S_{J^c}$: the sign of the exterior monomial
cancels the sign projector on $S_J$.

Since $q_a$ is invariant, its product with such a tensor
is a sum over $i\in J^c$ of terms
\[
 \Alt_{[n]}(\theta_i\theta_J\,a^{(i)}f).
\]
Fix one term and put $D=J^c\setminus\{i\}$.
The exterior factor has support $D^c=J\cup\{i\}$.
Inside the global sign projector, its coefficient can
be replaced by its symmetric average over $D^c$:
the exterior sign again cancels the sign projector.
But $f$ is alternating on $D\cup\{i\}$ and symmetric
on its complement.  Lemma~\ref{lem:complement-average}
therefore says that the averaged coefficient
$\Av_{D^c}(a^{(i)}f)$ is zero in $R_A$.
Every term vanishes.
\end{proof}

\begin{proof}[Proof of Theorem~\ref{thm:augmented-transfer}]
First quotient $S$ by the even generators $t_a$,
$a\in\mathfrak m$.  The result is
$R_A\otimes\Lambda N_n$.
Next quotient by $q_1=\theta_\Sigma$; the result is
$R_A\otimes\Lambda W_n$.
Since $A=K1\oplus\mathfrak m$, the only generators
of $J_A$ still to impose are $q_a$ with
$a\in\mathfrak m$.

Lemma~\ref{lem:odd-annihilation} kills these operators
on the sign part before imposing $q_1$.  Exactness of
the sign projector shows that they also kill the sign
part after imposing $q_1$.
Moreover multiplication by each $q_a$ is
$S_n$-equivariant.  The sign part of its image is
therefore its image on the sign part, and is zero.
The remaining ideal has no sign component.
Quotienting by that ideal preserves the sign submodule
and proves \eqref{eq:augmented-transfer}.
All the constructions used the canonical quotient
maps and averaging projectors, which also proves
the asserted naturality and grading compatibility.

For a homomorphism $K\to K'$, put $A'=A\otimes_KK'$.
The averaging and sign projectors are split idempotents, so their
images commute with arbitrary scalar extension.  The augmentation
on each invariant algebra is split by the unit; hence its kernel
also commutes with scalar extension.  Present $I_A$ and $J_A$ as
images of the multiplication maps from these kernels.  Right
exactness of tensor product then gives canonical isomorphisms
\[
 R_A\otimes_KK'\cong R_{A'},\qquad
 Q_A\otimes_KK'\cong Q_{A'}.
\]
The splitting $N_n=K\theta_\Sigma\oplus W_n$ and its finite-free
exterior powers commute with scalar extension as well.  Thus so
does \eqref{eq:augmented-transfer}.
\end{proof}

\begin{corollary}[Hook extraction on a pointed algebra]
\label{cor:pointed-hooks}
For $0\leq d\leq n-1$ there are natural isomorphisms
\[
 \Hom_{K[S_n]}(\sgn,Q_A[d])
       \cong
 \Hom_{K[S_n]}(S^{(d+1,1^{n-d-1})},R_A),
\]
where $[d]$ denotes exterior degree and $S^\mu$ is the base
change to $K$ of the integral Specht module.
They commute with arbitrary base change and every augmented
automorphism of $A$.
\end{corollary}

\begin{proof}
Theorem~\ref{thm:augmented-transfer} identifies the
left side with
$\Hom_{K[S_n]}(\sgn,R_A\otimes\Lambda^dW_n)$.
The standard pairing on $N_n$ restricts to a perfect pairing on
$W_n$: the complementary vector $\theta_\Sigma$ has norm $n$,
a unit.  Since $W_n$ is finite free and $\det W_n=\sgn$,
exterior duality and the integral hook Specht model
\cite[Section~5.5(b)]{DanzHofmann} give
\[
 \sgn\otimes\Lambda^dW_n
       \cong S^{(d+1,1^{n-d-1})}.
\]
Finite-free tensor-Hom adjunction gives the result.  Equivariant
Hom is the image of an averaging projector on ordinary Hom, so
these identifications also commute with arbitrary base change.
\end{proof}

\begin{example}
Let $A=K[z_1,\ldots,z_k]/L$, where every element of
$L$ vanishes at the origin.  Evaluation at the origin
defines an augmentation.
Theorem~\ref{thm:augmented-transfer} applies whether or
not the pointed affine scheme $\operatorname{Spec}A$
is smooth or reduced.  It applies, for example, to
$K[u,v]/(uv)$ and to $K[u_1,\ldots,u_k]/(u_1,\ldots,u_k)^m$.
The theorem also applies over $\Z[1/n!]$ and to every square-zero
extension $A=K\oplus M$, with $M^2=0$, for an arbitrary $K$-module
$M$.  In particular, $M$ need not be flat.  Over a field of
characteristic $p>n$, the same theorem holds.
The characteristic-zero Schur character identities in
the following subsections use the polynomial case and
additional representation-theoretic input.
\end{example}

\subsection{Supersymmetric specialization}
\label{sec:hook-identity}

Return to a field $K$ of characteristic zero and take
$A=K[z_1,\ldots,z_k]$ with augmentation at the origin.
Then $R_A=R_{n,k}$.
The algebra $Q_A$ has one exterior permutation alphabet in
addition to these bosonic alphabets.
For a row alphabet $\bq=(q_1,\ldots,q_k)$ and one odd
variable $u$, use the positive supersymmetric convention
\begin{equation}\label{eq:superschur}
 s_\lambda(\bq\mid u)
   =\sum_{\lambda/\rho\text{ a vertical strip}}
           s_\rho(\bq)\,u^{|\lambda|-|\rho|}.
\end{equation}
In particular,
$[u^d]s_\lambda(\bq\mid u)
     =(e_d^\perp s_\lambda)(\bq)$.

Lentfer's diagonal supersymmetry theorem
\cite[Theorem~1.2 and Corollary~1.3]{Lentfer}
states that the same universal coefficients
$c_{\lambda,\mu}^{(n)}$ describe the full coinvariant
characters in arbitrary bosonic and fermionic alphabets.
For the sign component with one odd alphabet, it says
\begin{equation}\label{eq:lentfer-sign}
 \ch\Hom_{S_n}(\sgn,Q_A)
       =\sum_\lambda
          c_{\lambda,(1^n)}^{(n)}s_\lambda(\bq\mid u).
\end{equation}

\begin{theorem}[All full hook identities]\label{thm:full-hook}
For $0\leq d\leq n-1$,
\begin{equation}\label{eq:full-hook}
 C_{(d+1,1^{n-d-1})}^{(n)}=e_d^\perp A_n.
\end{equation}
Equivalently,
\[
 c_{\rho,(d+1,1^{n-d-1})}^{(n)}
   =\sum_{\lambda/\rho\text{ a vertical }d\text{-strip}}
          c_{\lambda,(1^n)}^{(n)}.
\]
\end{theorem}

\begin{proof}
By Corollary~\ref{cor:pointed-hooks}, the coefficient
of $u^d$ on the left side of \eqref{eq:lentfer-sign}
is $C_{(d+1,1^{n-d-1})}^{(n)}(\bq)$.
By \eqref{eq:superschur}, the coefficient on the
right side is $(e_d^\perp A_n)(\bq)$.
The equality holds for every $k$, and hence holds
as an identity of stable symmetric functions.
\end{proof}

The full hook identity follows from augmented sign transfer and
diagonal supersymmetry.  The two-bosonic-alphabet transfer in
\cite{GLM} gives its two-variable specialization; the augmented
construction applies in every row dimension.  We next construct the
corresponding identity inside the operator module.

\begin{remark}[Order of skewing and specialization]
\label{rem:skew-specialization}
In \(k\) alphabets the formula is
\[
 C_{(d+1,1^{n-d-1})}^{(n)}(\bq)=(e_d^\perp A_n)(\bq).
\]
It is generally impossible to replace $A_n$ by the
sum of its terms of length at most $k$ before skewing.
A vertical strip can remove the last box in each of
$d$ rows, so source partitions of length $k+d$ can
contribute to a target of length $k$.
\end{remark}

\subsection{Hook transfer in the polarization module}
\label{sec:r2-operator-hooks}

Invariant differential forms give the operator counterpart.  An
odd differential trace can be realized by a polarization into a
fresh row.  Repeating this operation and alternating the fresh
rows lifts every hook vector to an alternating vector, which is
the required inverse to exteriorization.

We work over a field of characteristic zero.  For a finite row space
$V$, write $P(V)=\Sym(V\otimes N_n)$, where $N_n=K^n$ is the permutation
representation.  We use the polynomial row convention for the
variables.  All the subspaces below are polynomial row submodules
compatible with enlarging the row space, with taking homogeneous row
components, and with restriction obtained by setting auxiliary rows
to zero.  Let $\mathcal M(V)$
be a subspace of the diagonal harmonics, stable under constant
coefficient partial derivatives and generalized polarizations.
Its stable alternating character is denoted by $A^{\mathcal M}$.
For $\Harm$ these compatibilities follow directly from the
equations.  For $\Pol$ they follow from the commuting normal form
and its specialization statement,
Corollary~\ref{cor:r3-exact-specialization}.  

For a nonzero multiindex $\alpha$ on the old rows and a fresh row $z$,
the mixed polarization
\begin{equation}\label{eq:r2-mixed-polarization}
 E_{z,\alpha}=\sum_{i=1}^n z_i\partial_i^\alpha
\end{equation}
is an allowed operation.  Indeed, it is a linear combination of
conjugates of $\sum_i z_i\partial_{x_i}^{|\alpha|}$ under changes of
the old rows that fix $z$: pure powers of linear forms span the
homogeneous differential symbols of degree $|\alpha|$.

For $0\leq j\leq n-1$, add $j$ auxiliary rows $z_1,\ldots,z_j$.
Let $\mathscr E_j^{\mathcal M}(V)$ be the space of elements of
$\mathcal M(V\oplus K^j)^{\sgn}$ that are multilinear in these rows
and alternating under their permutations.  Such a polynomial has a
unique expression
\[
 F(X,z)=\sum_{i_1,\ldots,i_j}
       f_{i_1,\ldots,i_j}(X)z_{1,i_1}\cdots z_{j,i_j},
\]
with alternating coefficient tensor.  Set
\begin{equation}\label{eq:r2-exteriorization-map}
 \Psi_j(F)=\sum_{i_1,\ldots,i_j}
           f_{i_1,\ldots,i_j}(X)
                 \theta_{i_1}\cdots\theta_{i_j}.
\end{equation}
Here the $\theta_i$ generate $\Lambda N_n$.  Translation invariance in
each auxiliary row puts the coefficient tensor in $W_n^{\otimes j}$,
where $W_n=\ker(N_n\to K)$; hence the image is in $\Lambda^jW_n$.
Taking coefficients in the auxiliary rows is a sequence of partial
derivatives followed by setting those rows to zero, so every
$f_{i_1,\ldots,i_j}$ belongs to $\mathcal M(V)$.
Consequently
\begin{equation}\label{eq:r2-exteriorization-injection}
 \Psi_j:\mathscr E_j^{\mathcal M}(V)
       \lhook\joinrel\longrightarrow
          (\mathcal M(V)\otimes\Lambda^jW_n)^{\sgn}
\end{equation}
is a $\GL(V)$-equivariant injection.  Injectivity follows directly:
exteriorization is multiplication by $j!$ on alternating tensors.
The case $j=0$ is the identity.

Branching the stable alternating character from $V\oplus K^j$ to
$\GL(V)\times\GL_j$, and selecting the determinant representation of
$\GL_j$, gives
\begin{equation}\label{eq:r2-exteriorization-character}
 \ch\mathscr E_j^{\mathcal M}(V)
          =(e_j^\perp A^{\mathcal M})(V).
\end{equation}
Indeed, a multilinear alternating function of the auxiliary rows
transforms by their determinant, and the relevant
Littlewood--Richardson coefficients are
$c_{\rho,(1^j)}^\lambda$.

\begin{theorem}[Operator hook transfer]\label{thm:r2-operator-hook-transfer}
Let $\mathcal M$ be a compatible system of harmonic submodules generated
from alternating polynomials by constant coefficient derivatives and
generalized polarizations.  Then \eqref{eq:r2-exteriorization-injection}
is an isomorphism for every finite row space $V$ and every
$0\leq j\leq n-1$.  In particular,
\begin{equation}\label{eq:r2-operator-hook-character}
 C_{(j+1,1^{n-j-1})}^{\mathcal M}
                  =e_j^\perp A^{\mathcal M}.
\end{equation}
For $\mathcal M=\Pol$, this is Bergeron's Hook-Components Conjecture.
\end{theorem}

\begin{proof}
Put $\mathcal D=K[\partial_{x_i^{(a)}}]$.  By
Lemma~\ref{lem:r2-normal-order}, the evaluation map
\[
 \mathcal D\otimes\mathcal B(V)\longrightarrow\mathcal M(V),
 \qquad D\otimes g\longmapsto Dg,
\]
is surjective, and $S_n$ acts by its sign on $\mathcal B(V)$.
Exactness of averaging therefore shows that every element of the
right side of \eqref{eq:r2-exteriorization-injection} is a sum of
invariant differential forms in
$(\mathcal D\otimes\Lambda^jW_n)^{S_n}$ acting on elements of
$\mathcal B(V)$.

The invariant algebra $(\mathcal D\otimes\Lambda N_n)^{S_n}$ is generated by
the even traces and odd traces
\[
 p_\alpha=\sum_i\partial_i^\alpha
       \quad(|\alpha|>0),\qquad
 q_\alpha=\sum_i\theta_i\partial_i^\alpha
       \quad(|\alpha|\geq0).
\]
This is Lemma~\ref{lem:trace-generation}, applied to the polynomial
algebra of differential symbols.  Passing to $\Lambda W_n$ sets
$q_0=\sum_i\theta_i$ to zero.  Averaging is exact, so the images
of these generators still generate all invariants in the quotient.
The $p_\alpha$ commute with all the $q_\beta$ and annihilate the
harmonic polynomial on which the form acts.  A form of exterior
degree $j$ consequently reduces to a linear combination of
\begin{equation}\label{eq:r2-hook-odd-word}
             q_{\alpha_1}\cdots q_{\alpha_j}g,
       \qquad |\alpha_s|>0,\quad g\in\mathcal B(V).
\end{equation}

For this word, introduce the $j$ fresh rows and form
\[
 \widetilde F
       =E_{z_1,\alpha_1}\cdots E_{z_j,\alpha_j}g.
\]
All derivatives in these operators use only the old rows.
Thus the operators commute, their product is multilinear in the
fresh rows, and
\[
 \widetilde F
   =\sum_{i_1,\ldots,i_j}
        z_{1,i_1}\cdots z_{j,i_j}
        \partial_{i_1}^{\alpha_1}\cdots
        \partial_{i_j}^{\alpha_j}g.
\]
It belongs to $\mathcal M(V\oplus K^j)^{\sgn}$ by
\eqref{eq:r2-mixed-polarization}.  Let $F$ be its normalized
alternation under permutations of the fresh rows.
Then $F\in\mathscr E_j^{\mathcal M}(V)$.
Exteriorization already alternates these rows, so
\[
                  \Psi_j(F)
                       =q_{\alpha_1}\cdots q_{\alpha_j}g.
\]
Every expression \eqref{eq:r2-hook-odd-word} is therefore in the image
of $\Psi_j$, proving surjectivity.  The character identity follows from
\eqref{eq:r2-exteriorization-character} and
$\sgn\otimes\Lambda^jW_n\cong S^{(j+1,1^{n-j-1})}$.
\end{proof}

The construction also fixes the order of operations in
\eqref{eq:r2-operator-hook-character}: add the $j$ auxiliary rows,
select their exterior component, and then specialize the remaining
rows to $V$.  This is stable skewing followed by finite-row
specialization, as in Remark~\ref{rem:skew-specialization}.

\subsection{Simultaneous derivative embeddings}

Taylor expansion gives simultaneous maps behind the character
identities.  Derivative closure keeps every Taylor coefficient
inside the module.  Schur--Weyl duality then separates all
Littlewood--Richardson copies, so their contributions form an
actual direct sum.  This construction applies to any
derivative-closed polynomial module.

\begin{theorem}[The derivative coalgebra and Littlewood--Richardson embeddings]
\label{thm:derivative-coalgebra}
Let $V,W$ be finite-dimensional vector spaces over a field of
characteristic zero, let $G\subseteq\GL(W)$, and let
$\mathcal M\subseteq\Sym(V\otimes W)$ be a graded
$\GL(V)\times G$-submodule closed under all constant-coefficient
partial derivatives. Put
\[
 M_\lambda=\Hom_{\GL(V)}(\Schur_\lambda V,\mathcal M),
 \qquad M_\lambda=0\quad\text{if }\ell(\lambda)>\dim V.
\]
Then $\mathcal M$ is a subcoalgebra for the polynomial coproduct
$\Delta F(X,Y)=F(X+Y)$. For every pair of partitions $\tau,\rho$,
the coproduct induces a $G$-equivariant injection
\begin{equation}\label{eq:lr-coalgebra-injection}
 \bigoplus_\lambda M_\lambda^{\oplus c_{\tau,\rho}^{\lambda}}
       \lhook\joinrel\longrightarrow M_\tau\otimes M_\rho.
\end{equation}
All Littlewood--Richardson multiplicities occur independently.
The maps are compatible with inclusions of derivative-closed submodules.
\end{theorem}
\begin{proof}
Taylor expansion and derivative closure give
\[
 \Delta\mathcal M\subseteq
 [\Sym(V\otimes W)\otimes\mathcal M]
 \cap[\mathcal M\otimes\Sym(V\otimes W)]
 =\mathcal M\otimes\mathcal M.
\]
The counit is evaluation at zero. This proves the coalgebra assertion.

Put $a=|\tau|$ and $b=|\rho|$. Project the bidegree-$(a,b)$
coproduct to the row types $\tau$ and $\rho$ in its two factors.
The ambient Cauchy decomposition has column multiplicity
$\Schur_\lambda W$ in row type $\lambda$. Schur--Weyl duality
constructs this summand by the invariant pairing between two copies
of the Specht module $[\lambda]$.  Write
$L_{\tau,\rho}^{\lambda}$ for the multiplicity space of
$[\tau]\boxtimes[\rho]$ in its restriction to $S_a\times S_b$;
its dimension is $c_{\tau,\rho}^{\lambda}$.
The bidegree coproduct uses the restricted Specht pairing.
Over the reals, choose this pairing positive and invariant, so
its restriction to each isotypic component is nondegenerate.
On the two copies of $L_{\tau,\rho}^{\lambda}$ it is therefore a
perfect pairing.  In dual bases its tensor is the identity:
each row coefficient selects a distinct column coefficient.
We obtain exactly $c_{\tau,\rho}^{\lambda}$ independent ambient maps
\[
 \Schur_\lambda W\longrightarrow
                  \Schur_\tau W\otimes\Schur_\rho W.
\]
For different $\lambda$ their images lie in distinct isotypic
summands; for a fixed $\lambda$ they fill its independent
Littlewood--Richardson copies. Their direct sum is injective.
All relevant pairings and projectors can be chosen over $\Q$,
so the assertion holds over every characteristic-zero field.

The coalgebra inclusion places the projected coproduct inside
the row spaces tensored with $M_\tau\otimes M_\rho$. Taking its
independent row coefficients therefore restricts the ambient maps
to \eqref{eq:lr-coalgebra-injection}. The column group is used
only on the ambient polynomial ring; $\mathcal M$ need only be
$G$-stable. Fixed choices of ambient maps give compatibility
with inclusions.
\end{proof}

\begin{corollary}[Exterior derivative embeddings]
\label{thm:exterior-derivatives}
Under the hypotheses of Theorem~\ref{thm:derivative-coalgebra},
for every partition $\rho$ and $j\geq0$ there is a
$G$-equivariant injection
\begin{equation}\label{eq:exterior-injection}
 \bigoplus_{\lambda/\rho\text{ a vertical }j\text{-strip}}
       M_\lambda
 \lhook\joinrel\longrightarrow\Lambda^jW\otimes M_\rho.
\end{equation}
The maps arise from the exterior row component of the $j$th
derivatives and are compatible with inclusions.
\end{corollary}
\begin{proof}
Take $\tau=(1^j)$ in \eqref{eq:lr-coalgebra-injection} and use
the vertical Pieri rule and $M_{(1^j)}\subseteq\Lambda^jW$.
\end{proof}

\begin{corollary}[Derivative realization and exact hook defects]
\label{cor:hook-lower}
\label{cor:hook-realization}
\label{cor:hook-defects}
For $\mathcal M=\Harm$ or $\Pol$, let $M_\rho^{\mathcal M}$ denote
its stable row multiplicity space and put
$\mathcal A_\lambda^{\mathcal M}
=\Hom_{S_n}(\sgn,M_\lambda^{\mathcal M})$.
The exterior derivative maps give compatible isomorphisms
\begin{equation}\label{eq:hook-realization}
 \bigoplus_{\lambda/\rho\text{ a vertical }d\text{-strip}}
       \mathcal A_\lambda^{\mathcal M}
 \xrightarrow{\ \sim\ }
 \Hom_{S_n}(S^{(d+1,1^{n-d-1})},M_\rho^{\mathcal M}),
 \qquad 0\le d<n,
\end{equation}
and induce isomorphisms on the quotients for $\Pol\subseteq\Harm$.
Thus, writing $D_A=A_n^\Harm-A_n^\Pol$ and
$D_d=C_{(d+1,1^{n-d-1})}^\Harm-C_{(d+1,1^{n-d-1})}^\Pol$,
\begin{align}
 C_{(d+1,1^{n-d-1})}^{\mathcal M}
       &=e_d^\perp A_n^{\mathcal M},\label{eq:hook-lower}\\
 D_d&=e_d^\perp D_A.\label{eq:hook-defects}
\end{align}
Restricting source and target to $k$ rows gives
\[
 c_{\rho,(d+1,1^{n-d-1})}^{\mathcal M}
 \ge\sum_{\substack{\lambda/\rho\text{ vertical }d\text{-strip}\\
                    \ell(\lambda)\le k}}
            c_{\lambda,(1^n)}^{\mathcal M}.
\]
\end{corollary}
\begin{proof}
Apply Corollary~\ref{thm:exterior-derivatives} with $W=W_n$
and take $\Hom_{S_n}(\sgn,-)$.
Since $\sgn\otimes\Lambda^dW_n\cong S^{(d+1,1^{n-d-1})}$,
the target is the indicated hook multiplicity.
Theorems~\ref{thm:full-hook} and~\ref{thm:r2-operator-hook-transfer}
give equality of dimensions, so these compatible injections are
isomorphisms. They induce the quotient isomorphisms and both
character identities. Retaining only source partitions of length
at most $k$ gives the finite-row bound.
\end{proof}

Thus full alternating generation is equivalent to simultaneous
generation of every hook label. The derivative maps retain the
individual vertical-strip sources, and their quotients identify
every remaining hook defect.

\subsection{Two-column support and the alternating defect}
\label{sec:reconstruction}

After removing the linear moments, the labels form the standard
representation $W_n$ with its invariant inner product.  The
quadratic moments are its pairwise Gram products.  Two alternating
blocks whose total height exceeds $\dim W_n$ vanish modulo those
products.  Applied to the first two columns of a row partition,
this gives the support restriction used in the uniform generation proof.

\begin{lemma}[Two alternating blocks in the Gram ideal]
\label{lem:gram}
Let $E$ be a vector space of dimension $N$ over a
field of characteristic zero, equipped with a
nondegenerate symmetric bilinear form.
A multilinear polynomial on $E^{d+e}$ that is
alternating in the first $d$ arguments and in the
last $e$ arguments belongs to the ideal generated
by the inner products between the two blocks if $d+e>N$.
\end{lemma}

\begin{proof}
Identify $E$ with $E^*$ by the form. On the tensor product of exterior
algebras, using ordinary tensor signs, put
\[
 L=\sum_i\eps_i\otimes\eps_i,\qquad
 c=\sum_i\iota_i\otimes\iota_i.
\]
Invariantly the form and its inverse define these operators over $K$.
They satisfy $[c,L]=(N-d-e)\operatorname{id}$ in bidegree $(d,e)$.
Finite-dimensional $\mathfrak{sl}_2$ theory makes
\[
 L:\Lambda^{d-1}E\otimes\Lambda^{e-1}E
                 \longrightarrow\Lambda^dE\otimes\Lambda^eE
\]
surjective when $d+e>N$. Under the multilinear polynomial interpretation,
its image consists of alternating sums of a cross-block inner product
multiplied by a polynomial in the remaining arguments. This proves
membership in the stated ideal. An explicit finite right inverse is
recorded in \eqref{eq:gram-finite-certificate}.
\end{proof}

The next inequality is the admissibility condition for orthogonal
traceless tensors; see, for example, \cite{BGH}.  The Gram-ideal
argument proves it directly at the dimension $n-1$ relevant here.

\begin{theorem}[Two-column support]\label{thm:two-column}
For every coefficient of full diagonal harmonics,
\begin{equation}\label{eq:two-column}
 c_{\lambda,\mu}^{(n)}\ne0
       \quad\Longrightarrow\quad
       \lambda'_1+\lambda'_2\leq n-1.
\end{equation}
Here missing columns have length zero.
The same restriction holds in $\Pol$.
\end{theorem}

\begin{proof}
Quotient by the linear invariants first.  This
identifies the polynomial ring with
$\Sym(V\otimes W_n)$, where $V=K^k$ and
$\dim W_n=n-1$.
The quadratic invariants include all pairwise
inner products of the row vectors for the
nondegenerate form on $W_n$.
Thus the full coinvariant ring is a quotient of
the quotient by the Gram ideal.

The usual minor construction of the Cauchy
decomposition spans a row-$\lambda$ highest-weight
space by products of minors with column heights
$\lambda'_1,\lambda'_2,\ldots$.
If the sum of the first two heights exceeds
$n-1$, regard the vectors in these two minors
temporarily as independent.  Their product is
multilinear and alternating separately in the
two blocks.  Lemma~\ref{lem:gram} puts this product
in the Gram ideal.  Identifying the repeated
row vectors and multiplying by the remaining
minors keeps it in that ideal.
If only one column is present, a height greater
than $n-1$ is already excluded by the Cauchy
decomposition.
Every such highest-weight vector vanishes in
the quotient, proving \eqref{eq:two-column}.
The assertion for $\Pol$ follows by inclusion.
\end{proof}

\begin{corollary}[Support of the alternating defect]
\label{cor:r2-defect-support}
If the Schur coefficient of $s_\lambda$ in $D_A$ is nonzero, then
\begin{equation}\label{eq:r2-defect-support}
 \lambda_2\geq3,\qquad \lambda_3\geq2,\qquad
                 \lambda'_1+\lambda'_2\leq n-1.
\end{equation}
In particular, $\Pol$ contains every stable alternating component and
every stable hook-label component when $n\leq6$.
\end{corollary}

\begin{proof}
Haiman's two-row operator theorem gives $D_j(q_1,q_2)=0$ for every $j$.
If $\lambda_3\leq1$, delete one box from each nonempty row after the
second; the remaining partition has at most two rows, and the removed
boxes form a vertical strip.  Schur positivity of $D_A$ and the
vertical Pieri rule would then give a nonzero coefficient in some
$e_j^\perp D_A(q_1,q_2)$, a contradiction.

Every row-hook component of every label is generated, by
Proposition~\ref{prop:row-hook-generation}.  If $\lambda_2\leq2$, remove
the last box from every row of length two below the first.  The
remaining partition is a row hook.  Again positivity and
\eqref{eq:hook-defects} contradict generation of that
row-hook component.  The strip sizes used here are at most $n-1$,
since the row-length bound applies to the alternating label.

The final inequality is the two-column support theorem.  The first
two inequalities force the first two columns of $\lambda$ to have
at least three boxes each.  Thus $n\geq7$ for a nonzero defect.
Equation~\eqref{eq:hook-defects} gives the final assertion.
\end{proof}

\subsection{Finite-alphabet reconstruction}

The hook identities also recover stable alternating information
beyond the number of rows in which the hook characters are known.

\begin{theorem}[Reconstruction from finitely many alphabets]
\label{thm:reconstruction}
Let $\mathcal M=\Harm$ or $\Pol$, let $n\geq3$, and put
$k=\lfloor(n-1)/2\rfloor$.
The collection
\begin{equation}\label{eq:reconstruction-data}
 \left\{
 C_{(d+1,1^{n-d-1})}^{\mathcal M}(q_1,\ldots,q_k):
          0\leq d\leq n-1
 \right\}
\end{equation}
determines every Schur coefficient of the stable
alternating character $A_n^{\mathcal M}$.
\end{theorem}

\begin{proof}
Suppress the superscript $\mathcal M$, write
$A_n=\sum_\lambda a_\lambda s_\lambda$, and form
\[
 F(\bq,u)=\sum_{d=0}^{n-1}
 C_{(d+1,1^{n-d-1})}(\bq)u^d.
\]
The hook identity gives
\begin{equation}\label{eq:reconstruction-F}
 F(\bq,u)=\sum_\lambda a_\lambda
                         s_\lambda(\bq\mid u).
\end{equation}
If $\lambda_{k+1}\geq2$, then
$\lambda'_1+\lambda'_2\geq2k+2>n-1$.
Theorem~\ref{thm:two-column} excludes such a
partition.  Every supported partition therefore
satisfies $\lambda_{k+1}\leq1$.

At $u=0$, equation \eqref{eq:reconstruction-F}
determines all $a_\lambda$ with $\ell(\lambda)\leq k$,
because their Schur functions in $k$ variables
are linearly independent.  Subtract their entire
supersymmetric contributions from $F$.
Every remaining partition has the form
\[
 \lambda=(\lambda_1,\ldots,\lambda_k,1^r),
       \qquad r>0,\quad\lambda_k\geq1.
\]
The standard factorization of a supersymmetric
Schur function on the $(k,1)$ hook gives
\begin{equation}\label{eq:super-factorization}
 s_\lambda(\bq\mid u)
    =\prod_{i=1}^k(q_i+u)\,
       s_{(\lambda_1-1,\ldots,\lambda_k-1)}(\bq)\,u^r.
\end{equation}
It also follows directly from the vertical-strip
expansion \eqref{eq:superschur} and the Pieri rule.
Divide the residual polynomial by
$\prod_i(q_i+u)$.
The polynomials $s_\nu(\bq)u^r$ are linearly
independent, and the pair $(\nu,r)$ recovers
$\lambda=(\nu_1+1,\ldots,\nu_k+1,1^r)$.
This determines all remaining coefficients.
\end{proof}

Bergeron's reconstruction from hook information
\cite[Section~3.1]{Bergeron} is formulated under the corresponding
structural conjectures.  Theorem~\ref{thm:reconstruction} establishes
it for both the full and operator characters, with the number of
rows fixed by the two-column bound.  For seven labels, three-row
hook characters determine the entire stable alternating character.

\section{The sign component of the Theta conjecture at \texorpdfstring{$t=0$}{t=0}}
\label{sec:theta-sign}

We prove the sign-isotypic component of the $t=0$ specialization of
the Theta coinvariant conjecture of D'Adderio, Iraci and Vanden
Wyngaerd \cite[Conjecture~8.2]{DIV2021}, in every rank.  In particular,
the alternating coinvariants with one bosonic and two fermionic
alphabets have dimension $3^{n-1}$, proving the numerical prediction
in Bergeron's Table~3 \cite{BergeronSuper2020}.
Iraci, Nadeau and Vanden Wyngaerd
\cite[Theorem~5.6]{INV2024} computed the symmetric-function prediction;
the argument here identifies that prediction with the actual
coinvariant sign component.  Its ingredients are the full hook
identity, diagonal supersymmetry, and the Fields character theorem.

\subsection{Adjoining a fermionic alphabet}

For $k,j\geq0$, let
\begin{align*}
 P_n^{(k,j)}
   &=K[x_i^{(a)}:1\leq i\leq n,\ 1\leq a\leq k]
       \otimes\Lambda(\theta_i^{(b)}:1\leq i\leq n,\ 1\leq b\leq j),\\
 R_n^{(k,j)}
   &=P_n^{(k,j)}/\langle(P_n^{(k,j)})^{S_n}_+\rangle.
\end{align*}
Here $K$ has characteristic zero, all odd generators anticommute,
and $S_n$ acts by simultaneously permuting their label indices.
The subscript $+$ denotes zero constant term.  We record each
bosonic and fermionic alphabet separately, by
$Q=(q_1,\ldots,q_k)$ and $U=(u_1,\ldots,u_j)$, respectively.
Supersymmetric evaluation $f[Q\mid U]$ uses the positive convention
of \eqref{eq:superschur}; in particular, odd degrees have positive
coefficients in the Hilbert series.

Write $W_n=\ker(K^n\xrightarrow{\sum}K)$.  If
$\eta_1,\ldots,\eta_n$ is a new fermionic alphabet, its invariant
$\sum_i\eta_i$ is zero in $R_n^{(k,j+1)}$.  The quotient
$K^n/K\mathbf1$ identifies canonically with $W_n$ by
$z\mapsto z-n^{-1}(\sum_i z_i)\mathbf1$.
Consequently the quotient map factors through
$R_n^{(k,j)}\otimes\Lambda W_n$; the tensor product uses the Koszul
sign rule with the old odd variables.

\begin{theorem}[Natural sign transfer in superspace]
\label{thm:all-superspace-sign-transfer}
For all $n\geq1$ and $k,j\geq0$, the natural quotient induces an
isomorphism
\begin{equation}\label{eq:all-superspace-sign-transfer}
 \left(R_n^{(k,j)}\otimes\Lambda W_n\right)^{\sgn}
 \xrightarrow{\ \sim\ }
 \left(R_n^{(k,j+1)}\right)^{\sgn}.
\end{equation}
It preserves every old grading and the degree in the new fermionic
alphabet, and commutes with parity-preserving changes of the old
row alphabets.  In new fermionic degree $0\leq d\leq n-1$, the
induced isomorphism of multiplicity spaces is
\begin{equation}\label{eq:super-hook-multiplicity}
 \Hom_{S_n}\!\left(S^{(d+1,1^{n-d-1})},R_n^{(k,j)}\right)
 \cong
 \Hom_{S_n}\!\left(\sgn,
       (R_n^{(k,j+1)})_{\deg_\eta=d}\right).
\end{equation}
Both sides of \eqref{eq:all-superspace-sign-transfer} vanish in
new fermionic degree greater than $n-1$.
\end{theorem}

The case $(k,j)=(2,0)$ is the sign-transfer theorem of
Gonz\'alez, Lentfer and Mularczyk \cite[Theorem~4.21]{GLM}.
Here the full hook identity permits arbitrary old bosonic and
fermionic alphabets.

\begin{proof}
The quotient is surjective on sign components, since the sign
idempotent is exact in characteristic zero.  By Lentfer's diagonal
supersymmetry theorem \cite[Theorem~1.2 and Corollary~1.3]{Lentfer},
the character of the $S^\mu$ multiplicity space of $R_n^{(k,j)}$ is
$C_\mu^{(n)}[Q\mid U]$.  Theorem~\ref{thm:full-hook} gives
\[
 C_{(d+1,1^{n-d-1})}^{(n)}=e_d^\perp A_n
 \qquad(0\leq d\leq n-1).
\]
Super Schur branching in one additional odd variable $v$ yields
\begin{equation}\label{eq:super-sign-branching}
 A_n[Q\mid U,v]
   =\sum_{d\geq0}v^d(e_d^\perp A_n)[Q\mid U].
\end{equation}
The skew operator acts on the stable symmetric function before
the finite supersymmetric specialization is made.

The invariant dot product on $W_n$ identifies its exterior powers
with their duals, and
$\sgn\otimes\Lambda^dW_n\cong S^{(d+1,1^{n-d-1})}$.
Tensor--Hom adjunction therefore identifies the coefficient of
$v^d$ in the source sign multiplicity with
$C_{(d+1,1^{n-d-1})}^{(n)}[Q\mid U]$.
This equals the coefficient in \eqref{eq:super-sign-branching},
which is the target sign multiplicity.  Every multigraded piece is
finite-dimensional, so the surjection is an isomorphism in every
multidegree, proving both displayed isomorphisms.  The source has
no exterior degree above $n-1$; the target also has none, because
the relation $\sum_i\eta_i=0$ reduces its new exterior alphabet to
$W_n$.  Naturality follows from the quotient map itself.
\end{proof}

We use the usual grading convention when passing to harmonic
models: a harmonic polynomial and its coinvariant class have the
same multidegree.  More explicitly, fix an order on the odd
generators and use the apolar pairing
\[
 \langle x^\alpha\theta_I,x^\beta\theta_J\rangle
   =\alpha!\,\delta_{\alpha,\beta}\delta_{I,J}.
\]
The odd differential operators are taken in the reverse order to
the corresponding monomial.  The harmonic inverse system is the
orthogonal complement of the invariant ideal for this pairing.
The pairing is positive definite over the reals and is defined
over $\Q$, so the quotient projection identifies the harmonic
model and the coinvariant module in each multidegree over every
characteristic-zero field.  Thus no inversion of $q$, $u$, or $v$
occurs when using the Fields formula \eqref{eq:fields-character}.
On harmonic inverse systems, \eqref{eq:all-superspace-sign-transfer}
also gives
\[
 (\Harm_n^{(k,j+1)})^{\sgn}
   =(\Harm_n^{(k,j)}\otimes\Lambda W_n)^{\sgn}:
\]
the final harmonic space is contained in the intermediate one,
and the quotient isomorphism makes their sign components equal.

\subsection{The actual alternating character}

Put $[m]_q!=\prod_{r=1}^m(1+q+\cdots+q^{r-1})$, with
$[0]_q!=1$, and use the Gaussian coefficients
\[
 \genfrac{[}{]}{0pt}{}{N}{a,b,c}_q
    =\frac{[N]_q!}{[a]_q![b]_q![c]_q!}
 \quad(a+b+c=N),
 \qquad
 \genfrac{[}{]}{0pt}{}{N}{a}_q
    =\frac{[N]_q!}{[a]_q![N-a]_q!}.
\]

\begin{theorem}[The $t=0$ Theta sign component]
\label{thm:one-boson-two-fermion-sign}
For every $n\geq1$,
\begin{equation}\label{eq:one-boson-two-fermion-sign}
 \operatorname{Hilb}\!\left((R_n^{(1,2)})^{\sgn};q;u,v\right)
 =\sum_{\substack{a,b\geq0\\a+b\leq n-1}}
 u^av^b q^{\binom{n-a-b}{2}}
 \genfrac{[}{]}{0pt}{}{n-1}{a,b,n-1-a-b}_q.
\end{equation}
Here $q$ records the bosonic degree, and $u,v$ record the two
fermionic degrees.  This is the sign coefficient of
\cite[Conjecture~8.2]{DIV2021} at $t=0$, in every fermionic
bidegree.\footnote{Equation~(76) in the arXiv version of
\cite{DIV2021} prints $1\leq a+b<n$.  We include the ordinary
term $a=b=0$, whose sign character $q^{\binom n2}$ is classical.}
In particular,
\begin{equation}\label{eq:one-boson-two-fermion-sign-dimension}
 \operatorname{Hilb}\!\left((R_n^{(1,2)})^{\sgn};1;u,v\right)
 =(1+u+v)^{n-1},
 \qquad
 \dim(R_n^{(1,2)})^{\sgn}=3^{n-1}.
\end{equation}
\end{theorem}

\begin{proof}
Fix $0\leq b\leq n-1$ and set $N=n-1$, $d=N-b$.
Theorem~\ref{thm:all-superspace-sign-transfer} identifies the
coefficient of $v^b$ with the multiplicity of the hook
$\mu=(b+1,1^d)$ in $R_n^{(1,1)}$.
The Fields character theorem of Murai, Rhoades and Wilson
\cite[Theorem~5.2 and equation~(2.5)]{MRW}, in the form
\eqref{eq:fields-character}, evaluates this multiplicity as
\[
 \sum_{T\in\SYT(\mu)}q^{\maj(T)}
                  \prod_{r=1}^{\des(T)}(1+u q^{-r}).
\]
The hook-tableau enumeration used in \eqref{eq:hook-hook-formula}
gives $\operatorname{des}(T)=d$ and
$\sum_Tq^{\maj(T)}=q^{\binom{d+1}{2}}
 \genfrac{[}{]}{0pt}{}{N}{d}_q$.
It follows that
\begin{equation}\label{eq:theta-sign-fields-product}
 [v^b]\operatorname{Hilb}\!\left((R_n^{(1,2)})^{\sgn};q;u,v\right)
   =q^{\binom{d+1}{2}}\genfrac{[}{]}{0pt}{}{N}{b}_q
                \prod_{r=1}^{d}(1+u q^{-r}).
\end{equation}
For $0\leq a\leq d$, the Gaussian binomial theorem gives
\[
 [u^a]\prod_{r=1}^{d}(1+u q^{-r})
   =\sum_{\substack{J\subseteq[d]\\|J|=a}}q^{-\sum_{r\in J}r}
   =q^{-ad+\binom a2}\genfrac{[}{]}{0pt}{}{d}{a}_q.
\]
For $a>d$ the coefficient is zero.  Since
\[
 \binom{d+1}{2}-ad+\binom a2=\binom{d-a+1}{2}
      =\binom{n-a-b}{2},
 \qquad
 \genfrac{[}{]}{0pt}{}{N}{b}_q
 \genfrac{[}{]}{0pt}{}{N-b}{a}_q
     =\genfrac{[}{]}{0pt}{}{N}{a,b,N-a-b}_q,
\]
extracting $u^a$ from \eqref{eq:theta-sign-fields-product} proves
\eqref{eq:one-boson-two-fermion-sign}.

For comparison with the stated conjecture, let
$\Theta_f$ and $\nabla$ denote the symmetric-function operators
of \cite{DIV2021}.  The sign coefficient predicted there in
fermionic bidegree $(a,b)$, after setting $t=0$, is
\[
 \left\langle
     \left.\Theta_{e_a}\Theta_{e_b}\nabla e_{n-a-b}\right|_{t=0},
     e_n\right\rangle
   =q^{\binom{n-a-b}{2}}
      \genfrac{[}{]}{0pt}{}{n-1}{a+b}_q
      \genfrac{[}{]}{0pt}{}{a+b}{a}_q
 \qquad(a+b<n),
\]
by \cite[Theorem~5.6]{INV2024}; see also
\cite[Theorem~8.3]{Lentfer12}.  The two Gaussian binomials multiply
to the multinomial in \eqref{eq:one-boson-two-fermion-sign}.
Both series have zero sign coefficient for $a+b\geq n$.
Finally, putting $q=1$ in that equation and applying the
trinomial theorem proves
\eqref{eq:one-boson-two-fermion-sign-dimension}.
\end{proof}

\section{An explicit general linear action and Lefschetz maps}
\label{sec:hook-action}

Let $\mu=(c,\tau)\vdash n$, put $r=|\tau|=n-c>0$, and set
\[
 N=\binom n2,\qquad S=\frac{r(r+1)}2,\qquad L=N-S+1.
\]
Write $K_{a,j}$ for the primitive harmonic-form model of the row hook
$(a,1^j)$ from Lemma~\ref{lem:hook-bridge}, and put
\[
 \mathcal M_\mu=\bigoplus_{a\ge1}
          \operatorname{Hom}_{S_n}(S^\mu,K_{a,r-1}).
\]
Its grading is the first row length $a$. We construct the representation
behind Theorem~\ref{thm:maximal-hook-series}, using the bounded-exponent
model of the ordinary coinvariant algebra.
The spaces $K_{a,r-1}$ are the canonical harmonic-form models of
stable row-hook multiplicity spaces. The auxiliary action is defined on $\mathcal M_\mu$.
In first-row degree $a$, the ordinary $\operatorname{GL}_k$-isotypic
component is $\mathbb S_{(a,1^{r-1})}K^k\otimes\mathcal M_\mu[a]$,
with the row group acting on the first factor.

Let $C_n=K[x_1,\ldots,x_n]/(K[x_1,\ldots,x_n]^{S_n}_+)$.
For $1\le d<n$, let $J_{n,d}\subset C_n$ be generated by the
determinants $\det(x_{i_a}^{\,b})_{1\le a,b\le d}$ over all
$d$-element label sets. Write $\Omega_n$ for the coinvariant algebra
of differential forms, apolar dual to $\SH_n$.
Proposition~\ref{prop:standalone-eulerian-foundations} proves the
bounded-exponent complement and its Foulkes character,
$\operatorname{Ann}_{C_n}(J_{n,d})=J_{n,n-d}$, and the
prefix-contraction image used below.

Let $U_c=K[t]_{<c}$ with ordered basis $1,t,\ldots,t^{c-1}$, and
identify $B_c=U_c^{\otimes n}$ with the bounded-exponent subspace
\[
 B_c=\operatorname{span}\{x_1^{a_1}\cdots x_n^{a_n}:0\le a_i<c\}
                  \subset K[x_1,\ldots,x_n].
\]
The group $\operatorname{GL}(U_c)$ acts on the tensor factors and
commutes with $S_n$. Transporting this action will give the action on $\mathcal M_\mu$.

\begin{lemma}[The bounded-exponent realization at the minimum descent number]
\label{lem:maximal-hook-bounded-model}
The quotient projection induces a graded $S_n$-equivariant isomorphism
\[
 (B_c)_{\mu'}\xrightarrow{\ \sim\ }
                  (C_n/J_{n,c})_{\mu'},
\]
with the quotient and determinant ideal just defined.
\end{lemma}
\begin{proof}
The straightening argument in Appendix~\ref{app:eulerian-foundations}
gives the surjection
$B_c\twoheadrightarrow C_n/J_{n,c}$. Applying the $\mu'$ central
idempotent preserves surjectivity. Schur--Weyl duality identifies the
graded $S^{\mu'}$ multiplicity in the source with
$\mathbb S_{\mu'}U_c$, whose character is
$s_{\mu'}(1,q,\ldots,q^{c-1})$.

Equations~\eqref{eq:standalone-eulerian-complement} and
\eqref{eq:standalone-descent-filtration} identify the target
multiplicity character as
\[
 \sum_{\substack{T\in\operatorname{SYT}(\mu')\\
                              \operatorname{des}(T)<c}}q^{\operatorname{maj}(T)}.
\]
The shape $\mu'$ has exactly $c$ rows, so each tableau has at least
$c-1$ descents. Principal specialization on $c$ variables identifies
the displayed sum, consisting precisely of tableaux with $c-1$
descents, with $s_{\mu'}(1,q,\ldots,q^{c-1})$.
Thus the surjection has equal source and target dimensions in every
degree, proving injectivity as well.
\end{proof}

Use the harmonic-form operators
\[
 \Theta_j=\sum_i\theta_i\partial_{x_i}^j,
 \qquad \iota=\sum_i x_i\partial_{\theta_i},
 \qquad d=\dd=\Theta_1.
\]
For a homogeneous $b\in(B_c)_{\mu'}$ of degree $h$, define
\begin{equation}\label{eq:explicit-maximal-hook-map}
 \Phi_\mu(b)=\frac{1}{N-h-S+r}\,
        \iota\Theta_1\Theta_2\cdots\Theta_r\,b(\partial)\Delta_n.
\end{equation}
The denominator is positive whenever the source component is nonzero.
Indeed,
$\mathbb S_{\mu'}U_c=\det(U_c)\otimes\mathbb S_{\tau'}U_c$,
so its largest polynomial degree is at most
$\binom c2+r(c-1)=N-S$. Consequently
$N-h-S\ge0$ and $N-h-S+r\ge r>0$.

\Needspace{17\baselineskip}
\begin{samepage}
\begin{theorem}[An explicit Schur-module realization]
\label{thm:maximal-hook-action}
The map \eqref{eq:explicit-maximal-hook-map} is a sign-intertwining
isomorphism
\[
 (B_c)_{\mu'}\xrightarrow{\ \sim\ }
                     \bigoplus_{a\ge1}(K_{a,r-1})_\mu,
 \qquad h\longmapsto a=L-h.
\]
On multiplicity spaces, it constructs an auxiliary
$\operatorname{GL}_c$ action for which
\[
                    \mathcal M_\mu\cong\mathbb S_{\tau'}U_c
\]
as representations, with the displayed grading given by the
principal cocharacter
$\rho(z)=\operatorname{diag}(1,z,\ldots,z^{c-1})$:
weight $a-1$ is the harmonic multiplicity space $\mathcal M_\mu[a]$.
Thus the graded character has the shift by $q$ in
Theorem~\ref{thm:maximal-hook-series}; the group acts across first-row degrees.
The action is obtained from the fixed bounded-exponent model and the
displayed differential map; no choice of a graded basis in the harmonic
multiplicity spaces is involved.
\end{theorem}
\end{samepage}
\begin{proof}
The perfect ordinary coinvariant pairing is
\[
 C_n[h]\otimes C_n[N-h]\longrightarrow K\Delta_n.
\]
The complementary annihilator identity
\eqref{eq:standalone-eulerian-annihilator} states
$\operatorname{Ann}_{C_n}(J_{n,r})=J_{n,c}$.
Consequently the pairing restricts to a perfect sign-valued pairing
between $C_n/J_{n,c}$ and $J_{n,r}$, reversing degree around $N$.

On $\Omega_n$, put $\delta_j=\sum_i x_i^j\partial_{\theta_i}$.
The contraction $\delta_r\cdots\delta_1:\Omega_n^r\to C_n$ has
image $J_{n,r}$, and its apolar adjoint is
$\Theta_1\cdots\Theta_r$. On label $\mu$, the largest possible
descent number is $r$. The Fields character~\cite[Theorem~5.2]{MRW} in exterior degree $r$
therefore has only the tableaux with exactly $r$ descents and has
polynomial degrees $\operatorname{maj}(T)-S$.
The character of $(J_{n,r})_\mu$ consists of those same tableaux,
with degrees $\operatorname{maj}(T)$. Thus the prefix contraction
is an isomorphism on that label, with degree shift $S$.
Its adjoint identifies the perfect pairing above with the concrete map
\[
 b\longmapsto \Theta_1\cdots\Theta_r b(\partial)\Delta_n
\]
from $(C_n/J_{n,c})_{\mu'}$ onto $(\mathrm{SH}_n^r)_\mu$, with
the sign twist. Its kernel on $C_n$ is exactly the annihilator
$J_{n,c}$, so the map factors through the stated quotient.

There is no label-$\mu$ component in exterior degree $r+1$.
Hence every form in $(\mathrm{SH}_n^r)_\mu$ is killed by $d$.
The Euler identity $d\iota+\iota d=(p+r)\operatorname{id}$ in
polynomial degree $p$ identifies that space with primitive forms of
exterior degree $r-1$, by $\iota/(p+r)$, with inverse $d$.
Taking $p=N-h-S$ proves the explicit formula and the isomorphism,
using the preceding bounded-exponent lemma.

Apolarity reverses polynomial degree. Reversing the monomial basis
converts this to principal weights, and a determinant twist removes
the full first column of $\mu'$. Explicitly, let
$J:U_c\to U_c$ send $t^j$ to $t^{c-1-j}$.
Transport to the target the following action on the source multiplicity
space $\mathbb S_{\mu'}U_c$:
\[
       g\longmapsto \det(g)^{-1}\mathbb S_{\mu'}(JgJ^{-1}).
\]
Since $\mu'=(1^c)+\tau'$, this is the polynomial irreducible
representation $\mathbb S_{\tau'}U_c$. If
$g=\operatorname{diag}(1,z,\ldots,z^{c-1})$, its eigenvalue on source
degree $h$ is
\[
 z^{n(c-1)-h-\binom c2}=z^{L-h-1}=z^{a-1}.
\]
Here $N-S=\binom c2+r(c-1)$.
This proves the grading shift by one. Normalizing the Specht pairing
or Vandermonde rescales the intertwiner and does not change the
transported action.
\end{proof}

\begin{corollary}[Lefschetz maps on harmonic multiplicities]
\label{cor:maximal-hook-lefschetz}
There are explicitly defined operators $E,F,H$ on $\mathcal M_\mu$
satisfying the $\mathfrak{sl}_2$ relations. They have grading changes
$+1,-1,0$, and
\[
 H|_{\mathcal M_\mu[a]}=(2(a-1)-r(c-1))\operatorname{id}.
\]
In particular, whenever $a\le1+r(c-1)/2$,
\[
 E^{\,2+r(c-1)-2a}:\mathcal M_\mu[a]
       \xrightarrow{\ \sim\ }
                    \mathcal M_\mu[2+r(c-1)-a].
\]
\end{corollary}
\begin{proof}
On $B_c$, use the usual tensor action of the irreducible
$\mathfrak{sl}_2$ module $U_c=\operatorname{Sym}^{c-1}K^2$:
\[
 E_B=\sum_i((c-1)x_i-x_i^2\partial_{x_i}),\qquad
 F_B=\sum_i\partial_{x_i},\qquad
 H_B=2\sum_i x_i\partial_{x_i}-n(c-1).
\]
These preserve the bounded-exponent space and commute with $S_n$.
Basis reversal exchanges $E_B$ and $F_B$ and negates $H_B$.
Thus the desired operators are
\[
 E=\Phi_\mu F_B\Phi_\mu^{-1},\qquad
 F=\Phi_\mu E_B\Phi_\mu^{-1},\qquad
 H=-\Phi_\mu H_B\Phi_\mu^{-1},
\]
on the sign-adjusted multiplicity spaces. The inverse means the unique
bounded-exponent preimage under the explicit isomorphism above.
The Cartan eigenvalue follows from $a=L-h$. Decomposition into finite
dimensional $\mathfrak{sl}_2$ modules gives the stated Lefschetz
isomorphisms. Since all row hooks are Vandermonde-generated, these
operators act on the same actual subspaces in both modules.
\end{proof}

\section{Tensor coinvariants, Schur resolutions and finite jets}
\label{sec:tensor-applications}

The selected-label presentation retains the algebra of finite jets.
A canonical Schur cover then describes every label and all its moment
homology. Retain the augmented
algebra $B$, its augmentation ideal $\mathfrak m$, and $Q_n(B)$ from
Section~\ref{sec:selected-quotients}; initially, only invertibility of
$n!$ in $K$ is required.

\begin{theorem}[First parabolics are jets]
\label{thm:universal-first-parabolic-jet}
For every augmented commutative algebra $B$ as above, inclusion of
the first tensor factor gives a natural algebra isomorphism
\begin{equation}\label{eq:universal-first-parabolic-jet}
                  Q_n(B)^{S_{n-1}}\cong B/\mathfrak m^n.
\end{equation}
For $n\ge2$, let $V_n=S^{(n-1,1)}$ be the standard representation.
There is a natural isomorphism
\begin{equation}\label{eq:universal-standard-jet}
 \operatorname{Hom}_{S_n}(V_n,Q_n(B))
       \cong \mathfrak m/\mathfrak m^n
\end{equation}
of modules for every augmentation-preserving symmetry group of $B$.
If $B$ is graded and $\varepsilon$ preserves the grading, both
isomorphisms preserve the grading.
\end{theorem}

\begin{proof}
Take $r=1$ in Theorem~\ref{thm:universal-selected-trace-presentation}.
The relation is $b^n=0$ for every $b\in\mathfrak m$.
The ideal generated by these powers is exactly $\mathfrak m^n$:
the reverse inclusion follows by polarizing the $n$th power, or
directly from
\[
 n!\,b_1\cdots b_n
  =\sum_{J\subseteq\{1,\ldots,n\}}(-1)^{n-|J|}
       \left(\sum_{j\in J}b_j\right)^n.
\]
Moreover $Q_n(B)^{S_n}=K$, by Reynolds averaging and the definition
of the invariant-generated ideal.  The permutation module
$\operatorname{Ind}_{S_{n-1}}^{S_n}\mathbf1$ splits as
$\mathbf1\oplus V_n$.  Frobenius reciprocity therefore identifies
the standard multiplicity space with the augmentation kernel of
$Q_n(B)^{S_{n-1}}$, which is $\mathfrak m/\mathfrak m^n$ under
\eqref{eq:universal-first-parabolic-jet}.
Concretely, the reciprocity map evaluates a homomorphism on
$e_1-n^{-1}\sum_i e_i$, so the identification is natural.
\end{proof}

For $B=\Sym V$, the standard multiplicity space is
$\bigoplus_{d=1}^{n-1}\Sym^dV$. This polynomial specialization
recovers the bosonic standard-character formula of
Jiang--Lentfer~\cite[Proposition 6.5 and Theorem 6.6]{JiangLentfer2026};
their theorem also treats fermionic alphabets. The statement above
identifies the algebra and its augmentation kernel for arbitrary
augmented $B$.

\subsection{Canonical Schur covers and their relations}
\label{sec:schur-presentation}

For the remainder of this section $K$ has characteristic zero, and
$B=K\oplus\mathfrak m$ is any augmented commutative algebra. Put
$S_B=\Sym(\mathfrak m)$; its generator $T_b$, linear in $b\in\mathfrak m$,
acts on $P_n(B)=B^{\otimes n}$ by $p_n(b)=\sum_jb^{(j)}$.
The trace presentation in Section~\ref{sec:selected-quotients} identifies
$Q_n(B)$ with $P_n(B)/(S_B)_+P_n(B)$.

The augmentation decomposes $P_n(B)$ by occupied positions:
\[
 P_n(B)=\bigoplus_{j=0}^nP_n^{\langle j\rangle}(B),\qquad
 P_n^{\langle j\rangle}(B)=
       \bigoplus_{\substack{I\subseteq[n]\\|I|=j}}\mathfrak m^{\otimes I},
\]
with units elsewhere. If $e_\mu$ is the central label idempotent, define
the actual subspace
\[
 L_n(B)=\bigoplus_{\mu\vdash n}
          e_\mu P_n^{\langle n-\mu_1\rangle}(B)\subseteq P_n(B).
\]

\begin{lemma}[Minimum occupancy]
\label{lem:minimum-occupancy-cover}
The quotient map $L_n(B)\to Q_n(B)$ is surjective. For
$\mu=(t,\tau)\vdash n$, its label-$\mu$ source multiplicity is
$\Schur_\tau\mathfrak m$.
\end{lemma}
\begin{proof}
Put $s=|\tau|$. Selected elimination makes
$B^{\otimes s}\to Q_n(B)^{S_t}$ surjective. A selected tensor with
a unit occupies at most $s-1$ positions, so its global orbit is induced
from a subgroup containing $S_{t+1}$. Its label first row is therefore
at least $t+1$, and its $\mu$ projection is zero. Fully occupied
tensors cover the $S_t$-fixed part of label $\mu$, and its $S_n$-span
covers the whole isotype. For the trivial label this says
$Q_n(B)^{S_n}=K$.

Finally,
$P_n^{\langle s\rangle}=\operatorname{Ind}_{S_s\times S_t}^{S_n}
(\mathfrak m^{\otimes s}\otimes\mathbf1)$.
If $|\nu|=s$ and $\mu/\nu$ is horizontal, Pieri interlacing gives
$\nu_i\ge\tau_i$; equality of sizes forces $\nu=\tau$.
\end{proof}

For $\mu=(t,\tau)$ write $a=\tau_1$, with $a=0$ for the empty tail,
and put
\[
 \tau^+=(t+1,\tau_2,\tau_3,\ldots),\qquad
 d_1=|\tau^+|=n+1-a,
 \qquad M_\mu(B)=\Hom_{S_n}(S^\mu,Q_n(B)).
\]
Fix the usual one-dimensional Pieri identification lines.

\begin{theorem}[The two-Schur presentation]
\label{thm:two-schur-presentation}
There are natural exact presentations
\begin{align*}
 S_B\otimes\Schur_{\tau^+}\mathfrak m
 &\xrightarrow{R_{\mu,B}}
 S_B\otimes\Schur_\tau\mathfrak m
 \xrightarrow{\Phi_{\mu,B}}\Schur_\mu(B)\longrightarrow0,\\
 \Schur_{\tau^+}\mathfrak m
 &\xrightarrow{\beta_{\mu,B}}\Schur_\tau\mathfrak m
               \longrightarrow M_\mu(B)\longrightarrow0.
\end{align*}
Here $\Phi$ is actual multiplication by traces on minimum occupancy.
All maps commute with arbitrary augmented algebra maps, including
ungraded and noninjective maps, and preserve compatible gradings and
auxiliary symmetries. The left-hand maps need not be injective.
\end{theorem}

\begin{proof}
Filter $F_\mu=S_B\otimes\Schur_\tau\mathfrak m$ by input count:
$j$ formal traces and $|\tau|$ inputs have count $j+|\tau|$.
Filter the target by occupied positions. Collisions lower this count,
so the associated-graded map is the same cover for the square-zero
algebra $B_0=K\oplus\mathfrak m$. Giving $\mathfrak m$ degree one,
Lemma~\ref{lem:minimum-occupancy-cover} makes this cover surjective by
degree induction. Pieri and branching off the constant factor give
\begin{align*}
 F_\mu&=\bigoplus_{\nu/\tau\ \mathrm{horizontal}}\Schur_\nu\mathfrak m,\\
 \Schur_\mu(K\oplus\mathfrak m)
   &=\bigoplus_{\substack{\nu/\tau\ \mathrm{horizontal}\\\nu_1\le t}}
                                      \Schur_\nu\mathfrak m.
\end{align*}
Let $G_\mu$ be the sum of the retained source summands $\nu_1\le t$,
and $Z_\mu$ the omitted sum. These are canonical Schur isotypic
subfunctors. The square-zero map is zero on $Z_\mu$ and an isomorphism
on $G_\mu$. For $B$, the restriction
$\Psi_B=\Phi_{\mu,B}|_{G_\mu}$ has the same diagonal and only
lower-count terms. Retained counts lie between $|\tau|$ and $n$,
so its triangular inverse is finite. The entire kernel is the graph
\begin{equation}\label{eq:schur-relation-graph}
 Z_\mu\xrightarrow{\sim}\ker\Phi_{\mu,B},\qquad
 z\longmapsto z-\Psi_B^{-1}\Phi_{\mu,B}(z).
\end{equation}

Every omitted $\nu$ contains $\tau^+$ and $\nu/\tau^+$ is horizontal.
The symmetric Pieri composition lemma \cite[Lemma~2.1]{SamWeyman}
and multiplicity one therefore give
$Z_\mu=S_B\cdot\Schur_{\tau^+}\mathfrak m$ for the square-zero map.
Restrict \eqref{eq:schur-relation-graph} to this generating summand
and extend $S_B$-linearly. Subtracting lifted generators with a prescribed
leading term lowers count, so finite induction proves that they generate
the actual kernel. Reduction modulo $(S_B)_+$ gives the second
presentation, with the convention
$\beta_{\mu,B}=-\operatorname{ev}_{T=0}R_{\mu,B}$.

Schur projections, multiplication, and the inverse of the natural
isomorphism $\Psi_B$ commute with augmented maps. For infinite-dimensional
$\mathfrak m$, the square-zero constructions pass through finite-dimensional
vector subspaces; the formulas use finitely many products on each input.
\end{proof}

The graph is an identification of filtered vector spaces, not generally
of $S_B$-modules. For example,
\[
 M_{(2,2)}(B)=\Sym^2\mathfrak m\big/
 \Span\{ab\odot c+ac\odot b+bc\odot a:a,b,c\in\mathfrak m\}.
\]
The standard label gives $\mathfrak m/\mathfrak m^n$, as above.
If $\mathfrak m^2=0$, every $\beta$ is zero and
$M_{(t,\tau)}(B)=\Schur_\tau\mathfrak m$ in every rank.

\begin{proposition}[The first obstruction]
\label{prop:first-schur-obstruction}
Suppose $B$ is connected positively graded. The actual free cover on
label $\mu$ is an isomorphism in degrees less than $d_1=n+1-\tau_1$,
and every positive moment Tor group on this label vanishes there.
In degree $d_1$ there is a natural exact sequence
\begin{align*}
 0&\longrightarrow
 \Hom_{S_n}\bigl(S^\mu,\Tor_1^{S_B}(P_n(B),K)_{d_1}\bigr)
 \longrightarrow\Schur_{\tau^+}(B_1)\\
 &\xrightarrow{\beta_{\mu,B}}(\Schur_\tau\mathfrak m)_{d_1}
       \longrightarrow M_\mu(B)_{d_1}\longrightarrow0.
\end{align*}
All higher Tor groups on this label vanish in degree $d_1$.
More generally, $i$th Tor vanishes in degrees less than $d_1+i-1$.
\end{proposition}
\begin{proof}
Every omitted partition has size at least $d_1$, with equality only
for $\tau^+$. The graph construction identifies the first possible
kernel piece with $\Schur_{\tau^+}(B_1)$. Tensor the free-cover
sequence with $K$. Positive coefficients cannot reach that degree,
giving the exact sequence. The Koszul chains of the kernel in
homological degree $i-1$ have degree at least $d_1+i-1$.
\end{proof}

\subsection{A full resolution and its polynomial first degrees}
\label{sec:schur-resolution-statements}

The presentation extends in every homological degree. Extend the parts
of $\mu$ by zero and set
\begin{equation}\label{eq:schur-resolution-shapes}
 \rho_i=(\mu_1+1,\ldots,\mu_i+1,\mu_{i+2},\mu_{i+3},\ldots),
 \qquad d_i=|\rho_i|=n-\mu_{i+1}+i.
\end{equation}
Thus $\rho_0=\tau$ and $\rho_1=\tau^+$. The operation deletes row
$i+1$ and adds one box to each preceding row.

\begin{theorem}[A natural filtered Schur resolution]
\label{thm:full-schur-resolution}
For every augmented commutative $B$ in characteristic zero,
$\Schur_\mu(B)$ has a natural $S_B$-free resolution with terms
\[
                       F_i=S_B\otimes\Schur_{\rho_i}\mathfrak m.
\]
The augmentation and first differential are those of
Theorem~\ref{thm:two-schur-presentation}. Its associated graded for
input count is the classical pure Pieri resolution; the lifted
resolution need not be minimal. Modulo $(S_B)_+$ it gives a complex
\begin{equation}\label{eq:small-schur-moment-complex}
 \cdots\longrightarrow\Schur_{\rho_2}\mathfrak m
 \xrightarrow{\partial_2}\Schur_{\rho_1}\mathfrak m
 \xrightarrow{\partial_1}\Schur_\tau\mathfrak m
\end{equation}
computing $\Hom_{S_n}(S^\mu,\Tor_i^{S_B}(P_n(B),K))$.
The differential $\partial_i$ uses exactly
$e_i=\mu_i-\mu_{i+1}+1$ multiplications. It is natural for all
augmented algebra maps and preserves compatible gradings.
For $\dim\mathfrak m=q$, the terms with $i>q$ vanish.
\end{theorem}

Appendix~\ref{app:schur-resolution} constructs the maps by a natural
creation--annihilation contraction. Its Laplacian has a
dimension-independent content eigenvalue, and multiplication is a
finite perturbation. This constructs the full resolution without choices
of higher splittings.

\paragraph{Nilpotent augmentation ideals.}
If $\mathfrak m^L=0$ with $L\ge2$, then $\partial_i=0$ whenever
$d_i>(L-1)d_{i-1}$: every term partitions its $d_i$ inputs into
$d_{i-1}$ products, and one product must contain at least $L$ inputs.
In particular, the canonical presentation gives
\[
 M_{(n-|\tau|,\tau)}(B)\cong\Schur_\tau\mathfrak m
 \qquad\bigl(n\ge(L-1)|\tau|+\tau_1\bigr).
\]
This equality holds in every internal degree whenever $B$ is graded.
For square-zero algebras it holds at every rank supporting the label.

\begin{corollary}[Sharp universal homological support]
\label{cor:universal-schur-tor-support}
For connected positively graded $B$ and $i\ge1$, label $\mu$ in
degree-$D$ moment homology requires
\[
 D\ge n+i-\mu_{i+1},\qquad
 D\ge n+i-\left\lfloor\frac n{i+1}\right\rfloor.
\]
Both bounds are attained by square-zero algebras with sufficiently
large degree-one space.
\end{corollary}
\begin{proof}
The $i$th term of \eqref{eq:small-schur-moment-complex} begins in
internal degree $d_i$, and $\mu_{i+1}\le\lfloor n/(i+1)\rfloor$.
For square-zero multiplication all differentials are zero. A partition
maximizing its $(i+1)$st part attains the second threshold.
\end{proof}

For a polynomial algebra, the first possible degree behaves better.
The proof uses an ungraded specialization to orthogonal idempotents;
the full naturality in Theorem~\ref{thm:two-schur-presentation} is essential.

\begin{theorem}[Exactness in the first possible homological degree]
\label{thm:polynomial-schur-fronts}
Let $A=\Sym(V)$ and $S_A=\Sym(A_+)$. If $\mu\ne(n)$ and $i\ge1$,
then
\[
 \Hom_{S_n}\bigl(S^\mu,\Tor_i^{S_A}(A^{\otimes n},K)_D\bigr)=0
                         \qquad(D\le n+i-\mu_{i+1}).
\]
The same conclusion holds for $\mu=(n)$ and $i\ge2$.
The sole exceptional first degree is
\[
 \Hom_{S_n}\bigl(\mathbf1,
       \Tor_1^{S_A}(A^{\otimes n},K)_{n+1}\bigr)=\Sym^{n+1}V.
\]
Equivalently, every supported all-linear restriction
$\partial_i:\Schur_{\rho_i}V\to
(\Schur_{\rho_{i-1}}A_+)_{d_i}$ is injective, except for
$\mu=(n),i=1$, where it is zero.
\end{theorem}

\begin{proof}
Take the split augmented algebra $C=K^{q+1}$, with augmentation its
zeroth coordinate and augmentation ideal spanned by the orthogonal
idempotents $u_1,\ldots,u_q$. On a primitive tensor idempotent of
$C^{\otimes n}$, the moment $p_n(u_a)$ acts by the number of positions
of color $a$. Unless every position has color zero, some such count
is a positive integer, and its Koszul complex contracts. Consequently
\begin{equation}\label{eq:idempotent-moment-homology}
 Q_n(C)=K,\qquad
 \Tor_i^{\Sym(C_+)}(C^{\otimes n},K)
                      =\Lambda^i(C_+)\otimes\mathbf1.
\end{equation}

Suppose $\mu\ne(n)$ and $\Schur_{\rho_i}K^q\ne0$.
The complex \eqref{eq:small-schur-moment-complex} for $C$ is exact
in every degree on this label. Its term dimensions are those of the
square-zero pure resolution. Localizing that resolution at the fraction
field of $\Sym(K^q)$ also gives an exact complex: the square-zero
coefficient module has finite length. Exactness determines the ranks
of all differentials recursively from the term dimensions. Each
supported pure Pieri differential is nonzero, so $\partial_i$ for
$C$ is nonzero as well.

Set $q=\dim V$ and extend a linear isomorphism $V\to C_+$ to an
augmented algebra map $A\to C$. Naturality now makes the restriction
of $\partial_i$ to $\Schur_{\rho_i}V$ nonzero. It is $\GL(V)$-equivariant
with irreducible source, hence injective. If that Schur module is zero,
there is nothing to prove.

For $\mu=(n)$ and $i\ge2$, take $q=i$. The top small-complex term is
$\Schur_{(n+1,1^{i-1})}K^i$, of dimension
$\binom{n+i-1}{n}>1$, whereas its kernel has dimension one by
\eqref{eq:idempotent-moment-homology}. Thus its differential is nonzero.
The same specialization, followed by naturality in the row space,
proves injectivity in every supporting dimension. For $i=1$ the
target is $K$, so the positive-degree map is zero.

Finally, the degree-$d_i$ part of the $i$th Schur term uses only
$V=A_1$, and the next term starts strictly above $d_i$. The assertions
follow from the small complex and the preceding degree bounds.
\end{proof}

\begin{theorem}[The first character correction in every degree]
\label{thm:first-unstable-specht-wall}
For every nonempty $\tau\vdash s$, every $D>s$, and
\[
 n=D+\tau_1-1,\qquad
 \mu=(D-s+\tau_1-1,\tau),\qquad
 \widehat\tau=(D-s+\tau_1,\tau_2,\tau_3,\ldots),
\]
there is a natural exact sequence
\begin{equation}\label{eq:first-unstable-specht-wall}
 0\longrightarrow\Schur_{\widehat\tau}V
 \xrightarrow{\beta_{\mu,A}}(\Schur_\tau A_+)_D
 \longrightarrow\Hom_{S_n}(S^\mu,R_{n,k}[D])\longrightarrow0.
\end{equation}
All positive moment homology on this label vanishes through degree $D$.
Thus, when $\dim V\ge\ell(\tau)$, the stable threshold
$n\ge D+\tau_1$ is sharp in every degree $D>s$.
\end{theorem}
\begin{proof}
Here $D=d_1$ and $\rho_1=\widehat\tau$. Apply
Proposition~\ref{prop:first-schur-obstruction} and
Theorem~\ref{thm:polynomial-schur-fronts}. The correction has length
$\ell(\tau)$ and is nonzero under the stated dimension hypothesis.
\end{proof}

This identifies the canonical relation map at the first rank below
the stable threshold.

\subsection{Natural stable multiplicity spaces}
\label{sec:stable-specht-multiplicities}

Let $B=K\oplus\mathfrak m$ be a connected positively graded
commutative algebra over a field of characteristic zero. Retain $Q_n(B)$ from
Theorem~\ref{thm:universal-selected-trace-presentation}.
For $\tau\vdash s$ write $\tau[n]=(n-s,\tau)$ when this is a partition.
For the empty partition set $\tau_1=0$.

\Needspace{13\baselineskip}
\begin{samepage}
\begin{theorem}[Natural stable multiplicity spaces]
\label{thm:stable-augmented-multiplicities}
If $D\ge s$ and $n\ge D+\tau_1$, then
\[
 \operatorname{Hom}_{S_n}(S^{\tau[n]},Q_n(B)_D)
 \simeq (\mathbb S_\tau\mathfrak m)_D.
\]
Fixing the one-dimensional Pieri branching line makes the
isomorphism natural for graded augmented algebra maps and
auxiliary symmetries preserving the grading and augmentation.
For $D<s$ the left side is zero whenever $\tau[n]$ is a partition.
\end{theorem}
\end{samepage}
\begin{proof}
The canonical cover in Theorem~\ref{thm:two-schur-presentation} has
no relation in degrees below $n+1-\tau_1$, by
Proposition~\ref{prop:first-schur-obstruction}. Its quotient modulo
positive traces is therefore exactly $\Schur_\tau\mathfrak m$ in
the stated range. This is the actual minimum-occupancy quotient map,
so it has the asserted naturality. For $D<s$, occupancy gives zero.
\end{proof}

\begin{corollary}[Universal separation and its sharp range]
\label{cor:universal-augmented-separation}
For $D\ge1$ and $n\ge2D-1$, multiplication gives a canonical isomorphism
\[
 \bigl(P_n(B)^{S_n}\otimes L_n(B)\bigr)_{\le D}
                     \xrightarrow{\sim}P_n(B)_{\le D},
\]
and all positive moment homology vanishes through degree $D$.
The bound is sharp for arbitrary connected graded augmented algebras.
\end{corollary}
\begin{proof}
A degree-$d$ label has tail size at most $d$. At $n=2d-1$, a tail
with first part $d$ would have to be $(d)$ and would give the invalid
label $(d-1,d)$. Thus every occurring label satisfies the individual
bound $n\ge d+\tau_1$; larger $n$ are immediate. The trivial-label
presentation shows that $S_B$ and the actual invariant algebra agree
through degree $n$, giving the multiplication statement.

For sharpness let $D\ge2$, $n=2D-2=2m$, and $B=K[x]/(x^2)$.
The sole moment $p=\sum_{i=1}^{2m}x_i$ kills the nonzero element
$f=\prod_{j=1}^m(x_{2j-1}-x_{2j})$. Hence $f e_{T_x}$ is nonzero
first Koszul homology in degree $D$; a one-generator complex has no
boundaries there. The endpoint $D=n=1$ is included in the isomorphism.
\end{proof}

For a fixed tail, augmentation of the last factor gives compatible
surjections $M_{\tau[n]}(B)\to M_{\tau[n-1]}(B)$ whenever both labels
are defined. With compatible Pieri lines they are quotients of the
identity on $\Schur_\tau\mathfrak m$, so their kernels decrease with
$n$. Indeed, the occupancy-$s$ projection is onto, and among its source
labels of first row at least $n-s$, only $\tau[n]$ branches to
$\tau[n-1]$. Multiplicity monotonicity itself is known in the FI-module
setting \cite[Lemma~2.1.1]{Borisov2026}; here the common source and maps
are explicit for arbitrary augmented algebras.

\subsection{Polynomial algebras and positive plethysm}

Write
\[
 A=\operatorname{Sym}V,
 \qquad A_+=\bigoplus_{d\ge1}\operatorname{Sym}^dV,
 \qquad H(\mathbf q)=\prod_{a=1}^k(1-q_a)^{-1}.
\]
The grading on the Schur functor $\mathbb S_\tau(A_+)$
is induced by the usual polynomial grading of $A_+$.
In particular, every homogeneous part is finite dimensional.

\begin{corollary}[Explicit stable Specht multiplicities]
\label{cor:explicit-stable-specht-multiplicities}
For $D\ge1$ and $n\ge2D-1$, there is an isomorphism of
$\operatorname{GL}(V)\times S_n$-modules
\begin{equation}\label{eq:stable-specht-schur-functor}
 (R_{n,k})_D\cong
 \bigoplus_{\substack{|\tau|\le D\\\tau[n]\ \mathrm{a\ partition}}}
   \bigl(\mathbb S_\tau(A_+)\bigr)_D
       \otimes S^{(n-|\tau|,\tau)}.
\end{equation}
Thus the multigraded stable multiplicity of the label
$\tau[n]=(n-|\tau|,\tau)$ is the positive plethysm
\begin{equation}\label{eq:stable-specht-plethysm}
 s_\tau[H(\mathbf q)-1].
\end{equation}
For a fixed nonempty $\tau$ and
$D\ge|\tau|$, its degree-$D$ multiplicity is
already given by this formula when
\begin{equation}\label{eq:individual-stable-specht-range}
 n\ge D+\tau_1.
\end{equation}
The term for $\tau=\varnothing$ is the trivial label,
which occurs only in degree zero.
\end{corollary}

\begin{proof}
Apply Theorem~\ref{thm:stable-augmented-multiplicities} to
$B=\operatorname{Sym}V$. Corollary~\ref{cor:universal-augmented-separation} covers
every valid label when $n\ge2D-1$. The same theorem gives the sharper individual
range for every $D\ge1$. Schur functor characters give
\eqref{eq:stable-specht-plethysm}, in each multidegree.
\end{proof}

For example, the tails $\tau=(1),(2),(1,1)$ give
$A_+$, $\operatorname{Sym}^2(A_+)$ and $\bigwedge^2(A_+)$,
respectively, before expansion in row Schur functions.

\subsection{Finite jets and their maps}

The stable multiplicities depend on the graded augmentation space.
The compatible parabolic algebras retain more: they recover
multiplication, and also every map of the finite jets.  Here the
algebras need not be graded.  Put $J_n(B)=B/\mathfrak m_B^n$.

\begin{theorem}[Full faithfulness of finite jets]
\label{thm:finite-jet-full-faithfulness}
Let $n\ge2$, let $K$ be a field in which $n!$ is invertible, and
let $B,C$ be augmented commutative $K$-algebras.  There is a natural
bijection
\begin{equation}\label{eq:finite-jet-full-faithfulness}
 \Hom_{S_n\text{-alg}}(Q_n(B),Q_n(C))
 \cong
 \Hom_{\mathrm{aug}\text{-alg}}(J_n(B),J_n(C)).
\end{equation}
Maps on the left are unital; they automatically preserve augmentation.
Thus $Q_n$ is fully faithful on augmented algebras whose augmentation
ideals have $n$th power zero.  If the two jets are finite dimensional,
\eqref{eq:finite-jet-full-faithfulness} is an isomorphism of affine
schemes of homomorphisms, including the loci of isomorphisms.
\end{theorem}

\begin{proof}
By Theorem~\ref{thm:universal-first-parabolic-jet}, the first factor
map $B\to Q_n(B)$ kills $\mathfrak m_B^n$.  The same is true in
every factor, by permutation of labels.  Consequently
\begin{equation}\label{eq:finite-jet-truncation}
                         Q_n(B)\cong Q_n(J_n(B)).
\end{equation}
Indeed, the factor maps give the inverse to the evident quotient
map: the invariant ideal is generated by global traces of
augmentation elements, and these traces vanish under the factor
maps.

Reynolds averaging and the definition give $Q_n(B)^{S_n}=K$.
Its Reynolds projection therefore equals its augmentation, since
augmentation commutes with averaging and is the identity on $K$.
An equivariant unital algebra map commutes with this projection
and hence preserves augmentation.  Restrict it to the invariants
of $S_{n-1}$ fixing the first factor.  The first-parabolic theorem
gives an augmented map $J_n(B)\to J_n(C)$.  Conversely, apply
$Q_n$ to any such map and use \eqref{eq:finite-jet-truncation}.
These operations are inverse: $Q_n(B)$ is generated by the
$S_n$-orbit of its first factor.  They are natural in both algebras.

For finite-dimensional jets, repeat this argument after every
commutative $K$-algebra extension $R$.  Tensor powers, the
invariant-generated ideal, and Reynolds idempotents commute with
base change; the first parabolic is $J_n(B)\otimes_K R$.
The same inverse constructions therefore work over nonreduced
$R$ as well.  Yoneda's lemma identifies the affine homomorphism
schemes.  Applying the construction to a map and its inverse
identifies the isomorphism loci.
\end{proof}

Applying augmentation to the last tensor factor sends invariant traces
on $n$ factors to those on $n-1$ factors, defining $Q_n(B)\to Q_{n-1}(B)$.
On first parabolics this is the quotient map
$B/\mathfrak m_B^n\to B/\mathfrak m_B^{n-1}$.

\begin{corollary}[Formal towers]
\label{cor:formal-towers}
In characteristic zero, compatible equivariant algebra maps between
the towers $\{Q_n(B)\}_{n\ge2}$ and $\{Q_n(C)\}_{n\ge2}$ are
exactly the continuous augmented algebra maps
$\widehat B_{\mathfrak m_B}\to\widehat C_{\mathfrak m_C}$.
The transition maps augment the last tensor factor.
\end{corollary}

\begin{proof}
The theorem identifies compatible tower maps with compatible maps
of the jet quotients.  Taking inverse limits gives a continuous augmented
map of completions.  Conversely, such a map takes the closure of
$\mathfrak m_B^n$ into the closure of $\mathfrak m_C^n$: approximate
elements of the former by sums of products of $n$ augmentation
elements and use continuity.  It therefore induces the required
maps at every finite level.
\end{proof}

For example, take
\[
 B_1=K[x]/(x^3),\quad \deg x=1,
 \qquad
 B_2=K[u,v]/(u,v)^2,\quad \deg u=1,\quad\deg v=2.
\]
Their graded augmentation spaces are isomorphic: each has one
basis vector in degrees one and two.  After choosing this
identification, Theorem~\ref{thm:stable-augmented-multiplicities}
gives the same stable multiplicity spaces.  Their first parabolics
at $n=2$ differ, since $\dim J_2(B_1)=2$ and
$\dim J_2(B_2)=3$.  Thus the stable representation and the
parabolic algebra record different information.

\begin{corollary}[Equivariant automorphisms of polynomial coinvariants]
\label{cor:polynomial-jet-automorphisms}
Let $n\ge2$, let $K$ have characteristic zero, and let $\dim V=k$.
For $R_{n,k}=Q_n(\Sym V)$, the equivariant
automorphism group scheme is the group of origin-fixing coordinate
substitutions modulo degree $n$ with invertible linear part.  It is
smooth of dimension
\[
                 k\left(\binom{k+n-1}{k}-1\right),
\]
and its graded subgroup is $\GL(V)$.
\end{corollary}

\begin{proof}
An augmented endomorphism of $\Sym V/(V)^n$ is specified freely
by the image of $V$ in $(V)/(V)^n$.  It is invertible exactly
when its linear part is invertible, as an inverse can then be
constructed successively by degree.  Its scheme is the product,
as a variety, of $\GL(V)$ and the vector space of higher
coefficients.  Count these coefficients and apply
Theorem~\ref{thm:finite-jet-full-faithfulness}.
\end{proof}

\section{Stable characters in every number of alphabets}
\label{sec:stable-characters}

The full character is the Euler character of the moment complex once
positive homology vanishes. The Schur resolution controls its earliest
possible degrees, and idempotent specialization excludes their boundary.
Quadratic exactness gives the complementary label bound. These two
inputs leave a short list of first-homology shapes; bivariate and hook
exactness remove the routine cases, and one explicit cubic multiplication
resolves the last one. The finer trace bounds are retained separately in
Appendix~\ref{app:finer-character-bounds}.

\subsection{Quadratic characters and exactness}
\label{sec:quadratic-characters}

Let $n\geq2$, $V=K^k$, and let $W_n$ carry its invariant
nondegenerate symmetric form. Removing the linear moments leaves
$P^0=\Sym(V\otimes W_n)$. The quadratic moment map identifies
$\Sym^2V$ with a subspace of $P^0_2$; put
\[
 Q_{n,k}=P^0/(\Sym^2V).
\]
The orthogonal Littlewood complex computes this quotient in every
rank.  Its modification rule also makes the quadratic Koszul complex
exact in positive homological degrees through total degree $n$.
These two facts give exact initial characters and provide the
homological input for the character range proved below.

For a graded $\GL(V)\times S_n$ character, let $\Pi_n$ retain just
the row Schur terms indexed by partitions $\rho$ with
$\rho'_1+\rho'_2\leq n-1$. Here primes denote conjugate partitions.

\begin{theorem}[Orthogonal quadratic character; Sam--Snowden--Weyman]
\label{thm:quadratic-character-all-ranks}
In the degree-completed representation ring,
\begin{equation}\label{eq:quadratic-projected-euler}
 [Q_{n,k}]=\Pi_n\left(
  [\Sym(V\otimes W_n)]
       \sum_{j\geq0}(-1)^j[\Lambda^j(\Sym^2V)]\right),
\end{equation}
where $\Lambda^j(\Sym^2V)$ has polynomial degree $2j$ and trivial
$S_n$-action. Equivalently,
\begin{equation}\label{eq:quadratic-orthogonal-character}
 Q_{n,k}\simeq
 \bigoplus_{\substack{\ell(\rho)\leq k\\
                       \rho'_1+\rho'_2\leq n-1}}
 \Schur_\rho V\otimes
       \left.\Schur_{[\rho]}W_n\right|_{S_n},
\end{equation}
where $\Schur_{[\rho]}$ is the orthogonal Schur module.
\end{theorem}
\begin{proof}
The Koszul complex of the quadratic moment subspace is
$\Sym(V\otimes W_n)\otimes\Lambda^\bullet(\Sym^2V)$.
Its row-$\rho$ multiplicity complex is the orthogonal Littlewood
complex.  The orthogonal Littlewood theorem identifies its
zeroth homology with $\Schur_{[\rho]}W_n$ for
$\rho'_1+\rho'_2\leq\dim W_n$, and zero otherwise; moreover all
positive homology vanishes for these admissible partitions
\cite[\S4.2, Proposition~4.1 and Theorem~4.4]{SSW2013}.
Taking Euler characteristics only on those row components proves
\eqref{eq:quadratic-projected-euler}; taking zeroth homology gives
\eqref{eq:quadratic-orthogonal-character}. Restriction from the full
orthogonal group to $S_n\subset O(W_n)$ retains all label
multiplicities. The construction and the identities are defined
over $\mathbb Q$, so characteristic-zero scalar extension gives
the assertion over $K$.
\end{proof}

\Needspace{12\baselineskip}
\begin{samepage}
\begin{corollary}[Three initial layers of every label]
\label{cor:first-three-character-layers}
Let $\mu\vdash n$, and let $\delta_\mu(k)$ be the least degree
in which $S^\mu$ occurs in the ambient polynomial ring $P_{n,k}$.
For a row partition $\rho\vdash D$ with $\ell(\rho)\le k$, the
joint $(\rho,\mu)$ multiplicities in $R_{n,k}[D]$ and $Q_{n,k}[D]$
agree whenever
\[
 D-\delta_\mu(k)\le2\qquad\text{or}\qquad\rho_1\le2.
\]
Here the empty partition has first part zero. Thus
\eqref{eq:quadratic-projected-euler} computes the complete character
in each of the three degrees
\[
 \delta_\mu(k),\quad\delta_\mu(k)+1,\quad\delta_\mu(k)+2,
\]
for every label and every number of alphabets, in every rank $n\ge2$.
It also computes the entire two-column row sector in all degrees.
For $n=1$, the quotient is $K$ and all nonconstant components vanish.
\end{corollary}
\end{samepage}
\begin{proof}
A polarized power sum of degree $d$ has row type $\Sym^dV$ and trivial label
type. To affect the joint $(\rho,\mu)$ component, it must multiply
a polynomial component of type $(\sigma,\mu)$ in degree $D-d$,
with $\rho/\sigma$ a horizontal $d$-strip. Consequently
\[
 d\le\rho_1,\qquad D-d\ge\delta_\mu(k).
\]
Either hypothesis of the corollary excludes every $d\ge3$. Thus higher
moments have zero image in the indicated component, even before
quotienting by the linear and quadratic moments. The quadratic
quotient therefore has the full coinvariant character there, as
given by Theorem~\ref{thm:quadratic-character-all-ranks}.
\end{proof}

Theorem~\ref{thm:r5-complete-first} computes $\delta_\mu(k)$ and
generates its entire first layer. The corollary computes the next
two full harmonic layers as well; these degrees may lie outside
the uniform stable range. Operator generation in the additional
layers is governed by the criteria in Section~\ref{sec:linear-generation}.

\begin{theorem}[Quadratic Koszul exactness through degree $n$]
\label{thm:quadratic-koszul-degree-range}
For every $n\ge2$, every finite-dimensional row space $V$, and
every $0\le D\le n$,
\begin{equation}\label{eq:quadratic-koszul-degree-range}
 H_j\!\left(\Sym(V\otimes W_n)
                 \otimes\Lambda^\bullet(\Sym^2V)\right)_D=0
                 \qquad(j>0),
\end{equation}
where the exterior generators have total degree two and the
differential is multiplication by the quadratic moments.
\end{theorem}

\begin{proof}
Put $m=n-1$ and decompose the complex by its row partitions
$\lambda\vdash D$.  Its row-$\lambda$ multiplicity complex is
the orthogonal Littlewood complex.  Write
$c_i=\lambda'_i$ for the column lengths.  If
$c_1+c_2\le m$, the admissible case of the Littlewood theorem
gives zero positive homology.

For an inadmissible partition of size at most $m+1$, the inequalities
\[
          m<c_1+c_2\le |\lambda|\le m+1
\]
force $c_1+c_2=m+1$.  In the orthogonal modification rule,
the dot action of the first reflection is
\[
 s_0\mathbin{\boldsymbol\cdot}(c_1,c_2,c_3,\ldots)
       =(m+1-c_2,m+1-c_1,c_3,\ldots).
\]
It fixes this column sequence.  The partition is therefore singular,
and its Littlewood complex has zero homology in every degree
\cite[\S4.4 and Theorem~4.4]{SSW2013}.
Every row partition of size $D\le n$ is either admissible or singular,
which proves \eqref{eq:quadratic-koszul-degree-range}.
\end{proof}

Consequently, in total degrees at most $n$, the projection $\Pi_n$
in \eqref{eq:quadratic-projected-euler} can be omitted: the ordinary
Koszul Euler character already equals $[Q_{n,k}]$.  Every discarded
row is singular and has Euler character zero in those degrees.
This exactness holds uniformly in the row dimension because it is
proved one Schur functor at a time; it requires no global
regular-sequence assertion for the quadratic moments.

\subsection{The full moment complex}
\label{sec:explicit-stable-characters}

Let $V=K^k$ and $P=P_{n,k}=\operatorname{Sym}(V\otimes K^n)$.
Write $p_\alpha=\sum_{i=1}^nX_i^\alpha$ for the polarized moment,
and give
$\operatorname{Sym}^j V$ weight $j$ in
$U_D=\bigoplus_{j=1}^D\operatorname{Sym}^j V$.
The moment Koszul complex is
\[
 \mathcal K_{n,k,D}=P\otimes\bigwedge U_D,
 \qquad \partial(e_\alpha)=p_\alpha.
\]
Here $e_\alpha$ is the exterior generator corresponding to the
monomial of row weight $\alpha$ in $\Sym^{|\alpha|}V$.
Its total grading counts both polynomial degree and the assigned
weights of exterior generators. Its zeroth homology agrees with
$R_{n,k}$ in total degrees at most $D$.

Our target is positive-degree exactness for
$D\ge3$, $n\ge\max\{D,2D-6\}$. Throughout this proof $D\le n$,
$\lambda$ is a row partition, and $\mu$ a label partition.
The Schur relation and quadratic reduction eliminate all homological
degrees at least two. At $n=2D-5$ the only remaining source has at
most two rows. At $n=2D-6$, three labels remain; their only component
not covered by bivariate or row-hook exactness is the multiplicity-one
cubic source in Lemma~\ref{lem:cubic-witness}.

\subsection{Parabolic detection of low-degree homology}

\begin{lemma}[Parabolic detection of low-degree homology]
\label{lem:parabolic-koszul-detection}
Fix $3\le D\le n$. In positive homological degree, every
label partition $\mu$ occurring in
$H_j(\mathcal K_{n,k,D})_D$ satisfies
\begin{equation}\label{eq:moment-koszul-homology-strip}
 n-D+3j\le\mu_1\le D-2.
\end{equation}
In any total degree, if the moments through order $a\ge1$ form a
regular sequence in $P$, one also has the lower bound
$\mu_1\ge n-D+(a+1)j$.
\end{lemma}

\begin{proof}
Remove the linear moments and compute quadratic Koszul homology first.
By Theorem~\ref{thm:quadratic-koszul-degree-range}, in total degrees at
most $n$ the remaining complex is
\[
 Q_{n,k}\otimes\Lambda^\bullet
       \left(\bigoplus_{d\ge3}\Sym^dV\right),
 \qquad Q_{n,k}=\Sym(V\otimes W_n)/(\Sym^2V).
\]
A degree-$b$ coefficient has first label row at least $n-b$, because
each monomial occupies at most $b$ labels; the bound passes to quotients.
The exterior factors consume at least $3j$ degrees, giving the lower
bound. If moments through order $a$ are regular, quotient by them first
and use the minimum remaining exterior weight $a+1$.

The trivial label has no positive moment homology through degree $n$,
by Theorem~\ref{thm:polynomial-schur-fronts}. For any other label, that
same theorem gives the stronger restriction
\begin{equation}\label{eq:polynomial-higher-label-support}
                  \mu_{j+1}\ge n-D+j+1.
\end{equation}
In particular $\mu_2\ge n-D+2$, so $\mu_1\le D-2$.
More generally, summing the first $j+1$ label rows yields
\[
 \mu_1\le jD-(j-1)n-j(j+1),\qquad
 n\le(1+1/j)D-j-4,\qquad D\ge j(j+4).
\]
The latter two inequalities combine this estimate with the quadratic
lower bound and $D\le n$. For example, second homology is absent in
this range whenever $D<12$.
\end{proof}

\begin{proposition}[Intersection support]
\label{prop:moment-intersection-support}
If $\Schur_\lambda V\otimes S^\mu$ occurs in $H_i$ in degree $D$
with $i\ge1$,
then, with $\rho_i,d_i$ from \eqref{eq:schur-resolution-shapes},
\[
 |\lambda\cap\rho_i|\ge2d_i-D,\qquad
 D+\lambda_1\ge n+\mu_1-\mu_{i+1}+i+1.
\]
Here the intersection is that of Young diagrams.
\end{proposition}
\begin{proof}
Inside $\Schur_{\rho_i}(A_+)_D$, at most $D-d_i$ inputs have degree
at least two. Splitting $A_+=V\oplus\Sym^{\ge2}V$ gives a
linear-input shape $\nu\subseteq\rho_i$ with $|\nu|\ge2d_i-D$.
Every resulting row contains $\nu$, proving the first bound; also
$\nu_1\ge(\rho_i)_1-(D-d_i)$ and $(\rho_i)_1=\mu_1+1$.
\end{proof}

Appendix~\ref{app:finer-character-bounds} gives further row-sensitive
trace bounds and detects their closed upper boundary. They are not
needed for the uniform character theorem below.

\subsection{Separation through degree \texorpdfstring{$n$}{n} in two alphabets}

A Hilbert-series identity alone need not eliminate cancelling homology.
Here multiplication already gives a free cover. Equality of dimensions
then forces its kernel to begin in a higher degree.

\begin{lemma}[Truncated freeness from a Hilbert identity]
\label{lem:truncated-freeness}
Let $S$ be a connected positively graded commutative or
supercommutative algebra over a field, and let $M$ be a nonnegatively
graded $S$-module. Assume their homogeneous pieces are finite
dimensional, and put $R=M/S_+M$. If
\[
 H_M(t)=H_S(t)H_R(t)\pmod{t^{N+1}},
\]
then every homogeneous section $R\to M$ induces an isomorphism
\[
 (S\otimes R)_{\le N}\longrightarrow M_{\le N}.
\]
Moreover, $\Tor_j^S(M,K)_d=0$ for $j>0$ and $d\le N$.
For superalgebras use ordinary dimensions and parity-preserving
sections. Additional positive multigradings may be retained.
\end{lemma}

\begin{proof}
The multiplication map is surjective by degree induction: after
subtracting a lift of the quotient class, an element belongs to
$S_+M$ and is a sum of positive-degree coefficients times
elements of smaller degree. The Hilbert identity makes the source
and target have equal dimensions through degree $N$. Its kernel
$L$ therefore has no homogeneous piece in those degrees.
Choose a graded free resolution of $L$ using generators only
in degrees greater than $N$, and splice it with
$0\to L\to S\otimes R\to M\to0$. Every positive term of the resulting
free resolution of $M$ is zero in degrees at most $N$.
Tensoring with $K$ proves the Tor assertion. All arguments respect
the additional gradings.
\end{proof}

\begin{theorem}[Bivariate separation through degree $n$]
\label{thm:bivariate-separation}
Multiplication induces a canonical $\GL_2\times S_n$-equivariant
isomorphism
\[
 \bigl(P_{n,2}^{S_n}\otimes\Harm_{n,2}\bigr)_{\le n}
          \xrightarrow{\ \sim\ }(P_{n,2})_{\le n}.
\]
Every positive moment Koszul homology group vanishes in internal
degrees at most $n$. Consequently the full bivariate $S_n$ character
is given, in bidegree $(a,b)$ whenever $a+b\le n$, by
\begin{equation}\label{eq:bivariate-character}
 \operatorname{tr}\bigl(\sigma\mid (R_{n,2})_{a,b}\bigr)
 =[q^at^b]\!
   \prod_{\substack{i,j\ge0\\i+j>0}}(1-q^it^j)
   \prod_{m\ge1}\bigl((1-q^m)(1-t^m)\bigr)^{-X_m(\sigma)}.
\end{equation}
Here $X_m(\sigma)$ counts the $m$-cycles of $\sigma$.
Equivalently, every row partition of length at most two has zero
positive moment homology through degree $n$, in every row dimension.
\end{theorem}

\begin{proof}
For every fixed bidegree $(a,b)$, with $D=a+b$, Wang's dimension
theorem \cite[Theorem~5.1]{Wang2026Polynomiality} gives one polynomial in $n$
equal to $\dim(R_{n,2})_{a,b}$ throughout $n\ge D$.
The Euler dimension of the moment complex is also polynomial in
$n$. Corollary~\ref{cor:universal-augmented-separation} gives eventual
acyclicity, for example for $n\ge2D-1$.
The actual and Euler dimension polynomials therefore agree at all
sufficiently large integers and are identical. Degrees zero,
one, and two follow from the independent linear and quadratic
moments. Hence the quotient has the Euler Hilbert series in every
total degree at most $n$.

Let
\[
 F_n=\Sym\left(\bigoplus_{j=1}^n\Sym^jK^2\right)
\]
act on $P_{n,2}$ by the formal moment evaluation. Polarized moments of
degrees at most $n$ generate the invariant algebra, so
$R_{n,2}=P_{n,2}/(F_n)_+P_{n,2}$.
The dimension
identity just proved is
$H_{P_{n,2}}=H_{F_n}H_{R_{n,2}}$ through degree $n$.
Lemma~\ref{lem:truncated-freeness} gives freeness and positive
Koszul acyclicity in this range.

The fundamental trace presentation has its first relation in
degree $n+1$: each relation uses $n+1$ positive inputs
\cite{Domokos2017}. Thus $F_n$ maps isomorphically
onto the actual invariant algebra in degrees at most $n$.
Choose the canonical harmonic section of $R_{n,2}$. The free
cover becomes exactly the multiplication map in the statement.

Taking equivariant characters of the moment complex gives
\eqref{eq:bivariate-character}. Each $m$-cycle contributes
$(1-q^m)^{-1}(1-t^m)^{-1}$ to the polynomial trace, and each
invariant moment contributes its Euler factor. Finally, the
complex is a polynomial functor of the row space. In
characteristic zero its homology decomposes into Schur functors;
every partition of length at most two is detected on $K^2$.
\end{proof}

\subsection{Row-hook exactness}
\label{sec:hook-homology}

The Fields character has a finite principal-specialization form
that records both the exactness range and its first obstruction.
Write $(q;q)_m=\prod_{i=1}^m(1-q^i)$ and
$(q;q)_\infty=\prod_{i\ge1}(1-q^i)$.
Let $Z$ be the Frobenius alphabet for $S_n$ and set
$X=Z/(1-q)$ in plethystic notation. This means
$p_d[X]=p_d[Z]/(1-q^d)$.

Let $\mathscr H_n(q,u;Z)$ be the Frobenius series of the
$(1,1)$ bosonic--fermionic coinvariant ring: $q$ records polynomial
degree and $u$ exterior degree. Let $\mathscr E_n$ be the
moment Euler character after evaluation on a row superspace of
dimension $(1|1)$.

\begin{lemma}[Finite Fields identity and its first discrepancy]
\label{lem:fields-euler}
For $0\le b\le n$,
\begin{equation}\label{eq:finite-fields}
 [u^b]\mathscr H_n
  =\sum_{h=0}^b(-1)^{b-h}
       \frac{(q;q)_{n-h}}{(q;q)_{b-h}}e_h[X]h_{n-h}[X].
\end{equation}
The corresponding coefficient of $\mathscr E_n$ is obtained by
replacing $(q;q)_{n-h}$ with $(q;q)_\infty$. Hence the two
characters agree in every total degree at most $n$, and
\begin{equation}\label{eq:first-fields-discrepancy}
 [q^{n-b+1}u^b](\mathscr H_n-\mathscr E_n)
       =e_b[Z]h_{n-b}[Z]\qquad(0\le b\le n).
\end{equation}
\end{lemma}

\begin{proof}
The Fields theorem \cite[Theorem~5.2 and equation~(2.5)]{MRW} gives
\[
 \mathscr H_n
 =\sum_{\mu\vdash n}s_\mu[Z]
      \sum_{T\in\SYT(\mu)}
        q^{\maj(T)}
        \prod_{r=1}^{\des(T)}(1+uq^{-r}).
\]
Put
\[
 \mathcal A_n(t,q;Z)=
   \sum_{\mu\vdash n}s_\mu[Z]
       \sum_{T\in\SYT(\mu)}t^{\des(T)}q^{\maj(T)}.
\]
Principal specialization of Schur functions gives the
character-valued Eulerian identity
\[
 \frac{\mathcal A_n(t,q;Z)}
      {\prod_{i=0}^n(1-tq^i)}
  =\sum_{m\ge0}
       h_n\left[Z\frac{1-q^{m+1}}{1-q}\right]t^m.
\]
One obtains it by standardizing semistandard tableaux and summing
the weakly increasing entries with forced rises at descents.
The standardization convention using comajor index gives the same
identity, since tableau evacuation preserves the number of
descents and interchanges major and comajor indices.

Expand
\[
 h_n[X-q^{m+1}X]
   =\sum_{h=0}^n(-1)^h q^{(m+1)h}e_h[X]h_{n-h}[X].
\]
Multiplication by the denominator above and evaluation at
$t=q^{-h}$ leave just one summand, yielding
\[
 \mathcal A_n(q^{-h},q;Z)
 =q^{-\binom h2}(q;q)_h(q;q)_{n-h}e_h[X]h_{n-h}[X].
\]
The elementary-symmetric identity
\[
 e_b(q^{-1},\ldots,q^{-d})
  =\sum_{h=0}^b
       \frac{(-1)^{b-h}q^{\binom h2}}
            {(q;q)_h(q;q)_{b-h}}q^{-hd}
\]
follows from the finite $q$-binomial theorem. Apply it with
$d=\des(T)$ and use the preceding evaluation to obtain
\eqref{eq:finite-fields}.

The polynomial--exterior algebra has Frobenius series
\[
 \sum_{h=0}^n u^h e_h[X]h_{n-h}[X].
\]
The formal moment algebra contributes the reciprocal Euler factor
\[
 \frac{(q;q)_\infty}{\prod_{a\ge0}(1+uq^a)},\qquad
 \prod_{a\ge0}(1+uq^a)^{-1}
     =\sum_{j\ge0}\frac{(-u)^j}{(q;q)_j}.
\]
This proves the asserted formula for $\mathscr E_n$.

The difference $(q;q)_{n-h}-(q;q)_\infty$ begins in degree
$n-h+1$ with coefficient one. In fixed exterior degree $b$,
every term of the difference therefore has $q$-degree at
least $n-b+1$. At that first degree only $h=b$ contributes,
and the constant term of $e_b[X]h_{n-b}[X]$ is
$e_b[Z]h_{n-b}[Z]$. This proves both conclusions.
\end{proof}

\begin{theorem}[Row-hook moment acyclicity]
\label{thm:hook-acyclicity}
For every $n$, every row dimension $k$, every row hook
$\lambda\vdash D$, and every $D\le n$,
\[
 \Hom_{\GL_k}\!\left(
   \Schur_\lambda K^k,H_i(\mathcal K_{n,k,D})_D\right)=0
                  \qquad(i>0).
\]
Thus the full $S_n$ character in every row-hook sector is its
moment Euler character throughout $D\le n$.
\end{theorem}

\begin{proof}
Evaluate the polynomial row functors on $K^{1|1}$. The coefficient
algebra is
$P=K[x_1,\ldots,x_n]\otimes\Lambda(\theta_1,\ldots,\theta_n)$.
The free super moment algebra $S$ has even generators of
bidegree $(a,0)$ for $a\ge1$ and odd generators of bidegree
$(a,1)$ for $a\ge0$. They act by
$\sum_i x_i^a$ and $\sum_i\theta_i x_i^a$, respectively.
Its ordinary Hilbert series is
\[
 H_S(q,u)=\frac{\prod_{a\ge0}(1+uq^a)}{(q;q)_\infty}.
\]
Trace generation identifies $P/S_+P$ with the superspace
coinvariant ring. Lemma~\ref{lem:fields-euler}, after taking
dimensions, says exactly that
$H_P=H_SH_{P/S_+P}$ through total degree $n$.
Lemma~\ref{lem:truncated-freeness} forces positive Tor to vanish
there. The super Koszul resolution uses exterior powers on even
generators and symmetric powers on odd generators, so it is
precisely this evaluation of the moment complex.

To pass back to ordinary row spaces, express each homogeneous
polynomial functor as
$\bigoplus_\lambda\Schur_\lambda\otimes C_\lambda$.
In characteristic zero all maps act on the multiplicity spaces
$C_\lambda$, and this decomposition can be evaluated on a
super vector space with its signed tensor permutations.
It commutes with homology. A nonempty Schur functor is nonzero
on $K^{1|1}$ exactly when its partition is a hook, and
\[
 s_{(a,1^b)}[q\mid u]=(q+u)q^{a-1}u^b.
\]
Vanishing of the specialized homology therefore forces the
vanishing of every row-hook multiplicity, including its entire
$S_n$-module structure.
\end{proof}

\subsection{Exact Euler characters and character polynomials}
\label{sec:exact-euler-characters}

\begin{theorem}[Stable character polynomial formula]
\label{thm:explicit-stable-character-polynomials}
For $D\ge3$ and $n\ge\max\{D,2D-6\}$, the total-degree-$D$ part of
the moment Koszul complex is exact in positive degrees.
The same conclusion holds for $D=0,1$ and $n\ge1$, and
for $D=2$ and $n\ge2$.
Consequently, if $\mathbf b=(b_1,\ldots,b_k)$ has total
size $D$, then for every $\sigma\in S_n$,
\begin{equation}\label{eq:explicit-stable-character-polynomial}
 \operatorname{tr}\bigl(\sigma\mid(R_{n,k})_{\mathbf b}\bigr)
 = [\mathbf q^{\mathbf b}]
   \prod_{0<|\alpha|\le D}(1-\mathbf q^\alpha)
   \prod_{j=1}^{D}\prod_{a=1}^{k}
                 (1-q_a^j)^{-X_j(\sigma)},
\end{equation}
where $X_j(\sigma)$ is the number of $j$-cycles of $\sigma$.
This gives the character polynomial explicitly in every
multidegree and every number of alphabets.

For $3\le D\le n$, the same Koszul product computes every
label coefficient $S^\mu$ with $\mu_1\ge D-1$.
\end{theorem}

\begin{proof}
Suppose $H_j$ is nonzero on label $\mu$ in degree $D\le n$.
Quadratic reduction gives $\mu_1\ge n-D+3j$.
Theorem~\ref{thm:polynomial-schur-fronts} excludes the trivial label
and gives $\mu_2\ge\mu_{j+1}\ge n-D+j+1$.
Adding the first two rows therefore gives $n\le2D-4j-1$.
Thus $j\ge2$ is impossible when $n\ge2D-6$.
For $j=1$ the same inequality gives $n\le2D-5$.
This proves exactness for $n\ge\max\{D,2D-4\}$.

At $n=2D-5\ge D$, the only possible label is $(D-2,D-3)$.
Its first homology is a subquotient of
$\Schur_{(D-1)}(A_+)_D=\Sym^{D-2}V\otimes\Sym^2V$ by the
relation presentation. This has at most two rows and is removed by
Theorem~\ref{thm:bivariate-separation}.

At $n=2D-6\ge D$, put $b=D-3\ge3$. The bounds
$\mu_1\ge b$, $\mu_2\ge b-1$ leave only
\[
               (b+1,b-1),\qquad(b,b),\qquad(b,b-1,1).
\]
For $(b+1,b-1)$ the first relation source has degree-$D$ part
$\Sym^{b+1}V\otimes\Sym^2V$, again with at most two rows.
The other two labels first occur in coefficient degree $b$, with
multiplicities $\Sym^bV$ and $\Schur_{(b-1,1)}V$, respectively.
Positive invariant relations cannot change this first occurrence.
Thus their degree-$(b+3)$ first chain groups after quadratic reduction
are exactly $\Sym^3V$ tensored with these spaces. The first has only
two rows.

For the remaining label $(b,b-1,1)$, the relation source is
$\Schur_{(b+1,1)}(A_+)_{b+3}$, whose one-quadratic-input decomposition is
\[
       \bigl(\Schur_{(b,1)}V\oplus\Sym^{b+1}V\bigr)
                                      \otimes\Sym^2V.
\]
Pieri leaves only $(b+1,1,1)$ and $(b,2,1)$ after the two-row
constituents are removed. Row-hook acyclicity
(Theorem~\ref{thm:hook-acyclicity}) removes the first. The second occurs
once in the cubic source
$\Sym^3V\otimes\Schur_{(b-1,1)}V$, and its differential is nonzero
by Lemma~\ref{lem:cubic-witness}. It is therefore injective there.
When $\dim V<3$, both three-row components are absent.

For $D=3,4,5$, $n\ge D$ leaves only the earlier cases; the cubic
witness is first needed at $D=n=6$. Degrees zero and one use the
linear sequence. In degree two the remaining map is the inclusion
of the independent quadratic moments for $n\ge2$.
The alternating character of the exterior factor is
$\prod_{0<|\alpha|\le D}(1-\mathbf q^\alpha)$, and a $j$-cycle
has polynomial trace $\prod_a(1-q_a^j)^{-1}$. Their product gives
\eqref{eq:explicit-stable-character-polynomial}.
The labelwise assertion follows from Lemma~\ref{lem:parabolic-koszul-detection}.
\end{proof}

\begin{corollary}[Separation of invariants and harmonics]
\label{cor:uniform-separation}
Multiplication induces an equivariant isomorphism
\[
 (P_{n,k}^{S_n}\otimes\Harm_{n,k})_{\le D}
                         \xrightarrow{\ \sim\ }(P_{n,k})_{\le D}
\]
when $k\le2$ and $D\le n$, or for arbitrary $k$ when
$D\ge3$ and $n\ge\max\{D,2D-6\}$.
For arbitrary $k$, its degree-$D$ part is also an isomorphism on
every row-hook isotype and every isotype whose row partition has
length at most two, whenever $D\le n$.
\end{corollary}
\begin{proof}
The harmonic inclusion is an equivariant section of the coinvariant
quotient. Apply Lemma~\ref{lem:truncated-freeness} to the character
identity, using the formal moment algebra, which agrees with
$P_{n,k}^{S_n}$ through degree $n$. The bivariate case is
Theorem~\ref{thm:bivariate-separation}; all smaller degrees satisfy
the same bounds, with degrees zero, one and two immediate.

For the componentwise assertion, let $L$ be the kernel of the free
cover over the formal moment algebra $S$. Reduction modulo $S_+$
is an isomorphism, so $L/S_+L\cong\Tor_1^S(P_{n,k},K)$.
A Littlewood--Richardson constituent contains the partition of
each input factor. A kernel component of least degree among row
hooks and partitions with at most two rows cannot lie in $S_+L$:
its $L$ factor would have lower degree and belong to one of those classes.
It would therefore survive in $\Tor_1$, contrary to
Theorems~\ref{thm:bivariate-separation} and~\ref{thm:hook-acyclicity}
through degree $n$.
\end{proof}

The coefficient extraction in
\eqref{eq:explicit-stable-character-polynomial} is finite.
Expanding $(1-z)^{-X}$ as
$\sum_{v\ge0}\binom{X+v-1}{v}z^v$ shows directly that
it is a rational polynomial in the cycle counts, of weighted
degree at most $D$ when $\deg X_j=j$. At the identity it gives
\begin{equation}\label{eq:explicit-stable-dimension-polynomial}
 \dim(R_{n,k})_{\mathbf b}
 = [\mathbf q^{\mathbf b}]
    \prod_{0<|\alpha|\le D}(1-\mathbf q^\alpha)
    \prod_{a=1}^k(1-q_a)^{-n}.
\end{equation}
Its degree in $n$ is exactly $D$, with leading coefficient
$1/(b_1!\cdots b_k!)$.
For total degree $D$, the corresponding dimension is
\[
 [t^D](1-t)^{-k(n-1)}
       \prod_{j=2}^{D}(1-t^j)^{\binom{k+j-1}{j}},
\]
with leading term $k^Dn^D/D!$.
For $D\ge3$ and $n\ge2D-3$, Theorem~\ref{thm:linear-generation}
gives these same characters for the Vandermonde operator module.

This explicitly describes the eventual character and dimension
polynomials asked for in
\cite[Problems~1.10 and~1.12]{ChurchEllenbergFarb2015}.
The existence of eventual polynomiality and representation
stability is due to Church--Ellenberg--Farb; see also the
quantitative refinements of \cite{Bahran2024,Borisov2026}.
Outside this range, the same bounds locate every possible homological
correction.  The Euler product is the low-degree quotient proposed in
\cite[Section~4, equation~(33)]{Bergeron2013}.
Bergeron's conjecture asks for its validity throughout $D\le n$.

\subsection{Two-column homology, reflected cycles, and multiplication}
\label{sec:two-column-homology}

Put $W=W_n$ and $m=\dim W=n-1$. For $0\le q\le p$, let
$\lambda(p,q)=(2^q,1^{p-q})$, the partition with column heights $p,q$.
Its degree is $D=p+q$. Write $M_{p,q}$ for its stable row multiplicity
space in full diagonal harmonics, with $M_{0,0}=K$.
Only the invariant quadratic form survives in this sector; consequently
its multiplicity spaces and moment homology carry an $O(W)$ action
extending their $S_n$ action.

\begin{theorem}[Two-term reduction in every degree]
\label{thm:two-column-complex}
After eliminating the linear moments, the full moment Koszul complex
in row $\lambda(p,q)$ is
\begin{equation}\label{eq:two-column-two-term}
 \Schur_{\lambda(p-1,q-1)}W
                  \xrightarrow{\ L\ }
 \Schur_{\lambda(p,q)}W
\end{equation}
in homological degrees one and zero, with zero source when $q=0$.
The differential inserts the invariant symmetric form.
There is no $H_j$ for $j\ge2$ in any two-column row.
\end{theorem}
\begin{proof}
A moment of order at least three forces first row at least three
in every Littlewood--Richardson constituent. The same is true of
two or more quadratic Koszul generators, since
$\bigwedge^2\Sym^2V=\Schur_{(3,1)}V$ and
$\bigwedge^j\Sym^2V$ is a summand of
$\bigwedge^2\Sym^2V\otimes\bigwedge^{j-2}\Sym^2V$ for $j\ge2$.
Thus only zero or one quadratic generator survives. Cauchy
decomposition gives the target; horizontal Pieri removes one box
from each column to give the source once. The differential is
quadratic multiplication.
\end{proof}

Let $U=K^2$ record two exterior copies of $W$, and put
\[
 \mathcal E=\Lambda(U\otimes W),\qquad
 \omega=\sum_{i=1}^m\theta_i\wedge\phi_i,\qquad
 L=\omega\wedge-,\qquad
 c=\sum_{i=1}^m\partial_{\phi_i}\partial_{\theta_i}.
\]
The displayed orthonormal notation denotes the operators defined by
the form and its inverse, so no choice of an orthonormal basis over
$K$ is required. They satisfy
$[c,L]|_{\mathcal E_e}=(m-e)\operatorname{id}$. For
$P_e=\ker(c:\mathcal E_e\to\mathcal E_{e-2})$, $0\le e\le m$,
the exterior Lefschetz decomposition gives
\begin{align}
 \mathcal E&=\bigoplus_{e=0}^m\ \bigoplus_{j=0}^{m-e}L^jP_e,
 \label{eq:exterior-lefschetz-strings}\\
 cL^jf&=j(m-e-j+1)L^{j-1}f\qquad(f\in P_e).
 \label{eq:exterior-string-normalization}
\end{align}
These follow directly from finite-dimensional $\mathfrak{sl}_2$
theory; see also the exterior Lefschetz treatment in
\cite{KimRhoades}. The operator $L$ commutes with the row
$\mathfrak{sl}_2$ on $U$. On bidegree $(a,b)$, its formula in
$\Lambda^aW\otimes\Lambda^bW$ is
$(-1)^a\sum_i\eps_i\otimes\eps_i$. Exterior multiplication by
$\omega$ fixes these signs uniformly.

Write $\eta$ for the quadratic Koszul generator, of homological degree
one, internal degree two, and $GL(U)$ type $\det U$. The reduced
complex is $\eta\mathcal E\xrightarrow{\,L\,}\mathcal E$, with
$d(\eta f)=Lf$. For its algebra formulas we use the totalized
product, so $\eta$ anticommutes with $U\otimes W$; the complete
sign convention is given in Appendix~\ref{app:odd-derived}.

\begin{theorem}[All two-column homology and its reflection]
\label{thm:two-column-reflection}
For $0\le q\le p\le n$, set $(p^*,q^*)=(n-q,n-p)$.
\begin{enumerate}
\item If $p+q\le n-1$, the sole homology is
      $H_0=M_{p,q}=\Schur_{[\lambda(p,q)]}W$ as $O(W)$ modules.
\item If $p+q=n$, the complex is acyclic.
\item If $p+q\ge n+1$, the sole homology is
      \[
       H_1\cong M_{p^*,q^*}
          \cong M_{p^*,q^*}^{\,*}\otimes(\det W)^2
      \]
      as $O(W)$ modules, hence as $S_n$ modules.
      The reflected row has degree $2n-(p+q)$.
\end{enumerate}
If $p>n$, both terms vanish. These assertions hold in every row
dimension in which the original row occurs; in case~(3),
$p^*\le p-1$, so the reflected row occurs there as well.

More explicitly, if $f\in P_e$ and $k=m-e$, reflection is the cycle map
\begin{equation}\label{eq:exterior-reflected-cycle}
 \rho_e(f)=\frac{\eta L^kf}{k!},
 \qquad
 \rho_e^{-1}(\eta g)=\frac{c^kg}{k!}.
\end{equation}
Before extracting row multiplicity spaces, its source is
$P_e\otimes(\det U)^{k+1}$ and its target has internal degree $2n-e$.
\end{theorem}
\begin{proof}
Each string in \eqref{eq:exterior-lefschetz-strings} contributes its
initial vector to $\mathcal E/L\mathcal E$ and its terminal vector to
$\ker L$. Equation~\eqref{eq:exterior-string-normalization} gives
$L^{k+1}f=0$ and $c^kL^kf=(k!)^2f$, proving
\eqref{eq:exterior-reflected-cycle}.

The row-$\lambda(p,q)$ multiplicity space is the highest-weight kernel
of the transfer
$\Lambda^pW\otimes\Lambda^qW\to
\Lambda^{p+1}W\otimes\Lambda^{q-1}W$.
These kernels fix the Pieri normalization of
\eqref{eq:two-column-two-term}. Its low-degree primitive part is the
orthogonal traceless Schur module: within-column traces vanish by
antisymmetry, and the remaining trace is encoded by $c$.
This is also the admissible case of the orthogonal Littlewood theorem
\cite[Theorem~4.4]{SSW2013}.
If $f$ has $GL(U)$ highest weight $(a,b)$, then $e=a+b$ and
\eqref{eq:exterior-reflected-cycle} has weight
\[
       (a+k+1,b+k+1)=(n-b,n-a).
\]
Exterior Cauchy decomposition transposes this to the stated row
reflection. The strings leave no homology in internal degree $n$.
The dual formula follows from the exterior top-product pairing,
whose orientation after adjoining $\eta$ is
$(\det U)^n\otimes(\det W)^2$.
Finally, if $p>n$, both Schur terms have more than $\dim W$ rows.
\end{proof}

For the bivariate $S_n$-character $R_n(u,v)$ of the exterior $H_0$,
reflection and the Euler character give
\[
 R_n(u,v)-(uv)^nR_n(u^{-1},v^{-1})
    =(1-uv)\lambda_u(W_n)\lambda_v(W_n),
 \qquad
 \lambda_u(W_n)=\sum_j[\Lambda^jW_n]u^j.
\]
Products are tensor products in the representation ring. The dual
operation disappears because $S_n$-modules are self-dual.

\begin{corollary}[The full two-column first syzygies]
\label{cor:two-column-first-syzygies}
In degree $n+1$, for $1\le b\le\lfloor(n+1)/2\rfloor$, the row
$\lambda=(2^b,1^{n+1-2b})$ has
\[
 H_{1,\lambda}
 \cong\Schur_{[(2^{b-1},1^{n+1-2b})]}W_n\big|_{S_n}
 \cong\sgn\otimes\Schur_{[(2^{b-1})]}W_n\big|_{S_n}.
\]
In particular, for $n\ge4$,
\[
 H_{1,(2,2,1^{n-3})}
   \cong S^{(2,1^{n-2})}\oplus S^{(2,2,1^{n-4})}.
\]
\end{corollary}
\begin{proof}
At this boundary the reflected columns are obtained by subtracting
one from each original column. Orthogonal determinant duality
complements the first column in dimension $n-1$, giving the second
formula. For $b=2$ use
$\Schur_{[(2)]}W_n=\Sym^2W_n-K
 \cong S^{(n-1,1)}\oplus S^{(n-2,2)}$ for $n\ge4$.
\end{proof}

Thus at $n=3$ row $(2,2)$ reflects to $(2)$ and has
$H_1\cong S^{(2,1)}$; at $n=4$ row $(2,2,1)$ reflects to $(2,1)$
and has $H_1\cong S^{(2,1,1)}\oplus S^{(2,2)}$.
The exterior $H_0$ character is also obtained in
\cite{KimRhoades,Lentfer}, with explicit double-hook label formulas
in \cite{LentferSign}; the positive homology is the width-two case
of \cite[Theorem~4.4]{SSW2013}. The cycle map above will now identify
its multiplication as well.

\begin{proposition}[Multiplication on reflected homology]
\label{prop:two-column-module}
Put $R=\mathcal E/(\omega)$ and
$I=\eta\operatorname{Ann}_{\mathcal E}(\omega)$.
For $\beta\in U\otimes W$, set
$\iota_\beta=[c,\varepsilon_\beta]$. Then
\begin{equation}\label{eq:two-column-contraction-action}
 \beta\,\rho_e(f)=-\rho_{e-1}(\iota_\beta f),
 \qquad
 I=R\cdot(\eta P_m).
\end{equation}
In particular, the minimal $R$-generator space of $I$ is $\eta P_m$,
all in internal degree $n+1$, and
\[
 \dim(I/R_+I)
   =\binom{2m}{m}-\binom{2m}{m-2}
   =\frac{1}{n+1}\binom{2n}{n}.
\]
The last formula uses the usual zero convention for binomial
coefficients outside their range.
\end{proposition}
\begin{proof}
The identities $[L,\iota_\beta]=\varepsilon_\beta$ and
$[c,\iota_\beta]=0$ give
$\iota_\beta P_e\subseteq P_{e-1}$. With $k=m-e$,
\[
 L^{k+1}\iota_\beta f=(k+1)\beta L^kf
\]
because $L^{k+1}f=0$. Division by $(k+1)!$, followed by moving
$\beta$ across $\eta$, gives the first identity.
The contraction $U\otimes W\otimes P_e\to P_{e-1}$ is a nonzero
symplectic-equivariant map onto an irreducible primitive module.
Contracting $\theta_1\cdots\theta_e$ by $\phi_e$ proves nonvanishing.
It is therefore onto, and iteration from $e=m$ proves generation.
For $m=0$ the assertion is simply $I=K\eta$.
The generator quotient is its first nonzero degree; the dimension
formula is the primitive exterior-power formula.
\end{proof}

The Catalan formula of \cite[Corollary~4.4]{KimRhoades}
counts the minimal syzygy-generating layer after evaluation on two odd rows.
The identities are $Sp_{2m}$-equivariant. For $GL(U)\times O(W)$,
contraction contributes the extra $\det U$ needed between the source
twists of $\rho_e$ and $\rho_{e-1}$; on row multiplicity spaces they
are $O(W)$-equivariant.

\begin{corollary}[All two-column syzygies are generated at the first boundary]
\label{cor:two-column-syzygy-generation}
Write $H_{j,D}(V)$ for the degree-$D$ moment homology as a polynomial
functor of the ordinary row space $V$. For every $D\ge n+1$, the
cokernel of multiplication
\[
 H_{0,D-n-1}(V)\otimes H_{1,n+1}(V)
                         \longrightarrow H_{1,D}(V)
\]
has no constituent $\Schur_\lambda V$ with $\lambda_1\le2$.
This holds in every row dimension in which that constituent occurs.
\end{corollary}
\begin{proof}
Purely odd two-dimensional evaluation is exact in characteristic
zero and sends
$\Schur_\lambda(\Pi U)$ to
$\Pi^{|\lambda|}\Schur_{\lambda'}U$.
It kills precisely the shapes with $\lambda_1>2$ and detects every
remaining Schur multiplicity. It preserves the moment differential
and multiplication. The evaluated map is onto by
\eqref{eq:two-column-contraction-action}, so its original cokernel
has no surviving two-column shape.
\end{proof}

The quotient by shapes with more than two columns carries the induced
algebra structure, since these shapes form a tensor ideal by
Littlewood--Richardson containment.
Appendix~\ref{app:odd-derived} identifies the complete odd differential
algebra and shows that triple Massey products generate its positive
homology.

\begin{proposition}[A single-degree criterion for two-column generation]
\label{prop:two-column-middle-test}
For fixed $n\ge2$, Vandermonde generation in every two-column row and
every degree is equivalent to generation of the two-column rows in
degree $n-1$. It suffices to test the rows
\[
             (2^q,1^{n-1-2q}),
       \qquad 2\le q\le\lfloor(n-1)/2\rfloor.
\]
For $n\ge5$, these tests require only $n-3$ ordinary alphabets.
All two-column rows are generated for $n\le5$.
\end{proposition}
\begin{proof}
Let $G$ be the purely odd two-dimensional evaluation of the ordinary
operator subfunctor $\Pol$. This means evaluation of the entire
operator module, not of its Vandermonde seed alone.
The corresponding full harmonic space is
$\bigoplus_{e=0}^mP_e$.
Derivative closure is expressed by the natural polynomial coaction
\[
 \Pol_D(V)\longrightarrow
            (V\otimes W)\otimes\Pol_{D-1}(V).
\]
Exact odd evaluation therefore makes $G$ closed under every
constant-coefficient exterior contraction.
The map $(U\otimes W)^*\otimes P_e\to P_{e-1}$ is onto, by the
same nonzero symplectic-equivariant contraction used in
Proposition~\ref{prop:two-column-module}.
Thus $G_m=P_m$ propagates down to $G_e=P_e$ for every $e\le m$.
Odd evaluation detects every two-column shape, proving the
equivalence. The cases $q=0,1$ are row hooks and follow from
Proposition~\ref{prop:row-hook-generation}. Each remaining row has
$n-1-q\le n-3$ parts, so Lemma~\ref{lem:stability} and
Corollary~\ref{cor:r3-exact-specialization} make $n-3$ alphabets
detect every remaining multiplicity. The list is empty for $n\le4$;
for $n=5$ its only row is $(2,2)$, covered by Haiman's bivariate
operator theorem as in Theorem~\ref{thm:long-first-row}.
\end{proof}

The remaining task is to generate the displayed degree-$(n-1)$ rows.
The derivative coaction then propagates generation to every lower degree.

\subsection{The first row-hook moment syzygies}

\begin{theorem}[The first row-hook moment syzygies]
\label{thm:first-hook-syzygies}
In total degree $n+1$,
\[
 \left.H_1(\mathcal K_{n,k,n+1})_{n+1}
                 \right|_{\mathrm{row\ hooks}}
  \cong
  \bigoplus_{b=0}^{n-1}
      \Schur_{(n+1-b,1^b)}K^k\otimes S^{(n-b,1^b)}.
\]
Each summand occurs once. Row lengths exceeding $k$ give zero
Schur modules. Every $H_i$ with $i\ge2$ has zero row-hook part
in this total degree.
\end{theorem}

\begin{proof}
In the $(1|1)$ free cover from the proof of Theorem~\ref{thm:hook-acyclicity}, its kernel
begins in degree $n+1$. A positive term of the super Koszul
complex on that kernel has internal degree greater than
$n+1$. The Tor long exact sequence therefore makes every
$H_i$, $i\ge2$, zero in degree $n+1$.
Consequently the character discrepancy
\eqref{eq:first-fields-discrepancy} is exactly $H_1$.

Let $T_b$ be the label multiplicity space of the row
$(n+1-b,1^b)$ in this homology, and put $T_{-1}=0$.
The two monomials in its super Schur specialization give
\[
 \operatorname{ch}T_b+\operatorname{ch}T_{b-1}
       =e_b[Z]h_{n-b}[Z]\qquad(0\le b\le n).
\]
For $1\le b\le n-1$, Pieri's rule says
\[
 e_bh_{n-b}
    =s_{(n-b,1^b)}+s_{(n-b+1,1^{b-1})}.
\]
At the endpoints the right side is $h_n$ and $e_n$.
Induction gives $T_b\cong S^{(n-b,1^b)}$ for $0\le b\le n-1$,
and $T_n=0$. Exact Schur evaluation yields the stated
ordinary-row result and the higher-homology vanishing.
\end{proof}

For $n=1$, the theorem gives only $\Schur_{(2)}K^k$ with
trivial label in degree two. For $n=2$, it gives
$\Schur_{(3)}K^k\otimes S^{(2)}$ and
$\Schur_{(2,1)}K^k\otimes S^{(1,1)}$ in degree three.
These also follow by direct moment calculation.

\subsection{Superspace characters and the first correction}

Theorem~\ref{thm:first-unstable-specht-wall} identifies the first
unstable correction for every nonempty tail $\tau$ and every
$D>|\tau|$. Its idempotent specialization proves injectivity at
that wall directly, without an Euler-character comparison.

\begin{corollary}[Stable diagonal superspace characters]
\label{cor:stable-superspace-characters}
Let $H_+(X)=\sum_{d\ge1}h_d(X)$ and define
$\operatorname{sp}_{Q\mid U}(p_m)=p_m(Q)+(-1)^{m-1}p_m(U)$.
For $\tau\vdash s$, $D\ge s$ and $n\ge D+\tau_1$, the
degree-$D$ multiplicity character in any bosonic--fermionic superspace is
\[
 \operatorname{ch}\operatorname{Hom}_{S_n}
 (S^{\tau[n]},R_n^{(k,j)}[D])
 =\operatorname{sp}_{Q\mid U}\bigl([s_\tau[H_+(X)]]_D\bigr).
\]
The full character-polynomial formula of
Theorem~\ref{thm:explicit-stable-character-polynomials} likewise
transfers to arbitrary superspaces in its range
$D\ge3$, $n\ge\max\{D,2D-6\}$.
\end{corollary}
\begin{proof}
The bosonic identities hold in every row dimension and hence are
identities of universal row symmetric functions. Apply Lentfer's
diagonal supersymmetry theorem~\cite[Theorem~1.2 and Corollary~1.3]{Lentfer}.
The displayed specialization records parity; at module level the
corresponding functors are super Schur functors.
\end{proof}

\section{Disjoint-label products and first occurrences}
\label{sec:first}

The first occurrence of a label has a direct monomial description:
fill each row of its Young diagram with the least-weight one-label
monomials.  Each column then forms a lower monomial diagram.
We generate the alternants of these diagrams on disjoint label
sets and apply the Young symmetrizer.  This generates every
minimum-degree filling.

For $V=K^k$, define
\[
 b_j=\dim\Sym^jV=\binom{k+j-1}{k-1},\qquad
 N_j=\binom{k+j}{k},\qquad N_{-1}=0.
\]
Let $d_i=j$ for $N_{j-1}<i\leq N_j$, and pad all partitions with
zeros.  Thus $d_1,d_2,\ldots$ are the degrees of all monomials in $k$
variables, in nondecreasing order.

\begin{theorem}[Complete first occurrences]
\label{thm:r5-complete-first}
For every $n,k\geq1$ and every $\mu\vdash n$, put
\[
 \delta_\mu(k)=\sum_{i\geq1}\mu_i d_i.
\]
There is no $S^\mu$ in $P_{n,k}$ below this degree, and
\begin{equation}\label{eq:r5-first-equality}
 \Pol_{n,k}[\delta_\mu(k)]_\mu
 =\Harm_{n,k}[\delta_\mu(k)]_\mu
 =P_{n,k}[\delta_\mu(k)]_\mu.
\end{equation}
Its row multiplicity space is
\begin{equation}\label{eq:r5-first-character}
 \Hom_{S_n}\bigl(S^\mu,\Pol_{n,k}[\delta_\mu(k)]\bigr)
 \cong
 \bigotimes_{j\geq0}
 \Schur_{(\mu_{N_{j-1}+1},\ldots,\mu_{N_j})}(\Sym^j V).
\end{equation}
\end{theorem}

The degree and tensor-product formula come from the semistandard
construction of Schur functors; Romero records the ambient
plethystic character in \cite[Theorem~1.2 and \S2]{r5-Romero}.
For the alternating label, Wallach gives the determinant-twisted
exterior-power description of the first layer
\cite[Theorem~1]{r5-Wallach}.  The assertion here generates that
entire space.  Its proof rests on a product theorem for arbitrary
lower monomial diagrams.

\subsection{Curvilinear quotients and colored determinants}
\label{sec:r5-quot}

A finite nonempty set $\Gamma$ of monomials in $u_1,\ldots,u_k$ is a
\emph{lower diagram} if it is closed under divisors.
For labels $C$ with $|C|=|\Gamma|$, put
\[
 \Delta_\Gamma(X_C)=
       \det(X_i^\alpha)_{i\in C,\,u^\alpha\in\Gamma},
 \qquad
 X_i^\alpha=\prod_{a=1}^k(x_i^{(a)})^{\alpha_a}.
\]
An ordering of $C$ and $\Gamma$ fixes the irrelevant overall sign.

\begin{theorem}[Products of lower-diagram alternants]
\label{thm:r5-products}
Partition $\{1,\ldots,n\}$ into nonempty sets
$C_1,\ldots,C_r$.  For each $c$, choose any lower diagram
$\Gamma_c$ of $|C_c|$ monomials in $k$ variables.  Then
\begin{equation}\label{eq:r5-product}
       \prod_{c=1}^r\Delta_{\Gamma_c}(X_{C_c})
                         \ \in\ \Pol_{n,k}.
\end{equation}
The lower diagrams may be chosen independently.
The derivative and polarization module generated by this product is
also contained in $\Pol_{n,k}$.
\end{theorem}

The product theorem follows by placing all label blocks in one
quotient of a free module.  B\'erczi--Svendsen's theorem moves that
monomial quotient into the closure of the curvilinear locus, where
a jet determinant gives its operator realization.  We work over
$\C$ with the underlying reduced loci and then descend the resulting
linear-span statement to $K$.  The case $n=1$ is immediate.

For $S=\C[u_1,\ldots,u_k]$, the
punctual Quot scheme parametrizes length-$n$ quotients of $S^{\oplus r}$
supported at the origin.  Its curvilinear locus consists of quotients
whose underlying $S$-module is the structure module of an embedded
curvilinear length-$n$ subscheme.  In coordinates, such a quotient has
the form
\begin{equation}\label{eq:r5-curvilinear-quotient}
 S^{\oplus r}\longrightarrow A=\C[t]/(t^n),\qquad
 u_a\longmapsto f_a(t),\qquad e_c\longmapsto g_c(t),
\end{equation}
where the $f_a$ are pointed jets, at least one has nonzero linear
coefficient, and the $g_c$ generate $A$ as an $S$-module.  Equivalently,
at least one $g_c$ has nonzero constant term.

\begin{theorem}[B\'erczi--Svendsen]
\label{thm:r5-bs-input}
Every quotient
\[
               S/I_1\oplus\cdots\oplus S/I_r
\]
with each $I_c$ a monomial ideal of finite colength is contained in
the closure of the curvilinear locus in the length-$n$ punctual
Quot scheme, where $n=\sum_c\dim_\C S/I_c$.
\end{theorem}

This is \cite[Theorem~1.3, equivalently Theorem~5.2]{r5-BS}.
The same paper's Theorem~1.1 gives the $r=1$ Hilbert-scheme case.
The rank-$r$ Quot statement treats the direct sum of the independent
monomial quotients and is the input for \eqref{eq:r5-product}.

\subsubsection{The graph determinant}

For pointed jets $f=(f_1,\ldots,f_k)$ modulo $t^n$, define
\begin{equation}\label{eq:r5-jet-determinant}
 G_f(X)=\det\left([t^j]\exp\left(\sum_{a=1}^k
                                 x_i^{(a)}f_a(t)\right)\right)
                _{\substack{1\leq i\leq n\\0\leq j<n}}.
\end{equation}
This is the jet determinant of
Theorem~\ref{thm:r3-normal-form}, with the constant
$\prod_{j=0}^{n-1}j!$ omitted. On the graph chart
$f_1(t)=t$, $f_a(t)=\sum_s a_{a,s}t^s$, that theorem gives
\begin{equation}\label{eq:r5-graph-operator}
 \left(\prod_{j=0}^{n-1}j!\right)G_f
 =\exp\left(\sum_{a=2}^k\sum_{s=1}^{n-1}a_{a,s}Y_{a,s}\right)
       \Delta_n\in\Pol_{n,k},
 \qquad Y_{a,s}=\sum_i x_i^{(a)}\partial_{x_i}^{s}.
\end{equation}
Its inverse source reparametrization gives the same membership
throughout the dense regular-jet locus. We now retain this
membership under the colored specialization.

\subsubsection{The colored Pl\"ucker map}

To carry curvilinear closure into a linear span of polynomials,
we use one linear map on a fixed exterior power.  Its label colors
remember which lower diagram belongs to each block.  At a graph
quotient it becomes a differential operator on the jet determinant;
at a monomial quotient it becomes the desired product of alternants.

Every length-$n$ module supported at the origin is annihilated by
$\m^n$, where $\m=(u_1,\ldots,u_k)$.  To see this, the strict
nonzero steps in its $\m$-adic filtration each decrease dimension;
a repeated step would be zero by Nakayama's lemma.  All quotient
points therefore belong to the fixed Grassmannian of $n$-dimensional
quotients of
\[
                     W=(S/\m^n)^{\oplus r}.
\]
Let $c(i)=c$ for $i\in C_c$.  Define a linear map
\begin{equation}\label{eq:r5-colored-map}
 \begin{split}
 L_C:\bigwedge^nW^*&\longrightarrow P_{n,k},\\
 q_0\wedge\cdots\wedge q_{n-1}
   &\longmapsto
 \det\left(\sum_{|\alpha|<n}
        q_j(u^\alpha e_{c(i)})\frac{X_i^\alpha}{\alpha!}\right)_{i,j},
 \qquad \alpha!=\prod_a\alpha_a!.
 \end{split}
\end{equation}
Multilinearity and alternation in the $q_j$ make this a well-defined
linear map on the exterior power.  Its normalization does not depend
on the quotient point.  Its values have total degree at most
$n(n-1)$, so they lie in a fixed finite-dimensional polynomial
space $E=P_{n,k,\le n(n-1)}$.

For the graph quotient \eqref{eq:r5-curvilinear-quotient}, take
$q_j$ to be its coefficient-of-$t^j$ functional.  Then
\begin{align}
 L_C(q_0\wedge\cdots\wedge q_{n-1})
  &=\det\left([t^j]g_{c(i)}(t)
                \exp\left(\sum_a x_i^{(a)}f_a(t)\right)\right)_{i,j}
       \label{eq:r5-curv-image}\\
  &=\left(\prod_{i=1}^n g_{c(i)}(\partial_{x_i})\right)G_f
       \ \in\ \Pol_{n,k}.\nonumber
\end{align}
The second equality uses $f_1(t)=t$.  Each derivative affects just
one row of the determinant.  Constant terms of the $g_c$ are allowed
and act as scalar summands of these differential operators.

The graph chart is dense in the curvilinear Quot locus: the regular
jets form an irreducible open set, the condition that the $g_c$
generate is also open, and $f_1'(0)\ne0$ is a nonempty open condition.
Inverse reparametrization preserves the quotient and produces graph
form.  For $n=1$, the same assertion is immediate with all $f_a=0$.
Consequently the closed condition
\begin{equation}\label{eq:r5-closed-condition}
                        L_C(\ell)\subseteq\Pol_{n,k}
\end{equation}
holds on the entire curvilinear Quot component, where $\ell$ is the
Pl\"ucker line of a quotient.  More explicitly, compose $L_C$ with
the map $E\to E/(E\cap\Pol_{n,k})$.
Condition \eqref{eq:r5-closed-condition} says that this linear map
vanishes on $\ell$, so it is a closed linear condition on the
Grassmannian.  It includes the points where $L_C(\ell)=0$, and
therefore passes through every specialization in the closure.

\begin{proof}[Proof of Theorem~\ref{thm:r5-products}]
Let $I_c$ be the monomial ideal whose standard monomials are
$\Gamma_c$.  The monomial quotient $\bigoplus_c S/I_c$ has total
length $n$, so Theorem~\ref{thm:r5-bs-input} and
\eqref{eq:r5-closed-condition} apply.
Choose quotient functionals dual to the colored basis
$\{u^\alpha e_c:u^\alpha\in\Gamma_c\}$.
In \eqref{eq:r5-colored-map}, the corresponding column is zero
outside the label block $C_c$.  Thus, up to the ordering sign,
\[
 L_C(\ell)=
 K\left(
   \frac{\prod_{c=1}^r\Delta_{\Gamma_c}(X_{C_c})}
        {\prod_{c=1}^r\prod_{u^\alpha\in\Gamma_c}\alpha!}
     \right).
\]
This is a nonzero line.  It is contained in $\Pol_{n,k}$ by
\eqref{eq:r5-closed-condition}, proving \eqref{eq:r5-product} over
$\C$.  The product, the defining operators, and the Vandermonde all
have rational coefficients.  The assertion is membership in a
finite-dimensional rational linear span, so membership over $\C$
implies membership over $\Q$ and then over every characteristic-zero
field $K$.  Closure under derivatives and polarizations gives the
last assertion.
\end{proof}

Each factor in \eqref{eq:r5-product} is harmonic for its own label
set: a positive diagonal differential operator replaces a monomial
column by a divisor column already present, or by zero.  Since the
label sets are disjoint, the product is harmonic for all $n$ labels.
Theorem~\ref{thm:r5-products} places these products in the operator module.

\subsection{Multiplication on disjoint label sets}

The product theorem extends from lower-diagram alternants to the
whole family of operator modules.  Private row spaces let the
operators act independently on each label block; an invertible row
shear and exact specialization then identify all row spaces.

\begin{theorem}[Disjoint-label multiplication]
\label{thm:r5-disjoint-products}
Let $C_1,\ldots,C_r$ partition $\{1,\ldots,n\}$, with
$|C_c|=n_c>0$.  For arbitrary
$F_c\in\Pol_{n_c,k}$, relabelled on $C_c$, one has
\[
                         \prod_{c=1}^r F_c(X_{C_c})\in\Pol_{n,k}.
\]
These products are associative, commutative, and equivariant under
relabeling and the diagonal row action.  They define canonical
$S_n\times\GL_k$-equivariant maps
\begin{equation}\label{eq:r5-induced-product}
 \operatorname{Ind}_{S_{n_1}\times\cdots\times S_{n_r}}^{S_n}
       \left(\bigotimes_{c=1}^r\Pol_{n_c,k}\right)
                         \longrightarrow\Pol_{n,k}.
\end{equation}
\end{theorem}

\begin{proof}
Introduce $r$ private blocks of $k$ variable rows, with variables
$x_i^{(c,a)}$.  For the $c$th factor take the lower diagram
$\Gamma_c=\{1,u_{c,1},\ldots,u_{c,1}^{n_c-1}\}$ in these $kr$
one-label variables.  Theorem~\ref{thm:r5-products} gives
\[
                    \prod_c\Delta_{C_c}(x^{(c,1)})
                                      \in\Pol_{n,kr}.
\]
Generate $F_c$ using derivatives and polarizations whose source
and target rows both lie in its private block.  Such a global
polarization acts only on the $c$th factor: every other factor is
independent of its source row, and its derivative order is at least
one.  Label-specific partial derivatives have the same property.
By linearity, this produces
$\prod_c F_c(X_{C_c}^{(c)})$ inside $\Pol_{n,kr}$.

Apply the invertible row shear
\[
 x_i^{(c,a)}\longmapsto x_i^{(c,a)}+x_i^{(1,a)}\quad(c>1),
 \qquad x_i^{(1,a)}\longmapsto x_i^{(1,a)}.
\]
Then set all private rows with $c>1$ to zero.
The shear preserves $\Pol_{n,kr}$, and
Corollary~\ref{cor:r3-exact-specialization} sends the result into
$\Pol_{n,k}$.  The resulting polynomial is the desired product
in the common $k$ rows.  The associativity, commutativity, and
equivariance assertions follow from ordinary multiplication in
disjoint polynomial variables.  Frobenius reciprocity gives
\eqref{eq:r5-induced-product}.
\end{proof}

\begin{corollary}[Tensor products under restriction]
\label{cor:r5-product-restriction}
There is a natural injective map of
$\GL_k\times S_{n_1}\times\cdots\times S_{n_r}$-modules
\[
 \bigotimes_c\Pol_{n_c,k}
   \ \lhook\joinrel\longrightarrow\
 \operatorname{Res}_{S_{n_1}\times\cdots\times S_{n_r}}^{S_n}
                                     \Pol_{n,k}.
\]
Consequently, for label partitions $\mu_c\vdash n_c$ and a row
partition $\lambda$ of length at most $k$,
\begin{equation}\label{eq:r5-product-lr}
 \sum_{\nu\vdash n}
       \operatorname{LR}_{\mu_1,\ldots,\mu_r}^{\nu}
                         c_{\lambda,\nu}^{\Pol_{n,k}}
 \ \geq\
 \sum_{\lambda_1,\ldots,\lambda_r}
       \operatorname{LR}_{\lambda_1,\ldots,\lambda_r}^{\lambda}
                     \prod_{c=1}^r
                       c_{\lambda_c,\mu_c}^{\Pol_{n_c,k}}.
\end{equation}
Here $\operatorname{LR}$ denotes an iterated Littlewood--Richardson
coefficient, and the sum on the right is over row partitions of
length at most $k$.
\end{corollary}

\begin{proof}
Multiplication identifies the tensor product of polynomial rings on
the disjoint label sets with the full polynomial ring.  Its
restriction to the tensor product of the indicated subspaces is
therefore injective, and Theorem~\ref{thm:r5-disjoint-products}
places its image in $\Pol_{n,k}$.  Decompose this inclusion under
the row group and the label Young subgroup.  The
Littlewood--Richardson restriction and tensor product rules give
\eqref{eq:r5-product-lr}.
\end{proof}

For fixed label blocks, multiplication is injective as stated in
Corollary~\ref{cor:r5-product-restriction}.  The induced map
\eqref{eq:r5-induced-product} combines all relabelings of those
blocks and can have a kernel.

\subsection{Every positive-weight first layer}
\label{sec:r5-weighted}

The first-occurrence construction works for every positive row
weight.  A minimum-weight column contains all cheaper monomials,
so it remains a lower diagram even when many monomials have equal
weight.  Give the $a$th row a positive integer weight $w_a$.  Write
$P_{n,k}^{(w)}[d]$ for weighted degree $d$ and put
\[
 \begin{aligned}
 W_d(w)&=\Span_K\{u^\alpha:w\cdot\alpha=d\},
 &b_d(w)&=\dim W_d(w),\\
 N_d(w)&=\sum_{e=0}^d b_e(w),
 &N_{-1}(w)&=0.
 \end{aligned}
\]
Order all monomials by nondecreasing weight, breaking ties in any
fixed way.  Let $d_i(w)$ be the weight of the $i$th monomial and set
\[
 \delta_\mu(w)=\sum_i\mu_i d_i(w),\qquad
 \mu^{(d,w)}=
  (\mu_{N_{d-1}(w)+1},\ldots,\mu_{N_d(w)}).
\]
Empty blocks contribute the trivial factor.

\begin{theorem}[Complete generation on every lower weighted face]
\label{thm:r5-weighted-first}
For every positive integral weight vector $w$ and $\mu\vdash n$,
$\delta_\mu(w)$ is the first weighted degree containing $S^\mu$, and
\begin{equation}\label{eq:r5-weighted-equality}
 \Pol_{n,k}^{(w)}[\delta_\mu(w)]_\mu
 =\Harm_{n,k}^{(w)}[\delta_\mu(w)]_\mu
 =P_{n,k}^{(w)}[\delta_\mu(w)]_\mu.
\end{equation}
The multiplicity space of this layer is
\begin{equation}\label{eq:r5-weighted-character}
                  \bigotimes_{d\geq0}
                         \Schur_{\mu^{(d,w)}}(W_d(w)).
\end{equation}
This is an isomorphism for the row torus and for the subgroup of
$\GL_k$ preserving the weighted grading.
\end{theorem}

\begin{proof}
As an $S_n$-module with row action,
$P_{n,k}=(\Sym V)^{\otimes n}$.  Schur--Weyl duality identifies
the $S^\mu$ multiplicity space with $\Schur_\mu(\Sym V)$.
In its semistandard tableau construction, an entry in row $i$ is
at least the $i$th element of the ordered monomial basis.  Its
weight is therefore at least $d_i(w)$.  Filling row $i$ with that
$i$th monomial attains the bound $\delta_\mu(w)$.

Equality in the weight bound requires every box in row $i$ to
contain a monomial of weight exactly $d_i(w)$.  Inside a weight
block the fillings are the semistandard fillings for the partition
$\mu^{(d,w)}$ on $W_d(w)$.  Between different blocks, column
inequalities are automatic.  This proves
\eqref{eq:r5-weighted-character} and the claimed first degree.
Every positive $D_\alpha$ is $S_n$-equivariant and lowers weighted
degree by $w\cdot\alpha>0$.  It annihilates this whole first
isotype, proving equality of the $P$ and $\Harm$ terms in
\eqref{eq:r5-weighted-equality}.

It remains to generate the entire space.  Fix a standard label
tableau $T$ of shape $\mu$, with label sets $C_c$ in its columns.
Let $F$ be any semistandard monomial filling attaining
$\delta_\mu(w)$, and let $\Gamma_c$ be the monomial set in its
$c$th column.  Column strictness makes these monomials distinct.
For each weight strictly below the last weight in that column,
the column crosses the complete corresponding block of rows.
It must therefore contain every monomial of that weight; in the
last block it contains a subset.  Since all $w_a$ are positive,
a proper monomial divisor has strictly smaller weight.  Thus
every $\Gamma_c$ is a lower diagram.

Let $m_F$ be the monomial in $P_{n,k}$ obtained by putting the
monomial in each box of $F$ at the label in the same box of $T$.
If $b_T$ is the column antisymmetrizer, then
\[
                     b_Tm_F
                       =\prod_c\Delta_{\Gamma_c}(X_{C_c})
\]
up to the fixed column ordering signs.  Theorem~\ref{thm:r5-products}
places this polynomial in $\Pol_{n,k}$.  Applying the row
symmetrizer $a_T$ of the label tableau also preserves $\Pol_{n,k}$.
The standard
semistandard construction says that the vectors $a_Tb_Tm_F$, as
$F$ ranges over these minimum-weight fillings, span the
first-weight part of the image of the Young symmetrizer $a_Tb_T$.
That image contains one fixed Specht line tensored with the whole
multiplicity space \eqref{eq:r5-weighted-character}.
Its $S_n$-span is the complete $\mu$-isotype.  Since $\Pol_{n,k}$
is $S_n$-stable, this proves the remaining containment in
\eqref{eq:r5-weighted-equality}.
\end{proof}

\begin{proof}[Proof of Theorem~\ref{thm:r5-complete-first}]
Take $w=(1,\ldots,1)$ in Theorem~\ref{thm:r5-weighted-first}.
Then $W_j(w)=\Sym^jV$, $b_j(w)=b_j$, and $N_j(w)=N_j$.
\end{proof}

The proof generates monomial fillings one by one before assembling
their row representations.  It therefore carries every irreducible
constituent and every plethystic multiplicity in
\eqref{eq:r5-weighted-character} into the operator module.

\begin{corollary}[The whole positive-weight lower boundary]
\label{cor:r5-lower-boundary}
For $\mathcal M=\Pol,\Harm$, let
$\operatorname{ch}_\mu^{\mathcal M}(\bq)$ be the row-torus character
of the $S^\mu$ multiplicity space in $\mathcal M_{n,k}$.  For every
$w\in\N_{>0}^k$, the terms minimizing $w\cdot\beta$ in the two
characters are identical, including their coefficients.  They are
\begin{equation}\label{eq:r5-face-character}
 \prod_{d\geq0}
 s_{\mu^{(d,w)}}
     \left[\sum_{w\cdot\alpha=d}q^\alpha\right],
\end{equation}
and their common minimum is $\delta_\mu(w)$.
Consequently the generation defect has no component on any face
of the $\mu$-character support exposed by a strictly positive
weight vector.
\end{corollary}

\begin{proof}
Theorem~\ref{thm:r5-weighted-first} is equality of the corresponding
row-torus modules.  The character of a Schur functor of $W_d(w)$
is the plethystic factor in \eqref{eq:r5-face-character}.
Every face exposed by a strictly positive real vector is also
exposed by a strictly positive rational vector: the relevant
relatively open normal cone is defined by rational linear
equalities and strict inequalities.  Clearing denominators
reduces that face to the integral-weight case.
\end{proof}

Strict positivity makes every proper monomial divisor cheaper and
every positive invariant differential operator degree lowering.
These are the two properties that drive the weighted proof.

\section{Constructive generation in a linear range}
\label{sec:linear-generation}

The selected-label quotient now becomes an evaluation map for harmonic
polynomials.  Local generated lifts can be multiplied on disjoint label
sets, and a two-label kernel saves one auxiliary label.  These two
operations prove the long-row and odd-row criteria.  The hook identities
then reduce the uniform range to one explicit rank-five lift.

\begin{lemma}[Selected evaluation]\label{lem:selected-evaluation}
Let $V=K^k$, $B=\operatorname{Sym}V$, $\mathfrak m=B_+$, $\mu=(t,\tau)$,
$|\tau|=s$, and $n=t+s$.
If $D\ge s$ and $n\ge D+\tau_1$, setting the auxiliary variables to zero
identifies $\Harm_{n,k}[D]_\mu^{S_t}$ with
$(\mathfrak m^{\otimes s})_{D,\tau}$.
More generally, in any range, suppose the span of selected
restrictions of $\Pol_{n,k}[D]^{S_t}$ contains the whole
row-$\lambda$, selected-label-$\tau$ component of
$(\mathfrak m^{\otimes s})_D$.  Then the entire row-$\lambda$,
label-$\mu$ harmonic component is generated.
\end{lemma}
\begin{proof}
Let $A_s=B^{\otimes s}=\mathfrak m^{\otimes s}\oplus W$, where $W$
contains tensors with a unit position. The minimum-occupancy argument
of Lemma~\ref{lem:minimum-occupancy-cover} gives the selected surjection
$\pi:A_s\to R_{n,k}^{S_t}$, eliminates $W$ on label $\mu$, and
identifies $(S^\mu)^{S_t}$ with $S^\tau$.
For a selected polynomial $f$ and an $S_t$-invariant harmonic $g$,
\[
 \langle\pi(f),g\rangle
   =\langle f,g(X_1,\ldots,X_s,0,\ldots,0)\rangle.
\]
Restriction is therefore the Fischer adjoint of this cover and is
injective. The occupancy summands are orthogonal, so on label $\mu$
it takes values in the fully occupied summand. In the stable range
Theorem~\ref{thm:stable-augmented-multiplicities} makes it an isomorphism.

For the general assertion, suppose a coinvariant class in the
indicated global component annihilates all generated harmonics.
Take a nonzero $S_t$-fixed vector in this annihilator and choose
a selected representative $f$, projected to row type $\lambda$
and selected-label type $\tau$.  Pairing against the prescribed
restrictions shows that the fully occupied part of $f$ is zero.
Hence $f\in W$.  But $\pi(W)$ has no global label $\mu$, a
contradiction.  Nondegeneracy of the coinvariant--harmonic pairing
proves generation.
\end{proof}

We will build prescribed restrictions by multiplying harmonics on
disjoint sets of labels.  Each local lift uses auxiliary labels whose variables are later set to zero.

\begin{theorem}[Long first rows]
\label{thm:long-first-row}
For every degree $D\ge1$ and every label $\mu\vdash n$,
\[
 \mu_1\ge D-1
 \quad\Longrightarrow\quad
 \Pol_{n,k}[D]_\mu=\Harm_{n,k}[D]_\mu
\]
in every number $k$ of alphabets.
\end{theorem}
\begin{proof}
For $D\le2$, every row partition has length at most two.
Lemma~\ref{lem:stability} reduces its highest-weight vectors to two
alphabets, where Haiman's operator theorem applies. The same reduction
will generate the two-row local lifts below. Assume $D\ge3$ and put
$t=D-1$. If $n<t$
there is no label in the stated range; if $n=t$, only the trivial
label qualifies, and its positive-degree harmonic part is zero.
We may therefore take $s=n-t\ge1$ selected labels.

The first selected relations have degree $t+1=D$.  In that degree
all their inputs are linear, so their row type is $(D)$.
The selected quotient, and hence its dual harmonic restriction,
is an isomorphism on every other row type.  We will generate all
of these row components with exactly $D-1$ auxiliary labels.

Fix an occupied selected multidegree
$d_1+\cdots+d_r=D$, with each $d_i>0$, and write
\[
 T=\bigotimes_{i=1}^r\Sym^{d_i}V,\qquad
 K_{ij}=\ker\bigl(\Sym^{d_i}V\otimes\Sym^{d_j}V
                         \longrightarrow\Sym^{d_i+d_j}V\bigr).
\]
The multiplication map gives
\begin{equation}\label{eq:pair-kernels-symmetric-quotient}
 T\Big/\sum_{i<j}
       K_{ij}\otimes\bigotimes_{h\ne i,j}\Sym^{d_h}V
                    \ \cong\ \Sym^D V.
\end{equation}
Indeed, realize $T$ as the quotient of $V^{\otimes D}$ by
permutations within its blocks of sizes $d_i$.  Quotienting by
$K_{ij}$ imposes symmetry on the union of the $i$th and $j$th
blocks.  These block-pair permutations generate all of $S_D$.
Thus the pair kernels span every row component of $T$ other than
$(D)$.

Put $d=d_i+d_j$.  By Pieri, $K_{ij}$ has at most two rows and
contains no copy of $(d)$.  At local rank $d+1$, take two selected
labels and $d-1$ auxiliary labels.  Its first selected relation has degree
$d$ and row type $(d)$, so restriction has a harmonic lift of
every vector in $K_{ij}$.  Haiman's theorem generates these lifts.
Each remaining selected factor $\Sym^{d_h}V$ has a generated lift
with $d_h$ auxiliary labels: in one row, differentiate the Vandermonde on
$x,y_0,\ldots,y_{d_h-1}$ by
$\prod_{q=0}^{d_h-1}\partial_{y_q}^{q}$ and set the auxiliary labels to
zero.  The result is a nonzero multiple of $x^{d_h}$; row
polarization gives the whole symmetric power.

The total number of auxiliary labels is
\[
                  (d-1)+\sum_{h\ne i,j}d_h=D-1.
\]
Multiply these lifts on disjoint label sets using
Theorem~\ref{thm:r5-disjoint-products}, and average over the
auxiliary labels; unused selected labels contribute constants.
Restriction sets every auxiliary label to zero, so averaging
preserves the prescribed selected polynomial.  Equation
\eqref{eq:pair-kernels-symmetric-quotient}, summed over selected
multidegrees, now generates every non-one-row harmonic component
fixed by $S_t$.  The ordinary operator theorem and row polarization
generate the one-row component.

Injectivity of harmonic restriction gives
$\Pol_{n,k}[D]^{S_t}=\Harm_{n,k}[D]^{S_t}$.
Every $S^\mu$ with $\mu_1\ge t$ has a nonzero $S_t$-fixed vector
by branching.  Exactness of invariants therefore gives the stated
equality on every such label.
\end{proof}

Alternating pairs save one auxiliary label each.  Placing such pairs in the
columns of a label tail makes the savings compatible with its
Young symmetrizer.

\begin{theorem}[Odd-row generation criterion]
\label{thm:odd-row-generation}
For every $\mu\vdash n$,
\[
 D\le\mu_1+\mu_3+\mu_5+\cdots
 \quad\Longrightarrow\quad
 \Pol_{n,k}[D]_\mu=\Harm_{n,k}[D]_\mu .
\]
\end{theorem}
\begin{proof}
Write $\mu=(t,\tau)$ and pair adjacent boxes vertically in every
column of $\tau$. The number of pairs is
$p=\tau_2+\tau_4+\cdots$.
An alternating pair of positive-degree polynomials of total degree $d$
has a generated lift on $d+1$ labels, with two selected labels.
Indeed, its selected-sign component has no degree-$d$ relation:
the first selected relations have degree $d$ and are $S_2$-invariant.
Its row types have length at most two, so Haiman's operator theorem
generates the harmonic lift. Such a pair costs $d-1$ auxiliary labels;
a singleton of degree $d$ costs $d$.

The disjoint products therefore use $D-p$ auxiliary labels. Applying the
Young symmetrizer for $\tau$ spans its whole fully occupied selected
component: each prealternated pair multiplies the column
antisymmetrization by two. Averaging over the auxiliary labels and Lemma~\ref{lem:selected-evaluation} prove generation whenever $t\ge D-p$, as asserted.
\end{proof}

\Needspace{9\baselineskip}
\begin{lemma}[Four universal tail degrees]
\label{lem:four-tail-degrees}
The alternating component is generated in every degree $d\le n+2$.
Every hook label of tail $s$ is generated in degrees $d\le s+3$.
More generally, if $\mu=\tau[n]$, $s=|\tau|$, then
\[
 0\le D-s\le3,\qquad n\ge D+\tau_1
 \quad\Longrightarrow\quad
 \Pol_{n,k}[D]_\mu=\Harm_{n,k}[D]_\mu .
\]
\end{lemma}
\begin{proof}
The ambient sign multiplicity is
$\bigwedge^n(\operatorname{Sym}V)$.
Through degree $n+2$, its nonlinear factors, after removing exterior
powers of $V$ and the possible constant polynomial, are among
\[
 \operatorname{Sym}^2V,\quad \operatorname{Sym}^3V,\quad
 \bigwedge^2\operatorname{Sym}^2V,\quad
 \operatorname{Sym}^4V,\quad
 \operatorname{Sym}^3V\otimes\operatorname{Sym}^2V,\quad
 \bigwedge^3\operatorname{Sym}^2V .
\]
The last two decompose as
$\Schur_5V\oplus \Schur_{41}V\oplus \Schur_{32}V$ and
$\Schur_{411}V\oplus \Schur_{33}V$, respectively;
$\bigwedge^2\operatorname{Sym}^2V=\Schur_{31}V$.
Vertical Pieri shows that none of these, after tensoring with an
exterior power of $V$, has both $\lambda_2\ge3$ and $\lambda_3\ge2$.
Corollary~\ref{cor:r2-defect-support} therefore rules out a sign
defect. The exact hook-defect identities give the hook assertion.

Now consider a column of $\tau$ of height $c$ and degree $d$.
Since the total excess is at most three, $c\le d\le c+3$.
The stable evaluation statement at local rank $d+1$ identifies
its fully occupied alternating component with the hook
$(d+1-c,1^c)$, which is generated by the assertion just proved.
It costs $d+1-c$ auxiliary labels. Summing over columns costs
$D-s+\tau_1$. Disjoint products, Young symmetrization and Lemma~\ref{lem:selected-evaluation} prove the final claim.
\end{proof}

\begin{theorem}[Uniform linear generation]
\label{thm:linear-generation}
For $D\ge3$ and $n\ge2D-3$, in every number of alphabets,
\[
 \Pol_{n,k}[d]=\Harm_{n,k}[d]\qquad(0\le d\le D).
\]
\end{theorem}
\begin{proof}
First prove the assertion in degree $D$.
Write $\mu=(t,\tau)$, $s=|\tau|$. Occupancy gives $s\le D$.
If $t\ge D-1$, apply Theorem~\ref{thm:long-first-row}.
Otherwise $t\le D-2$, and $n\ge2D-3$ forces $s\ge D-1$.
For $s=D$, Lemma~\ref{lem:four-tail-degrees} applies because
$n=D+t\ge D+\tau_1$.  For $s=D-1$, necessarily
$t=D-2$ and $n=2D-3$.  The same lemma applies if
$\tau_1\le D-3$.  The sole remaining label is therefore
\[
                         \mu=(D-2,D-2,1).
\]

Put $s=D-1$.  We reach its one exceptional row by combining a
three-selected-label harmonic at rank five with linear lifts.
The general part of Lemma~\ref{lem:selected-evaluation} lets us
work in the fully occupied selected space: spanning that space
by restrictions suffices at this rank as well.
The degree-$s+1$ part of $\Schur_{(s-1,1)}(\mathfrak m)$ has one
quadratic slot. Its only row type of length greater than two which
is not a row hook is $(s-2,2,1)$, with multiplicity one; it occurs
only for $s\ge4$. This follows by removing one corner of $(s-1,1)$
and applying the horizontal Pieri rule with $S_2V$.
For $s<4$, Haiman's theorem and row-hook generation therefore
finish the proof.  Assume $s\ge4$ and $\dim V\ge3$, so the
exceptional row is present.

At rank five, a local selected-three vector of selected type
$(2,1)$ and row type $(2,1,1)$ has a harmonic lift: the degree-four
relations have at most two rows.  Indeed, they come either from
a degree-three relation multiplied by a linear form, with row
source $\Sym^3V\otimes V$, or from a relation on two linear
inputs and one quadratic input, with row source
$\Sym^2V\otimes\Sym^2V$.  Neither source contains $(2,1,1)$.
Proposition~\ref{prop:row-hook-generation} generates this lift.
Multiply it with $s-3$ linear lifts. The number of auxiliary labels is
$2+(s-3)=s-1=D-2$.

The local row projection and selected Young projection must both
survive. In the realization
$\Schur_{(s-1,1)}\mathfrak m
=\ker(\operatorname{Sym}^{s-1}\mathfrak m\otimes\mathfrak m
\to\operatorname{Sym}^s\mathfrak m)$, the selected projection is
\[
 h\otimes\bigl(w\otimes(a\wedge b)\bigr)
             \longmapsto (hwa)\otimes b-(hwb)\otimes a.
\]
It symmetrizes the first row and preserves the distinguished last box.
The quadratic factor $w\in\Sym^2V$ varies independently, so
one can check its $x^2$ component separately.
Set
$W=\operatorname{Sym}^2V$ and choose
independent $x,y,z\in V$. In $W\otimes\bigwedge^2V$ put
\[
 v=x^2\otimes(y\wedge z)-xy\otimes(x\wedge z)
                         +xz\otimes(x\wedge y).
\]
This is a highest vector of type $(2,1,1)$, and
\[
 E_{21}v=xy\otimes(y\wedge z)-y^2\otimes(x\wedge z)
                         +yz\otimes(x\wedge y).
\]
Thus
\[
 x^{s-4}y\otimes v-x^{s-3}\otimes E_{21}v
\]
is highest of type $(s-2,2,1)$.
The image of the displayed highest vector has nonzero $w=x^2$
component
\[
 x^2\otimes
 (x^{s-4}y^2\otimes z-x^{s-4}yz\otimes y).
\]
The second summand has no such quadratic factor, so it cannot cancel.
The target is $\operatorname{Sym}^2V\otimes\operatorname{Sym}^{s-2}V\otimes V$;
its first factor records the quadratic factor.
Thus the selected Young projection is nonzero on the generated
products.  Row polarization reaches the whole multiplicity-one
component, and averaging over auxiliary labels preserves its selected restriction.
Lemma~\ref{lem:selected-evaluation} gives generation of this last
row-and-label component.
All remaining rows are generated by Haiman or row-hook generation.

This proves degree $D$.  For $3\le d\le D$ the same hypothesis
gives $n\ge2d-3$; degrees zero, one and two follow from the
ordinary and bivariate operator theorems.
\end{proof}

The labelwise criteria also locate every possible remaining defect.
If $(\Harm_{n,k}[D]/\Pol_{n,k}[D])_\mu\ne0$, then
\[
 n-D+1\le\mu_1\le D-2,
 \qquad \mu_1+\mu_3+\mu_5+\cdots\le D-1.
\]
Indeed, occupancy gives $n-\mu_1\le D$. Equality would make degree
$D$ a first occurrence, which is generated by
Theorem~\ref{thm:r5-complete-first}; the other bounds are
Theorems~\ref{thm:long-first-row} and~\ref{thm:odd-row-generation}.

\appendix
\section{The cubic multiplication at the character boundary}
\label{app:cubic-witness}

The proof of Theorem~\ref{thm:explicit-stable-character-polynomials}
leaves one possible multiplicity-one cubic source. We exhibit a
highest vector whose image survives the quadratic quotient. The
dual-number specialization below makes all unwanted moments vanish
and leaves only one term of that image.

\begin{lemma}[A nonzero cubic moment multiplication]
\label{lem:cubic-witness}
For $a\ge3$, put $n=2a$, $D=a+3$, and
\[
 \tau=(a-1,1),\qquad\mu=(a,a-1,1),\qquad\lambda=(a,2,1).
\]
For $\dim V\ge3$, the cubic multiplication map
\[
 \operatorname{Sym}^3V\otimes
       \operatorname{Hom}_{S_n}(S^\mu,Q_{n,k}[a])
      \longrightarrow\operatorname{Hom}_{S_n}(S^\mu,Q_{n,k}[a+3])
\]
is nonzero, and hence injective, on its multiplicity-one
$\mathbb S_\lambda V$ summand. Here $Q_{n,k}$ denotes the quotient
by linear and quadratic moments.
\end{lemma}
\begin{proof}
It is enough to work in three rows $x,y,z$, because the map is a
natural transformation of polynomial functors. Set $m=a-1$ and
choose disjoint label blocks of sizes $3,2,\ldots,2,1$, with $a-2$
blocks of size two. Define
\[
 F=\det\begin{pmatrix}1&x_0&y_0\\1&x_1&y_1\\1&x_2&y_2\end{pmatrix}
       \prod_{j=1}^{a-2}(x_{2j+1}-x_{2j+2}).
\]
The unused label is $2a-1$. This polynomial has row highest weight
$(m,1)$, since $\sum_i x_i\partial_{y_i}$ kills its determinant.
It belongs to label $\mu$. To justify this last assertion directly,
its column antisymmetry forces every label constituent to be
 dominated by $\mu$, hence to have first row at most $a$; degree $a$
forces first row at least $a$. At first row exactly $a$, the
row multiplicity in degree $a$ is $\mathbb S_\nu V$ for label
$(a,\nu)$, by the fully occupied, linear part of Schur--Weyl and
Pieri branching (equivalently, Theorem~\ref{thm:stable-augmented-multiplicities}
in its linear degree). The row highest weight $(a-1,1)$ therefore selects
exactly $\mu$. All invariant differential operators annihilate this
first occurrence, so $F$ survives the quadratic quotient.

Write $p_{ijk}=\sum_s x_s^iy_s^jz_s^k$, and put
\[
 A=\sum_s y_s\partial_{x_s},\qquad
 B=\sum_s z_s\partial_{x_s},\qquad
 C=\sum_s z_s\partial_{y_s}.
\]
The row weight $(a,2,1)$ and the two simple raising equations
determine the following combination of cubic moments and lowered
copies of $F$:

\begin{align*}
 G={}&p_{300}(AB-A^2C)F-(m-2)p_{210}BF
 +(2m-1)p_{210}ACF\\
 &-(m+1)p_{201}AF-(m^2-1)p_{120}CF
 +(m^2-1)p_{111}F.
\end{align*}
It has weight $(a,2,1)$ and is highest. For completeness, if
$X=\sum_sx_s\partial_{y_s}$ and
$Z=\sum_sy_s\partial_{z_s}$, then
\[
 XAF=(m-1)F,\quad XBF=-CF,\quad XACF=mCF,
\]
\[
 X(AB-A^2C)F=(m-2)BF-(2m-1)ACF,
\]
while $ZBF=ZACF=AF$, $ZCF=F$, $ZAF=0$, and
$Z(AB-A^2C)F=0$. The corresponding identities on the cubic moments
show $XG=ZG=0$. These are the simple raising operators of
$\operatorname{GL}_3$. Every operator used preserves the label
$\mu$, and Pieri's rule gives the indicated row summand once.

We show that $G$ survives the quadratic quotient by an exact
specialization over $K(\zeta)[\epsilon]/(\epsilon^2)$, where $\zeta$
is a primitive $2a$th root of unity. The frequencies $1,a-1,a$
sum to $2a$ for $xyz$. Every linear or quadratic frequency is
nonzero modulo $2a$, except that of $z^2$, whose coefficient
vanishes because $\epsilon^2=0$.
Set, for $0\le i<2a$,
\[
 x_i=\zeta^i,\qquad y_i=\zeta^{(a-1)i},\qquad
 z_i=\epsilon\zeta^{ai}.
\]
Every linear or quadratic moment is zero. For quadratic moments
containing two $z$'s this follows from $\epsilon^2=0$; every other
vanishing is a geometric sum with exponent not divisible by $2a$.
The same observation gives
\[
 p_{300}=p_{210}=p_{201}=p_{120}=0,
                      \qquad p_{111}=2a\epsilon.
\]
All factors of $F$ are nonzero. Its first determinant is
$\zeta^2-\zeta^{-2}$, which is nonzero for $a\ge3$, and each
remaining factor is a difference of distinct roots of unity.
Consequently
\[
                           G=2a(m^2-1)\epsilon F\ne0.
\]
The specialization kills the quadratic ideal, so $G$ is nonzero
in $Q_{n,3}$. Characteristic-zero scalar descent and functoriality
prove the assertion for all $V$.
\end{proof}

\section{Eulerian ideals and bounded exponents}
\label{app:eulerian-foundations}

We prove the coinvariant algebra statements used in
Section~\ref{sec:hook-action}.  The classical descent filtration and
the Gram-ideal lemma give a short proof of the complementary
annihilator identity inside the ordinary coinvariant algebra.

Put $P_n=K[x_1,\ldots,x_n]$, $p_s=\sum_i x_i^s$,
$I_n=(p_1,\ldots,p_n)$ and $C_n=P_n/I_n$.
For a label set $I$, ordered increasingly, write
\[
 x_I=\prod_{i\in I}x_i,\qquad
 \Delta_I=\det(x_i^{a-1})_{1\le a\le |I|,\ i\in I}.
\]
Thus the ideal $J_{n,d}$ of Section~\ref{sec:hook-action} is
generated by $g_I=x_I\Delta_I$ for $|I|=d$.
Also put $J_{n,0}=C_n$ and $J_{n,n}=0$.
For $d\ge1$, let
\[
 B_{n,d}=\operatorname{span}\{x^\alpha:0\le\alpha_i<d\},
 \qquad F_{n,d}=\operatorname{image}(B_{n,d}\longrightarrow C_n).
\]

\begin{proposition}[Eulerian complements and annihilators]
\label{prop:standalone-eulerian-foundations}
For $1\le d\le n$, the quotient projection induces a graded
$S_n$-equivariant isomorphism
\begin{equation}\label{eq:standalone-eulerian-complement}
 F_{n,d}\xrightarrow{\ \sim\ } C_n/J_{n,d}.
\end{equation}
Moreover, for $0\le d\le n$,
\begin{align}
 \Frob(J_{n,d};q)
 &=\sum_{\lambda\vdash n}s_\lambda
       \sum_{\substack{T\in\SYT(\lambda)\\\des(T)\ge d}}
                        q^{\maj(T)},
       \label{eq:standalone-eulerian-character}\\
 \operatorname{Ann}_{C_n}(J_{n,d})&=J_{n,n-d}.
       \label{eq:standalone-eulerian-annihilator}
\end{align}
On the coinvariant differential forms
\[
 \Omega_n=
 \frac{P_n\otimes\Lambda(\theta_1,\ldots,\theta_n)}
 {(p_1,\ldots,p_n,\dd p_1,\ldots,\dd p_n)},
 \qquad \delta_j=\sum_i x_i^j\partial_{\theta_i}\quad(j\ge1),
\]
the map $D_d=\delta_d\cdots\delta_1:\Omega_n^d\to C_n$
has image $J_{n,d}$ for $0\le d<n$; $D_0$ means the identity.
Its apolar adjoint is $\Theta_1\cdots\Theta_d$.
\end{proposition}

\begin{proof}
\emph{Determinants and straightening.}
For $1\le d\le n$, the two determinant presentations agree:
\begin{equation}\label{eq:standalone-determinant-presentations}
 I_n+(\Delta_I:|I|=d+1)
       =I_n+(x_I\Delta_I:|I|=d).
\end{equation}
Expansion along the constant row gives the inclusion from left
to right.  Conversely, append $(1,x_j,\ldots,x_j^d)^t$ to the
$d$ columns indexed by a fixed $I$ and sum the resulting
determinants over $j\in[n]$.  Each summand is zero or a signed
$(d+1)$-label Vandermonde.  The last column of the sum is
$(n,p_1,\ldots,p_d)^t$, so modulo $I_n$ the sum is
$\pm n x_I\Delta_I$.  This proves the reverse inclusion, and
also proves $g_{[n]}=0$ in $C_n$ when $d=n$.

We claim that every polynomial is congruent to an element of
$B_{n,d}$ modulo
$L_{n,d}=I_n+(\Delta_I:|I|=d+1)$.
Induct on $d$ for all $n$, and then descend on monomial support
for fixed $(n,d)$.  If $d=1$, the presentation above puts every
$x_i$ in $L_{n,1}$.  If $n\le d$, the Artin basis, with
exponents $0\le\alpha_i<i$, suffices.
Otherwise write a monomial of nonempty support $I$ as $x_I f$.
Apply the induction hypothesis for $(|I|,d-1)$ in $K[x_i:i\in I]$:
\[
 f=b+\sum_{s=1}^{|I|}p_s(I)h_s
       +\sum_{\substack{T\subseteq I\\|T|=d}}\Delta_T u_T,
 \qquad p_s(I)=\sum_{i\in I}x_i^s,
\]
where every exponent in $b$ is at most $d-2$ and all coefficients
use only variables in $I$.  After multiplication by $x_I$, the
first term is bounded as required and every determinant term
belongs to $L_{n,d}$ by
\eqref{eq:standalone-determinant-presentations}.  Finally,
\[
 x_Ip_s(I)h_s\equiv
       -\sum_{j\notin I}x_Ix_j^s h_s\pmod{I_n}.
\]
Each monomial on the right has larger support; for $I=[n]$ the
sum is zero.  The empty-support monomial is $1$.
This completes both inductions.  Taking homogeneous components
throughout proves that $F_{n,d}\to C_n/J_{n,d}$ is a graded
$S_n$-equivariant surjection.

\emph{The descent filtration.}
The classical descent straightening theorem identifies
\[
 F_{n,d}=\operatorname{span}
   \left\{b_w=\prod_{i\in\Des(w)}
                   (x_{w(1)}\cdots x_{w(i)}):\des(w)<d\right\}.
\]
These descent monomials belong to the Garsia--Stanton basis of
$C_n$, and the descent-representation refinement gives
\begin{equation}\label{eq:standalone-descent-filtration}
 \Frob(F_{n,d};q)=\sum_{\lambda\vdash n}s_\lambda
       \sum_{\substack{T\in\SYT(\lambda)\\\des(T)<d}}
                       q^{\maj(T)}.
\end{equation}
See \cite[Theorems~3.12 and~4.1]
{AdinBrentiRoichman2005} for the descent representations and
\cite[Section~4, Theorem~4.3 and the paragraph following it]
{GnedinGorinKerov2013} for the maximum-exponent realization.
These statements are defined
over $\mathbb Q$ and hence hold over $K$.
Put
\[
 r_d=\dim J_{n,d},\qquad
 R_d=\#\{w\in S_n:\des(w)\ge d\}.
\]
The surjection just proved gives $r_d\ge R_d$.

\emph{The Gram ideal and complementary dimensions.}
Newton's identities give $p_s=0$ in $C_n$ for every $s>0$.
The standard inner product is nondegenerate on
$W_n=\{(z_i)\in K^n:\sum_i z_i=0\}$, since $n\ne0$ in $K$.
Thus
\[
 v_s=(x_1^s,\ldots,x_n^s)\in C_n\otimes W_n,
 \qquad \langle v_s,v_t\rangle=p_{s+t}=0
 \quad(s,t\ge1).
\]
For $1\le d<n$, a generator $g_I$ of $J_{n,d}$ is a coordinate minor of
$v_1,\ldots,v_d$.  The product of such a minor and a generator
of $J_{n,n-d}$ is multilinear and alternating in two blocks of
vector arguments of total size $n>\dim W_n$.  Apply
Lemma~\ref{lem:gram} with the arguments independent, and then
substitute these prefix vectors.  Every Gram generator becomes
zero.  The polynomial ideal identity remains valid under this
substitution in $C_n$, even when vectors occur in both blocks.
Consequently
\begin{equation}\label{eq:standalone-prefix-orthogonality}
                 J_{n,d}J_{n,n-d}=0.
\end{equation}
The perfect Frobenius pairing of $C_n$ now gives
$r_d+r_{n-d}\le n!$.  On the other hand, complementing every
permutation value replaces its descent number by $n-1-\des(w)$,
so $R_d+R_{n-d}=n!$.  Hence
\[
 n!=R_d+R_{n-d}\le r_d+r_{n-d}\le n!,
\]
which forces $r_d=R_d$ and $r_{n-d}=R_{n-d}$.
The cases $d=0,n$ are immediate.

The graded surjection is therefore an isomorphism, proving
\eqref{eq:standalone-eulerian-complement}; subtracting
\eqref{eq:standalone-descent-filtration} from the ordinary
coinvariant character proves
\eqref{eq:standalone-eulerian-character}.
Equation~\eqref{eq:standalone-prefix-orthogonality} gives the
actual inclusion $J_{n,n-d}\subseteq\operatorname{Ann}(J_{n,d})$,
and the same dimension equality makes it an equality.

\emph{The prefix map.}
The contractions descend to $\Omega_n$, since
$\delta_j(\dd p_s)=s p_{s+j-1}=0$ in $C_n$.
Applying $D_d=\delta_d\cdots\delta_1$ to
$\theta_{i_1}\cdots\theta_{i_d}$ gives $g_I$.
The map is $C_n$-linear, so its image is exactly $J_{n,d}$.
Under the standard apolar pairing, $\delta_j^*=\Theta_j$;
reversing the composition order gives
$D_d^*=\Theta_1\cdots\Theta_d$.
\end{proof}

\subsection{A finite certificate for the Gram lemma}

Let $Z$ be a nondegenerate quadratic space of dimension $N$. On
$\Lambda^dZ\otimes\Lambda^eZ$, use the raising and contraction
operators $L,c$ with the ordinary tensor signs of Lemma~\ref{lem:gram}.
For $t=d+e-N>0$, the operator
\begin{equation}\label{eq:gram-finite-certificate}
 R_{d,e}=\sum_{j\ge0}
 \frac{(-1)^j c^{j+1}L^j}
      {(j+1)!\,t(t+1)\cdots(t+j)}
 \qquad\text{satisfies}\qquad LR_{d,e}=\operatorname{id}.
\end{equation}
Indeed,
$Lc^{j+1}L^j=c^{j+1}L^{j+1}+(j+1)(t+j)c^jL^j$,
so the sum telescopes; it is finite because exterior raising is
nilpotent. Interpreting its terms as alternating multilinear
polynomials gives an actual certificate in the ideal of
cross-block inner products.

\section{The canonical Schur resolution}
\label{app:schur-resolution}

We prove Theorem~\ref{thm:full-schur-resolution}. Retain
$B=K\oplus\mathfrak m$, $S=S_B$, and the shapes $\rho_i$, degrees $d_i$
and differences $e_i=\mu_i-\mu_{i+1}+1$ from
Section~\ref{sec:schur-resolution-statements}. The essential point is a
natural contraction for square-zero multiplication. Its count filtration
makes the perturbation for an arbitrary product finite.

\subsection{A dimension-independent contraction}

Put $E=K1\oplus\mathfrak m$ and give
$B_0=K\oplus\mathfrak m$ the square-zero product. On
\[
                        C_i^0=\Lambda^i\mathfrak m\otimes E^{\otimes n}
\]
let $d_0$ be its moment differential: remove a wedge input and place
it in each empty coefficient position. Define $c:C_i^0\to C_{i+1}^0$
by removing each occupied coefficient $x$ in turn, replacing it by
$1$, and sending the exterior factor $\omega$ to $x\wedge\omega$.
These formulas are natural for every linear map of $\mathfrak m$;
they require no bilinear form. Both operators preserve the count
\[
           k=\text{number of occupied positions}+\text{exterior degree}.
\]
They satisfy $d_0^2=c^2=0$. With any temporary real orthonormal basis,
$c$ is the adjoint of $d_0$ for the standard tensor and exterior
inner products.

\begin{lemma}[The content Laplacian]
\label{lem:schur-content-laplacian}
Let $\cont(\alpha)=\sum_{(r,s)\in\alpha}(s-r)$. On label
$\mu\vdash n$ and polynomial Schur type $\lambda\vdash k$ in $C_i^0$,
the operator $L=d_0c+cd_0$ is the scalar
\begin{equation}\label{eq:schur-content-laplacian}
           \cont(\mu)-\cont(\lambda)+n-\binom{n-k+1}{2}.
\end{equation}
Here the binomial denotes its polynomial value, also for negative
arguments. There is a natural contraction onto $\mathcal H=\ker L$,
with inclusion $j$, projection $p$, and homotopy $h$ satisfying
\[
 d_0h+hd_0=\id-jp,\qquad h^2=0,\qquad hj=ph=0.
\]
On label $\mu$, the harmonic space in homological degree $i$ is
$\mathcal H_i\cong\Schur_{\rho_i}\mathfrak m$, in count $d_i$.
\end{lemma}

\begin{proof}
First let $q=\dim\mathfrak m<\infty$. On the coefficient module
$\Schur_\mu(E)$, write $E_{ab}$ for the matrix-unit action, with
index $0$ denoting the unit. Exterior creation and contraction give
\[
 d_0=\sum_{a=1}^q E_{a0}\iota_a,
 \qquad c=\sum_{a=1}^qE_{0a}\varepsilon_a,
 \qquad
 L=\sum_aE_{a0}E_{0a}+iE_{00}-\sum_{a,b}E_{ab}\varepsilon_b\iota_a.
\]
Let $\Omega_q=\sum_{a,b}E_{ab}E_{ba}$ be the quadratic Casimir.
Its eigenvalue on $\Schur_\lambda K^q$ is
$2\cont(\lambda)+q|\lambda|$. The cross-Casimir identity rewrites
the preceding formula as
\[
 2L=\Omega_{q+1}^{\mu}-\Omega_q^{\mathrm{total}}
       +\Omega_q^{\Lambda^i}
       -E_{00}^2+(2i-q-1)E_{00}+n.
\]
Now $E_{00}=n-k+i$ and
$\Omega_q^{\Lambda^i}=i(q+1-i)$. Substitution cancels $q$ and $i$
and gives \eqref{eq:schur-content-laplacian}.

Over the real numbers the adjoint identity implies
$\ker L=\ker d_0\cap\ker c$. Let $L^+$ invert every nonzero
eigenvalue and be zero on the kernel. Then $h=cL^+$ and the spectral
projection $jp$ give the asserted contraction. Its formulas are
rational: the eigenvalues in \eqref{eq:schur-content-laplacian}
are integers and the Schur projections are natural. They therefore
hold over every characteristic-zero field and commute with arbitrary
linear maps, including maps that lower dimension.

To identify the harmonic spaces, apply the pure-resolution theorem
\cite[Theorem~2.2]{SamWeyman} to the square-zero presentation, with
$\alpha=\tau$ and $\beta=\tau^+$. Its successive generator shapes
are precisely $\rho_i$, and its degree shifts, including the initial
shift $|\tau|$, are $d_i$. Hence the square-zero Koszul homology is
$\Schur_{\rho_i}\mathfrak m$ in count $d_i$. This harmonic copy
occurs once: its coefficient branching shape deletes row $i+1$
of $\mu$, and its exterior Pieri strip adds one box to each of the
first $i$ rows. The case $i=0$ is minimum occupancy.

For infinite-dimensional $\mathfrak m$, each fixed-count construction
passes through finite-dimensional vector subspaces. The scalar formula
is independent of their dimension, so the contractions are compatible
under inclusion and define the same natural operators on their union.
\end{proof}

The diagram convention in \cite{SamWeyman} is transposed, while its
highest-weight tuples agree with those used here. This argument uses
its full pure-resolution theorem; the earlier presentation used only
the Pieri composition lemma.

\subsection{Multiplication as a finite perturbation}

\begin{proof}[Proof of Theorem~\ref{thm:full-schur-resolution}]
We first transfer the collision differential to compute moment homology.
We then adjoin formal trace variables and transfer once more to obtain
the free resolution; setting those variables to zero recovers the first
complex. Both transfers terminate because they lower coefficient count.

For the original algebra, write the moment differential as
$d_B=d_0+\delta$. The collision operator $\delta$ removes a wedge
input and multiplies it into each occupied position. It lowers count
by one, while $h$ preserves count. The perturbation lemma, with our
convention $d_0h+hd_0=\id-jp$, gives the transferred differential
\[
 p\delta(1+h\delta)^{-1}j
       =\sum_{r\ge0}(-1)^r p\delta(h\delta)^rj.
\]
Every series is finite on each vector. This is the standard formula
of \cite[\S2, especially~2.4]{Crainic}, with its homotopy replaced by
$-h$. A term from $\mathcal H_i$ to $\mathcal H_{i-1}$ must lower
count from $d_i$ to $d_{i-1}$, so only one term survives:
\begin{equation}\label{eq:schur-small-differential}
 \partial_i=(-1)^{e_i-1}
          p_{i-1}\delta(h\delta)^{e_i-1}j_i.
\end{equation}
Thus the small complex computes all moment homology, with no choices
of higher splittings.

To obtain an actual free resolution, use
$C_{\mathrm{free}}=S\otimes\Lambda^\bullet\mathfrak m
                    \otimes\Schur_\mu(E)$.
Let $d_T$ remove a wedge input $b$ and multiply the $S$ factor by
$T_b$. Subtracting $d_T$ enforces the action $T_b=p_n(b)$.
The differential $d_B-d_T$ resolves $\Schur_\mu(B)$ under
$f(T)\otimes v\mapsto f(p_n)v$. Indeed, these are polynomial-variable
elimination operators $p_n(b)-T_b$; eliminating any finite set is
exact by successive substitution, and their directed union handles
infinitely many variables.

Extend the preceding contraction $S$-linearly, starting with $d_0$,
and perturb it by $\mathcal A=\delta-d_T$. Both terms lower
\emph{coefficient} count by one; the increase of formal polynomial
degree under $d_T$ is irrelevant to this count. Therefore the
transferred free differential is the finite formula
\begin{equation}\label{eq:schur-free-differential}
 D_i=(-1)^{e_i-1}
       p_{i-1}\mathcal A(h\mathcal A)^{e_i-1}j_i.
\end{equation}
It gives a chain homotopy equivalent free complex with the stated
terms, and hence a resolution. Setting $T=0$ gives
\eqref{eq:schur-small-differential}.

For the resolution filtration, give a formal trace count one and
the generator of $F_i$ count $d_i$. A term with no collisions has
formal degree $e_i$ and preserves this total count; every collision
lowers it. At $\delta=0$ the complex is a minimal equivariant pure
resolution, so its leading maps are the classical nonzero Pieri maps.

The corrected harmonic inclusion is $(1+h\mathcal A)^{-1}j$.
On degree zero, $\mathcal A j_0=0$, so its augmentation is exactly
the minimum-occupancy map. Choose the Schur identification line
for $\mathcal H_1$ so that the leading part of $D_1$ equals the
previously fixed Pieri inclusion. Every lower-count term lies below
$d_1$, where the omitted Schur summands are absent. The unique lift
of that leading relation is therefore the graph lift~\eqref{eq:schur-relation-graph}, so
$D_1=R_{\mu,B}$ literally with this normalization. The earlier
quotient map used $\beta_{\mu,B}=-\operatorname{ev}_{T=0}R_{\mu,B}$;
accordingly $\partial_1=-\beta_{\mu,B}$. The kernels and cokernels
are unchanged. For example, before this line normalization the raw
adjoint convention at $n=1$, $\mathfrak m=K a$, $a^2=ca$ gives
$D_1=cT_a-T_a^2$.

All operators are built from natural Schur projections, unit and
augmentation, and multiplication. This proves the claimed functoriality
and compatibility with gradings. It gives algebra-automorphism
equivariance for general $B$; full $\operatorname{GL}(\mathfrak m)$
equivariance belongs to the square-zero specialization.
\end{proof}

\section{Finer row-sensitive character ranges}
\label{app:finer-character-bounds}

The uniform character theorem needs only the Schur relation, quadratic
reduction, bivariate and hook exactness, and one cubic witness. This
appendix retains the separate row-sensitive trace bounds and their
closed boundary, which can give better ranges for individual rows.

\subsection{Row-sensitive support of moment Koszul homology}

The first row of the row partition gives a second measure of how
early a tail relation can occur.  A relation on $t+1$ positive
inputs has many linear inputs until its degree becomes large.
Their symmetric product forces a long first row.  Comparing this
trace relation with the free moment complex gives the following
joint restriction on row and label.

\begin{theorem}[A row-and-label strip for positive moment homology]
\label{thm:row-sensitive-moment-homology-strip}
Suppose $\mathbb S_\lambda V\otimes S^\mu$ occurs in
$H_j(\mathcal K_{n,k,D})_D$, where $j\ge1$ and
$\lambda\vdash D$.  Write $a=\lambda_1$.  Then
\begin{equation}\label{eq:row-sensitive-moment-homology-strip}
 n-D+2j\le\mu_1\le
       \left\lfloor\frac{D+a}{2}\right\rfloor-j.
\end{equation}
When $3\le D\le n$, the sharper strip is
\begin{equation}\label{eq:row-sensitive-quadratic-moment-strip}
 n-D+3j\le\mu_1\le
 \min\left\{\left\lfloor\frac{D+a}{2}\right\rfloor-j,\ D-2\right\}.
\end{equation}
More generally, if all moments of orders at most $m\ge1$ form a
regular sequence in $P$, then
\begin{equation}\label{eq:row-sensitive-regular-moment-strip}
 n-D+(m+1)j\le\mu_1\le
 \left\lfloor\frac{D+a-(m+1)(j-1)-2}{2}\right\rfloor.
\end{equation}
\end{theorem}

\begin{proof}
For the general lower bound, remove the linear moments.  Every
remaining exterior generator has weight at least two, and a
polynomial of degree at most $D-2j$ has first label row at least
$n-D+2j$.  In the range $D\le n$, quadratic exactness raises the
minimum remaining weight to three, as in
Lemma~\ref{lem:parabolic-koszul-detection}.

We prove the new upper bound.  For $1\le t\le n$, let $S_t$
act on the last $t$ labels, put $r=n-t$, and write
\[
 A_r=\operatorname{Sym}(V\otimes K^r),\qquad
 A=\operatorname{Sym}V,\qquad
 F=\operatorname{Sym}(A_+),\qquad
 T=(A^{\otimes t})^{S_t}.
\]
The map $F\to T$ sends the formal trace $[b]$ to
$\sum_{i=1}^t b(X_i)$.  By the trace presentation used in
Theorem~\ref{thm:universal-selected-trace-presentation}, its kernel
$I_t$ is generated by the equivariant image of
$\operatorname{Sym}^{t+1}(A_+)$.  Explicitly, the relation on
$b_1,\ldots,b_{t+1}$ is
\[
 \sum_{\pi\in\Pi_{t+1}}
 (-1)^{t+1-|\pi|}
 \prod_{B\in\pi}(|B|-1)!
 \prod_{B\in\pi}\left[\prod_{i\in B}b_i\right].
\]
The internal degree is $\sum_i\deg b_i$; the same formula is
extended multilinearly to arbitrary inputs in $A_+$.

A degree-$g$ generator has at least $2t+2-g$ linear inputs.
Splitting $\operatorname{Sym}^{t+1}(A_+)$ by the number $\ell$
of linear inputs therefore expresses its relevant summands as
$\operatorname{Sym}^{\ell}V\otimes W$, with
$\ell\ge2t+2-g$ and $W$ polynomial.  The
Littlewood--Richardson rule shows that every constituent has first
row at least $2t+2-g$.  Multiplication by any polynomial
representation, or tensoring with one, cannot decrease that first
row.  

First eliminate the global linear moments in this comparison.
They form a regular sequence both in $A_r\otimes F$ and in
$A_r\otimes T=P^{S_t}$.  In the first ring their formal trace
variables are independent.  In the second, splitting the label
representation into the invariant diagonal and the centered
subspace gives a polynomial factor on the global linear moments;
taking $S_t$-invariants preserves this factorization.
Consequently the exact sequence with kernel $A_r\otimes I_t$
remains exact after quotienting by those moments.  Its new kernel
$\overline I_t$ is the image of that ideal, generated by the same
equivariant trace-relation sources.  Their first-row bounds are
therefore unchanged.

The source after elimination is the polynomial algebra
\[
       A_r\otimes\operatorname{Sym}(A_{\ge2}).
\]
In total degrees at most $D$, form its Koszul complex
$\overline{\mathcal L}$ on the remaining global moments
$[X^\alpha]+\sum_{i\le r}X_i^\alpha$, $2\le|\alpha|\le D$.
A polynomial change of the independent formal trace variables
makes this complex exact in positive degrees.  It surjects onto
the corresponding reduced moment complex of $P^{S_t}$, whose
homology is that of $\mathcal K_{n,k,D}^{S_t}$ because the removed
linear sequence is regular.  Its kernel complex has terms
\[
 \mathcal N_i=\overline I_t\otimes
           \bigwedge^i\!\left(\bigoplus_{d=2}^D
                                      \operatorname{Sym}^dV\right).
\]
All maps are equivariant.  Using the full free trace algebra and
then taking the stated degrees gives the same terms: variables
and relations of weight greater than $D$ cannot contribute.

Positive homology of the finite-label complex is detected one
homological degree earlier in this trace-relation kernel.
If a row $\lambda\vdash D$ occurs in $\mathcal N_{j-1}$,
it comes from a trace-relation generator of some degree $g$,
polynomial factors, and $j-1$ exterior generators of weight at
least two.  Hence $g\le D-2j+2$.  By the first-row bound and
Littlewood--Richardson containment,
\[
        a\ge2t+2-g\ge2t+2j-D.
\]
The long exact sequence of
$0\to\mathcal N\to\overline{\mathcal L}\to
\overline{\mathcal K_{n,k,D}^{S_t}}\to0$ identifies positive homology in
degree $j\ge2$ with $H_{j-1}(\mathcal N)$, and injects degree-one
homology into $H_0(\mathcal N)$.  Thus its row-$\lambda$ part
vanishes whenever $D+a<2t+2j$.

Now take $t=\mu_1$.  Branching gives a nonzero $S_t$-fixed
vector in $S^\mu$, and exactness of invariants identifies
$H_j(\mathcal K_{n,k,D})^{S_t}$ with the homology just considered.
Therefore $D+a\ge2\mu_1+2j$, which is the desired upper bound.
This also permits $r=0$: the selected-label polynomial factor is
then $K$.

Combining this estimate with
Lemma~\ref{lem:parabolic-koszul-detection} gives
\eqref{eq:row-sensitive-quadratic-moment-strip}.

If the moments through order $m$ are regular in $P$, their Koszul
complex stays exact after taking $S_t$-invariants, so they are
regular in $P^{S_t}$ as well.  The same comparison can therefore
eliminate all their independent formal traces first.  The
trace-relation sources retain their row bound, while every
remaining exterior generator has weight at least $m+1$.
Thus $g\le D-(m+1)(j-1)$, giving the upper bound in
\eqref{eq:row-sensitive-regular-moment-strip}; removal of the same
regular sequence gives its lower bound by monomial support.
\end{proof}

The \emph{row-$\lambda$ Euler character} is the
$S_n$-representation
\[
 \sum_{j\ge0}(-1)^j
 \operatorname{Hom}_{\operatorname{GL}(V)}
       \bigl(\mathbb S_\lambda V,
              (P\otimes\bigwedge^jU_D)_D\bigr).
\]
Its evaluation at $\sigma\in S_n$ is obtained by projecting the
total-degree-$D$ Koszul product in
\eqref{eq:explicit-stable-character-polynomial} onto the Schur
character $s_\lambda$.  Thus exactness can be tested for one row
partition while other rows in the same multidegree retain homology.

\subsection{Detection of the first row-sensitive trace relations}

At equality in the trace-relation bound, all positive inputs are
linear or quadratic.  This leaves a multiplicity-one source whose
highest vector can be written with bordered minors.  A factorial
evaluation of that vector proves that the first possible relation
survives negative-moment substitution.  It therefore contributes
no positive moment homology.

\begin{lemma}[The boundary trace is detected by selected labels]
\label{lem:boundary-trace-detection}
Let $\lambda\vdash D$, write $a=\lambda_1$, and suppose
\begin{equation}\label{eq:boundary-trace-equality}
                         D+a=2t+2.
\end{equation}
Set $q=t+1$, $r=n-t$, and use the free tail-trace algebra
$F=\operatorname{Sym}(A_+)$, its kernel $I_t$, and the selected-label
polynomial algebra $A_r$ from the proof of
Theorem~\ref{thm:row-sensitive-moment-homology-strip}.
The degree-$D$, row-$\lambda$ part of $A_r\otimes I_t$ has
multiplicity at most one.  It is zero unless
\begin{equation}\label{eq:boundary-even-tail-shape}
       \lambda=(a,2\kappa),\qquad
       |\kappa|=\frac{D-a}{2},\qquad 2\kappa_1\le a.
\end{equation}
If \eqref{eq:boundary-even-tail-shape} holds and
$r\ge\ell(\lambda)$, the negative-moment specialization
\[
 \operatorname{ev}_{-}:A_r\otimes F\longrightarrow A_r,
 \qquad [b]\longmapsto-\sum_{i=1}^r b(X_i),
\]
is nonzero, hence injective, on that single irreducible summand.
Consequently the row-$\lambda$ part of
$H_j(\mathcal K_{n,k,D})^{S_t}_D$ vanishes for every $j\ge1$
if the even-tail condition fails, or if $r\ge\ell(\lambda)$.
A row of length greater than $\dim V$ is absent.
\end{lemma}

\begin{proof}
The trace-relation generators have source
$\operatorname{Sym}^{q}(A_+)$.  A degree-$g$ generator has first
row at least $2q-g$.  Under \eqref{eq:boundary-trace-equality}, a
row-$\lambda$ constituent of total degree $D$ therefore requires
$g=D$, with no positive-degree multiplier.  Moreover its source
must have exactly $a$ linear inputs and
$m=(D-a)/2$ quadratic inputs.  A larger number of linear inputs
would force a larger first row, and any input of degree at least
three would force more linear inputs.  The entire possible source
is consequently
\begin{equation}\label{eq:boundary-trace-source}
     \operatorname{Sym}^{a}V\otimes
          \operatorname{Sym}^{m}(\operatorname{Sym}^{2}V).
\end{equation}

The symmetric-matrix Cauchy identity gives
\[
 \operatorname{Sym}^{m}(\operatorname{Sym}^{2}V)
       =\bigoplus_{\kappa\vdash m}\mathbb S_{2\kappa}V.
\]
Pieri's rule now says that a constituent of
\eqref{eq:boundary-trace-source} with first row exactly $a$ must
be $(a,2\kappa)$, once: its horizontal strip gives
$2\kappa_i\ge\lambda_{i+1}$, while both sides sum to $D-a$.
This proves the asserted multiplicity bound and shape restriction.

It remains to show that this possible summand survives the map
from trace-relation sources to selected polynomials.  Multiplicity
one reduces injectivity to a single nonzero image.  We will choose
a highest vector whose negative-moment image can be evaluated
without cancellation.
Realize
\eqref{eq:boundary-trace-source} using vector generators
$u_1,u_2,\ldots$ of weight one and symmetric-matrix generators
$M_{ij}=M_{ji}$ of weight two.  For $j\ge1$ let
\begin{equation}\label{eq:boundary-bordered-minor}
 p_j(u,M)=
 \det\begin{pmatrix}
  0&u_1&\cdots&u_{j+1}\\
  u_1&M_{11}&\cdots&M_{1,j+1}\\
  \vdots&\vdots&&\vdots\\
  u_{j+1}&M_{j+1,1}&\cdots&M_{j+1,j+1}
 \end{pmatrix}.
\end{equation}
It has $u$-degree two, $M$-degree $j$, and highest weight
$(2^{j+1})$: each raising operator acts by an elementary row
operation together with the corresponding column operation on a
leading principal block.  Set $\kappa_{\ell(\kappa)+1}=0$.
Then
\begin{equation}\label{eq:boundary-highest-vector}
 P_\kappa(u,M)=u_1^{a-2\kappa_1}
       \prod_{j\ge1}p_j(u,M)^{\kappa_j-\kappa_{j+1}}
\end{equation}
is a highest-weight vector of weight $(a,2\kappa)$, with
$u$-degree $a$ and $M$-degree $m$.  For $\kappa=\varnothing$
the formula means $P_\varnothing=u_1^a$.

Introduce auxiliary variables $\xi_1,\ldots,\xi_r$ and the
linear functional
\[
       \mathcal E(\xi_1^{\alpha_1}\cdots\xi_r^{\alpha_r})
                       =\prod_{i=1}^r\alpha_i!.
\]
For the multilinear fundamental trace relation $\Psi_q$, its
negative-moment specialization is exactly
\begin{equation}\label{eq:trace-factorial-moment-identity}
 \operatorname{ev}_{-}\Psi_q(b_1,\ldots,b_q)
  =(-1)^q\mathcal E\!\left[
              \prod_{s=1}^{q}\left(\sum_{i=1}^r
                                  \xi_i b_s(X_i)\right)\right].
\end{equation}
Indeed, all set-partition terms in $\Psi_q$ acquire the same
sign $(-1)^q$.  For a fixed assignment of the inputs to the
selected labels, the sum over set partitions refining each fiber
is the product of the fiber factorials: its blocks are the cycles
of a permutation.  These are precisely the coefficients given
by $\mathcal E$.  Thus the normalization in
\eqref{eq:trace-factorial-moment-identity} is fixed, without any
choice of polarization convention.

Write $\widehat\Psi_q:\operatorname{Sym}^q(A_+)\to I_t$ for
the linear map induced by the multilinear trace relation.  The
actual relation used here is $\widehat\Psi_q(P_\kappa)$.
Applying \eqref{eq:trace-factorial-moment-identity} to it amounts
to substituting
\[
             u=\sum_i\xi_iX_i,\qquad
             M=\sum_i\xi_iX_iX_i^{\mathsf T},
\]
then applying $(-1)^q\mathcal E$.  Since
$r\ge\ell(\lambda)$, specialize $X_i$ to the $i$th coordinate
vector for $1\le i\le\ell(\lambda)$ and set all remaining
$X_i$ equal to zero.  Now $u_i=\xi_i$ and $M$ is diagonal, so
\[
 p_j(u,M)=-\left(\prod_{i=1}^{j+1}\xi_i\right)
                    \left(\sum_{i=1}^{j+1}\xi_i\right).
\]
The resulting $P_\kappa$ is a nonzero polynomial all of whose
nonzero coefficients have sign $(-1)^{\kappa_1}$.
Its image under $\mathcal E$ is therefore a nonzero integer.
This proves nonvanishing over every field of characteristic zero.
The case $\kappa=\varnothing$ gives the same conclusion directly
from $\mathcal E(\xi_1^a)=a!$.

Finally compare the free Koszul complex on
$[X^\alpha]+\sum_{i\le r}X_i^\alpha$ in $A_r\otimes F$
with the moment complex of $P^{S_t}$.  The free complex is exact
in positive degrees; its zeroth homology in the degrees under
consideration is $A_r$.  Its kernel complex has terms
\[
       (A_r\otimes I_t)\otimes\bigwedge^j U_D.
\]
Under \eqref{eq:boundary-trace-equality}, none of its terms with
$j\ge1$ can contain row $\lambda$ in total degree $D$:
an exterior factor would require a trace generator of degree
strictly less than $D$, whose first row is greater than $a$.
All homology in homological degrees at least two is therefore zero
in this row.  The remaining part of the long exact sequence is
\[
 0\longrightarrow H_1(\mathcal K_{n,k,D})^{S_t}_{D,\lambda}
      \longrightarrow
 (A_r\otimes I_t)_{D,\lambda}
             \xrightarrow{\operatorname{ev}_{-}}
                         (A_r)_{D,\lambda}.
\]
The source is zero in the excluded shapes, and in the remaining
case it is a single irreducible on which the displayed map is
nonzero.  This proves the homology assertion.
\end{proof}

\begin{corollary}[The closed upper boundary and row characters]
\label{cor:boundary-row-character-range}
Let $3\le D\le n$, let $\lambda\vdash D$, and put
$a=\lambda_1$.  A joint row-and-label constituent of positive
moment homology satisfies
\begin{equation}\label{eq:boundary-refined-moment-strip}
 n-D+3j\le\mu_1\le
 \min\left\{
       \left\lfloor\frac{D+a}{2}\right\rfloor-j,
       \left\lfloor\frac{D+a-3}{2}\right\rfloor
       \right\}.
\end{equation}
In particular the row-$\lambda$ Koszul Euler character is the
actual row-$\lambda$ character of $(R_{n,k})_D$ whenever
\begin{equation}\label{eq:boundary-row-character-range}
 n\ge\max\left\{D,\;
          D-2+\left\lfloor\frac{D+a-3}{2}\right\rfloor\right\}.
\end{equation}
For a single joint row-and-label component, equality of the Euler
and actual characters already holds, for $3\le D\le n$, if
\begin{equation}\label{eq:boundary-row-label-character-range}
                  \mu_1\ge
                      \left\lfloor\frac{D+a-1}{2}\right\rfloor.
\end{equation}
\end{corollary}

\begin{proof}
The $j$-dependent bounds are
Theorem~\ref{thm:row-sensitive-moment-homology-strip}.
If $D+a$ is odd, they already give the additional upper bound
in \eqref{eq:boundary-refined-moment-strip}.  If $D+a$ is even,
the only additional case to exclude is $j=1$ and
$t=\mu_1=(D+a-2)/2$.  Branching gives an $S_t$-fixed vector.
At that value, Lemma~\ref{lem:boundary-trace-detection} excludes
the row unless $\lambda=(a,2\kappa)$.  In the latter case,
with $m=|\kappa|=(D-a)/2$, the hypothesis $n\ge D$ gives
\[
       r=n-t\ge D-t=m+1\ge\ell(\kappa)+1=\ell(\lambda).
\]
The same lemma therefore excludes the remaining case as well.
This proves the strip.

Write $C=\lfloor(D+a-3)/2\rfloor$.  Under
\eqref{eq:boundary-row-character-range}, positive homology would
require $\mu_1\ge n-D+3\ge C+1$ and $\mu_1\le C$,
a contradiction.  Thus zeroth homology equals the Koszul Euler
character in that row.  The final assertion follows directly
from the upper bound and $C+1=\lfloor(D+a-1)/2\rfloor$.
\end{proof}

For the widest row $a=D$, the strip gives
$n\ge\max\{D,2D-4\}$; for bounded $a$ its threshold is
$\frac32D+O_a(1)$, and balanced three-row partitions have threshold
$\frac53D+O(1)$. These finer rowwise conclusions are independent of
the shorter proof of the uniform theorem.

\section{The odd moment differential algebra}
\label{app:odd-derived}

The two-column reduction has a simple differential-algebra model.
It makes the ring duality explicit and identifies the minimal
syzygy-generating layer of Proposition~\ref{prop:two-column-module}
as values of triple Massey products. Throughout, $m=n-1$ and
$\mathcal E$, $\omega$, $L$, $c$, and $P_e$ have the meanings of
Section~\ref{sec:two-column-homology}.

\subsection{Totalization and the Heisenberg algebra}

On a purely odd two-dimensional row space, the only moment generators
are
$u=\sum_i\theta_i$, $v=\sum_i\phi_i$, and
$q=\sum_i\theta_i\phi_i$.
Let $\xi,\zeta,\eta$ be their Koszul generators. For internal degree
$D$ and homological degree $j$, assign cohomological degree $D-j$
and twist the product by
\begin{equation}\label{eq:odd-totalization}
                       a\star b=(-1)^{D(a)j(b)}ab.
\end{equation}
The original commutation sign
$(-1)^{D(a)D(b)+j(a)j(b)}$ becomes the cohomological commutation
sign, and the differential satisfies the corresponding Leibniz rule
without modification.

Now $\xi,\zeta$ have cohomological degree zero, with
$d\xi=u$ and $d\zeta=v$. Centering the odd variables and making the
$GL(U)$-equivariant change
\[
 \omega=q-\frac{uv}{n},\qquad
 \eta'=\eta-\frac{\xi\star v-\zeta\star u}{2n}
\]
gives $d\eta'=\omega$. The two factors
$(K[\xi]\otimes\Lambda(u),d\xi=u)$ and
$(K[\zeta]\otimes\Lambda(v),d\zeta=v)$ have only constant cohomology
in characteristic zero. Their joint augmentation is a
quasi-isomorphism onto
\begin{equation}\label{eq:odd-heisenberg-model}
 A_m=\bigl(\Lambda(\theta_1,\ldots,\theta_m,
                  \phi_1,\ldots,\phi_m,\eta),\,d\eta=\omega\bigr).
\end{equation}
This is the Chevalley--Eilenberg algebra of the $(2m+1)$-dimensional
Heisenberg Lie algebra, up to the sign convention for its central
generator. Its cohomology is classical; see \cite{CairnsJambor}.
Here it is the actual complete odd specialization of the moment
complex, with its product.

Put $R=\mathcal E/(\omega)$ and
$I=\eta\operatorname{Ann}_{\mathcal E}(\omega)$. Cycles and
boundaries in \eqref{eq:odd-heisenberg-model} give the ring identity
\begin{equation}\label{eq:odd-cohomology-double}
 H^*(A_m)=R\oplus I,\qquad I^2=0.
\end{equation}
The exterior top-product pairing identifies $I^j$ with the dual of
$R^{2m+1-j}$, tensored with the orientation line
$(\det U)^{m+1}\otimes(\det W)^2$.
The orientation class is $\eta\omega^m/m!$, of cohomological degree
$2n-1$ and internal degree $2n$. Thus \eqref{eq:odd-cohomology-double}
is a Poincar\'e-duality trivial extension. By
Proposition~\ref{prop:two-column-module} and this pairing,
$\operatorname{Soc}(R)=R_m=P_m$; its type is the Catalan number
$\frac1{n+1}\binom{2n}{n}$, so $R$ is not Gorenstein when $n\ge2$.
Also,
\[
 \dim R_e=\binom{2m}{e}-\binom{2m}{e-2}\quad(0\le e\le m),
 \qquad
 \dim H^*(A_m)=\binom{2n}{n}.
\]

Fixing $\omega$ gives an action of $Sp_{2m}$ on the differential
algebra. Here $P_0=K$ and $P_e$ is the primitive fundamental
representation for $1\le e\le m$. This larger symmetry mixes row
components; a fixed $GL(U)$ highest-weight space need not be stable.

\subsection{Massey products and all positive homology}

For $\beta\in U\otimes W$, set
$\iota_\beta=[c,\varepsilon_\beta]$ and
$T_\beta=\varepsilon_\beta\iota_\beta$.
In Darboux coordinates,
$\iota_{\theta_i}=\partial_{\phi_i}$ and
$\iota_{\phi_i}=-\partial_{\theta_i}$.
The operator $T_\beta$ is the rank-one symplectic infinitesimal
transformation, acting as a derivation on $\mathcal E$.

\begin{theorem}[Massey generation of positive moment homology]
\label{thm:odd-massey-generation}
For $m\ge1$, $\alpha\in P_m$, and $\beta\in U\otimes W$, the triple
Massey product is defined, has zero indeterminacy, and equals
\begin{equation}\label{eq:odd-massey-action}
 \langle[\alpha],[\beta],[\beta]\rangle
       =\{(-1)^m[\eta T_\beta\alpha]\}.
\end{equation}
Its values span the minimal generator space $\eta P_m$ and generate
all positive moment homology as an $R$-module.
\end{theorem}
\begin{proof}
Since $L\alpha=0$, the identity
$[L,\iota_\beta]=\varepsilon_\beta$ gives
$L\iota_\beta\alpha=\beta\alpha$.
Use the convention that $dx=ab$, $dy=bc$ represents the triple
product by $ay+(-1)^{|a|+1}xc$.
Taking $x=(-1)^m\eta\iota_\beta\alpha$ and $y=0$ proves
\eqref{eq:odd-massey-action}. Its indeterminacy is
$[\alpha]H^1(A_m)+H^m(A_m)[\beta]=0$:
these groups have no $\eta$ part, and their indicated products lie
in $R^{m+1}=0$.

Polarizing $T_\beta$ identifies their span with
$\mathfrak{sp}_{2m}\cong\Sym^2(U\otimes W)$.
The middle primitive module $P_m$ is nontrivial and irreducible, so
$\mathfrak{sp}_{2m}P_m=P_m$.
Proposition~\ref{prop:two-column-module} now gives the final assertion.
\end{proof}

For example,
\[
 \left\langle[\theta_1\cdots\theta_m],[\phi_1],[\phi_1]\right\rangle
       =\{[\eta\theta_2\cdots\theta_m\phi_1]\}
\]
is nonzero. Consequently the odd moment differential algebra is
nonformal for every $n\ge2$.
For $n=1$ the reduced algebra is $\Lambda(\eta)$ and is formal.
These assertions concern this odd specialization, not every
bosonic specialization.

\subsection{The complete ternary operation}

Define the Lefschetz partial inverse $J$ by
\[
 J(L^jf)=L^{j-1}f\quad(j>0),\qquad J(f)=0\quad(f\in P_e),
 \qquad \pi=1-LJ,\quad Q=1-JL.
\]
The operators $\pi$ and $Q$ project to initial and terminal string
vectors, respectively. Identify $R$ with $\pi(\mathcal E)$.

\begin{proposition}[A minimal model with only binary and ternary operations]
\label{prop:odd-ternary-model}
There is a minimal $A_\infty$ model on $R\oplus I$, with its usual
cohomology product, for which
\begin{equation}\label{eq:odd-ternary-formula}
 m_3(a,b,c)=\eta Q\bigl(J(ab)c-aJ(bc)\bigr)
                                   \qquad(a,b,c\in R).
\end{equation}
It has $m_3=0$ if an input lies in $I$, and $m_r=0$ for every $r\ge4$.
\end{proposition}
\begin{proof}
There is an internal-degree-preserving deformation retract with
$h(a)=\eta J(a)$, $h(\eta b)=0$ and
$p(a+\eta b)=\pi(a)+\eta Q(b)$.
If $i$ includes the chosen representatives, then $dh+hd=1-ip$.
Apply homotopy transfer, with the cohomological conventions of
\cite[Section~3]{KellerAInfinity}:
\[
 \sum_{r+s+t=N}(-1)^{r+st}
 m_{r+1+t}
   (\operatorname{id}^{\otimes r}\otimes m_s
                  \otimes\operatorname{id}^{\otimes t})=0.
\]
The tensor products of maps carry their Koszul signs.
The two ternary transfer terms are
$p(h(ab)c-(-1)^{|a|}a h(bc))$; moving $a$ across $\eta$ gives
\eqref{eq:odd-ternary-formula}.
Every transferred operation preserves internal degree. Since
$|m_r|=2-r$, its output moment-homological index is
$\sum j_{\mathrm{in}}+r-2$. Only indices zero and one survive,
which proves the vanishing assertions.
\end{proof}

In this convention
$m_3(\alpha,\beta,\beta)=-\eta T_\beta\alpha$ for $\alpha\in P_m$.
The Massey convention in \eqref{eq:odd-massey-action} differs by
$(-1)^{m+1}$, as it must. Formula~\eqref{eq:odd-ternary-formula}
also specifies the operation when an ordinary Massey product is
not defined. The cutoff concerns this chosen bigraded minimal
model, not all possible higher ordinary Massey products.

\section*{Acknowledgment}
The author used ChatGPT for mathematical exploration and proof development, literature searches, computational experiments, and code for numerical certificates. Proofs developed with this assistance were carefully checked, adapted, and rewritten by the author, who also rewrote the exposition. The author takes full responsibility for the mathematical content, computations, and any errors.

\end{document}